\documentclass[letter,12pt]{article}%
\usepackage{amsmath,bm}
\usepackage{amsfonts}
\usepackage{amssymb}
\usepackage{graphicx}
\usepackage{rotating}
\usepackage[height=9in,left=1in,right=0.75in,bottom=1in]{geometry}
\usepackage[breaklinks=true,bookmarksopen=true,colorlinks=true,citecolor=blue]%
{hyperref}
\usepackage[authoryear]{natbib}
\usepackage{color}
\usepackage{colortbl}
\usepackage{paralist}
\usepackage{amsmath}%
\providecommand{\U}[1]{\protect\rule{.1in}{.1in}}
\RequirePackage{hypernat}

\catcode`~=11
\newcommand{\urltilde}{\kern -.15em\lower .7ex\hbox{~}\kern .04em}
\catcode`~=13
\def \@seccntformat#1{\csname the#1\endcsname.\quad}
\makeatother
\numberwithin{equation}{section}
\hypersetup{pdftitle={Escanciano}, pdfsubject={Robust Two-Step Semiparametric Estimation}, pdfauthor={Juan Carlos Escanciano}, pdfkeywords={Two-step estimation; Stochastic Equicontinuity;
Semiparametric inference} }
\begin{document}

\title{Locally Robust Semiparametric Estimation}
\author{Victor Chernozhukov\\\textit{MIT}
\and Juan Carlos Escanciano\\\textit{Universidad Carlos III de Madrid}
\and Hidehiko Ichimura\\\textit{University of Arizona}
\and Whitney K. Newey\\\textit{MIT and NBER}
\and James M. Robins\\\textit{Harvard University}}
\date{July 2020}
\maketitle

\begin{abstract}
Many economic and causal parameters depend on nonparametric or high
dimensional first steps. We give a general construction of locally
robust/orthogonal moment functions for GMM, where moment conditions have zero
derivative with respect to first steps. We show that orthogonal moment
functions can be constructed by adding to identifying moments the
nonparametric influence function for the effect of the first step on
identifying moments. Orthogonal moments reduce model selection and
regularization bias, as is very important in many applications, especially for
machine learning first steps.

We give debiased machine learning estimators of functionals of high
dimensional conditional quantiles and of dynamic discrete choice parameters
with high dimensional state variables. We show that adding to identifying
moments the nonparametric influence function provides a general construction
of orthogonal moments, including regularity conditions, and show that the
nonparametric influence function is robust to additional unknown functions on
which it depends. We give a general approach to estimating the unknown
functions in the nonparametric influence function and use it to automatically
debias estimators of functionals of high dimensional conditional location
learners. We give a variety of new doubly robust moment equations and
characterize double robustness. We give general and simple regularity
conditions and apply these for asymptotic inference on functionals of high
dimensional regression quantiles and dynamic discrete choice parameters with
high dimensional state variables.

Keywords: Local robustness, orthogonal moments, double robustness,
semiparametric estimation, bias, GMM.

\begin{description}
\item[JEL classification:] C13; C14; C21; D24

\end{description}
\end{abstract}

\section{Introduction}

Many economic and causal parameters depend on nonparametric or high
dimensional first steps. Examples include dynamic discrete choice, games,
average consumer surplus, and treatment effects. This paper shows how to
construct moment functions for GMM estimators that are locally robust,
referred to henceforth as orthogonal, where moment conditions have a zero
derivative with respect to first steps. We show that such moment functions can
be constructed by adding to identifying moment functions the nonparametric
influence function from the effect of the first step on identifying moments.
This construction follows practice where identifying moment conditions are
often derived first from economic or causal models and then orthogonal moment
functions constructed. Also the nonparametric influence function is entirely
determined by the identifying moment functions and first step which motivates
automatic ways to estimate the nonparametric influence function. In
constructing sample moment functions we also cross-fit, a form of sample
splitting where the moment function for each observation is evaluated at first
step estimators that only use other observations, which further reduces bias.
A GMM estimator based on orthogonal moment functions with cross-fitting used
in construction of sample moments is referred to here as a debiased GMM estimator.

Debiased GMM has several advantages over plug-in GMM where only the
identifying moment functions are used. First, standard confidence intervals
for plug-in GMM are invalid under local alternatives when there is model
selection in the first step while confidence intervals with orthogonal moment
functions remain valid. Thus GMM with orthogonal moments is preferred over
plug-in GMM in the many applications with first step model selection. Second,
with regularized first steps debiased GMM will be root-n consistent under
conditions where plug-in GMM is not. Model selection and/or regularization is
often an important feature of machine learning, which is useful for estimating
econometric models with many regressors or state variables, making debiased
GMM especially important with machine learning first steps. Third, orthogonal
moment functions will be doubly robust when the orthogonal moment functions
are affine in the first step. We give this double robustness characterization
and use it to derive new classes of doubly robust estimators. Fourth, in some
important settings debiased GMM is known to have faster remainder rates or
smaller 2nd order mean square error than plug-in GMM. In addition regularity
conditions for debiased GMM are general and simple relative to those for
plug-in GMM. We show asymptotic normality for debiased GMM for any first step
estimator where certain mean square consistency conditions hold and either one
(under double robustness) or two (more generally) mean square rates hold. We
also show that conditions for plug-in GMM do not share this generality and
simplicity due to an additional remainder that is specific to the first step
estimator and complicated.

Debiased GMM is computationally more complicated than plug-in GMM in requiring
estimation of additional unknown functions on which the nonparametric
influence function depends. Also plug-in GMM has been innovative and useful in
a variety of settings, including Powell, Stock, and Stoker (1989), Hotz and
Miller (1993), Newey (1994a), Shen (1997), Chen and Shen (1998), Ai and Chen
(2003, 2007), and many others. The advantages of debiased GMM discussed here
motivate its use as an alternative to plug-in GMM, especially for machine
learning or any other first steps involving model selection and/or regularization.

Machine learning is useful for estimating economic and causal models where
there are high dimensional covariates or state variables, e.g. as in Belloni
et al. (2012), Belloni, Chernozhukov, and Wei (2013), Robins et al. (2013),
Belloni, Chernozhukov, and Hansen (2014), Farrell (2015), Kandasamy et al.
(2015), Belloni, Chernozhukov, and Kato (2015), Belloni et al. (2017), and
Athey, Imbens, and Wager (2018). Machine learning methods that are useful for
these purposes include Lasso, Dantzig, neural nets, boosting, and others.
Orthogonal moment functions reduce model selection and/or regularization
biases which are common for machine learning first steps. Cross-fitting for
debiased GMM avoids the need for Donsker conditions, which do not hold for
many machine learning first steps, and reduces own observation bias. The large
sample theory given here imposes only mean square convergence properties which
will hold for a variety of machine learning first steps. The advantages of
debiased GMM make it preferred to plug-in GMM for many machine learning first steps.

Previous to this paper debiased GMM with machine learning first step was based
on orthogonal moment functions constructed in various ways. Constructing
orthogonal moment functions by adding the nonparametric influence function to
identifying moments opens the way to debiased GMM for many objects of
interest. We illustrate by constructing debiased GMM estimators for
functionals of conditional quantiles and for parameters of dynamic discrete
choice models. The dynamic discrete choice estimator is based on machine
learners of conditional choice probabilities allowing for high dimensional
state variables. The estimator incorporates a novel Lasso estimator of
conditional value function differences where the Lasso left-hand side variable
is a function of a machine learner of the conditional choice probability. The
estimator and the results we give provide a prototype for using machine
learning for dynamic structural models. The conditional quantile estimator
allows for high dimensional regressors.

This paper gives automatic estimators of the additional functions on which the
nonparametric influence function depends. The approach uses just the
orthogonal moment functions and the first step to construct estimators of the
additional functions. The conditional quantile estimator employs this
automatic estimator. This approach generalizes the automatic method in
Chernozhukov, Newey, and Singh (2018).

We show that the nonparametric influence function has a useful robustness
property. The robustness property is that the expected value of the
nonparametric influence function is zero when the additional functions are not
equal to the truth but the first step is. Consequently, the estimator of these
additional functions is not required to converge faster than $n^{-1/4}$. We
also show orthogonality using a standard Gateaux derivative characterization
of the nonparametric influence function and some regularity conditions.

Orthogonal estimators for functionals of a density constructed by adding the
nonparametric influence function have previously been given by Hasminskii and
Ibragimov (1978), Pfanzagl and Wefelmeyer (1981), and Bickel and Ritov (1988).
Newey, Hsieh, and Robins (1998, 2004) suggested the construction of orthogonal
moment functions from adding the nonnparametric influence function. This
construction was considered in Chernozhukov et al. (2016), Chernozhukov et al.
(2018), and Bravo, Escanciano, and van Keilegom (2020). This paper innovates
by showing robustness of the nonparametric influence function to the
additional unknown functions on which it depends, so that these additional
functions need not be estimated at a $n^{-1/4}$ rate. Also, regularity
conditions are given for any first step to have no first order effect on
expected orthogonal moments, the precise orthogonality condition for used in
the asymptotic theory. These results are obtained using a standard Gateaux
derivative characterization of the influence function. None of the Theorems in
this paper appear in previous work. The relationship of the orthogonalization
results in this paper to previous literature is discussed more fully in
Section 4.

The Robinson (1988) semiparametric regression and Ichimura (1993) index
regression estimators have first order conditions that are orthogonal moment
functions. The object of interest in those papers minimizes an objective
function that is an expectation that is also minimized by the first step. The
objects of interest we consider are much more general in including many
economic and causal parameters that do not minimize the same expectation as
the first step.

Doubly robust moment functions have been constructed by Robins, Rotnitzky, and
Zhao (1994, 1995), Robins and Rotnitzky (1995), Scharfstein, Rotnitzky, and
Robins (1999), Robins, Rotnitzky, and van der Laan (2000), Robins and
Rotnitzky (2001), Graham (2011), and Firpo and Rothe (2019). This paper
innovates by deriving large classes of new doubly robust moment functions,
including affine functionals of nonparametric regressions, functions
satisfying other conditional moment restrictions, and density estimators, and
by characterizing double robustness, all based on adding the nonparametric
influence function. We also give related, partial robustness results where
original moment conditions are satisfied even when the first step is not equal
to the truth.

Targeted maximum likelihood, Van der Laan and Rubin (2006), based on machine
learners has been considered by Van der Laan and Rose (2011) and large sample
theory given by Luedtke and Van Der Laan (2016), Toth and Van der Laan (2016),
and Zheng et al. (2016). Here we directly target parameters of interest via
GMM based on adding the nonparametric influence function, with automatic
estimation of additional unknown parameters and general and simple regularity
conditions for asymptotic inference.

Recent work on debiased machine learning by Chernozhukov et al. (2018),
Chernozhukov, Newey, and Robins (2018), and Chernozhukov, Newey, and Singh
(2018) is partly based on and is also generalized by this paper. The
construction of orthogonal moments here was described in Chernozhukov et al.
(2018), which cited this paper for that construction and contains no results
from this paper. The asymptotic theory in this paper uses the orthogonal
moment construction here to improve on the asymptotic theory Chernozhukov et
al. (2018), as described in Section 8. The doubly robust moment conditions
considered in Chernozhukov, Newey, and Robins (2018) and Chernozhukov, Newey,
and Singh (2018) were derived in this paper and the asymptotic theory in those
other papers uses theory given in this paper. The automatic machine learner of
the additional unknown functions given here generalizes that in Chernozhukov,
Newey, and Singh (2018). In addition Newey and Robins (2017) and Hirshberg and
Wager (2018) are concerned with linear functions of a regression that are
formulated here. Furthermore, Bonhomme and Weidner (2018) have shown the
importance of orthogonal moment functions in specification analysis, Foster
and Srygkanis (2019) in deriving rates of convergence for machine learners,
Chernozhukov, Hausman, and Newey (2019) for demand analysis with endogenous
total expenditure, Semenova (2019) for machine learning for partially
identified models, Singh and Sun (2019) for machine learning of complier
effects, and Chernozhukov, Semenova, and Newey (2019) for machine learning of
weighted average value functions in dynamic structural models.

There are other sources of bias arising from nonlinearity of moment conditions
in the first step. Cattaneo and Jansson (2018) and Cattaneo, Jansson, and Ma
(2018) give useful bootstrap and jackknife methods that reduce nonlinearity
bias. Newey and Robins (2017) show that one can also remove this bias by cross
fitting in some settings. We use cross-fitting in this paper.

To summarize the contributions of this paper, we consider GMM estimation with
nonparametric first steps, with orthogonal moment conditions constructed by
adding the nonparametric influence function to identifying moment functions.
We give novel such estimators of functionals of high dimensional conditional
quantiles and of dynamic discrete choice parameters with high dimensional
state variables, including a novel Lasso estimator of conditional value
function differences. We show that that the nonparametric influence function
is robust to additional unknown functions on which it depends, so that
$n^{-1/4}$ consistency is not required there. We show orthogonality of the
constructed moment functions using the standard Gateaux derivative
characterization of the influence function, including regularity conditions.
We give examples showing that plug-in GMM is severely biased by model
selection and/or first step regularization whereas debiased GMM is not. We
give a general approach to estimating additional unknown functions in the
nonparametric influence function and use it to automatically debias estimators
of functionals of high dimensional conditional location learners, including
regression quantiles. We give a variety of new doubly robust moment equations
and characterize double robustness. We give general and simple regularity
conditions that improve on previous conditions and apply these for asymptotic
inference on functionals of high dimensional regression quantiles and dynamic
discrete choice parameters with high dimensional state variables.

Section 2 describes orthogonal moments and debiased GMM, and gives the
conditional quantile estimator. Section 3 gives the dynamic discrete choice
estimator and reports results of a Monte Carlo study. Section 4 shows
orthogonality and the robustness of the nonparametric influence function.
Section 5 compares the properties of debiased and plug-in GMM estimators.
Section 6 gives automatic estimators of the additional functions. Section 7
gives novel classes of doubly robust moment functions and characterizes double
robustness. Section 8 provides general and simple asymptotic theory for
debiased GMM.

\section{Debiased GMM}

The subject of this paper is GMM estimators of parameters identified by moment
functions that depend on a first step unknown function. In this Section we
describe this type of estimator and give examples.

\subsection{The Estimator}

To describe such an estimator let $\theta$ denote a finite dimensional
parameter vector of interest, $\gamma$ an unknown function, and $W$ a data
observation. We assume that there is a vector $g(w,\gamma,\theta)$ of known
functions of a possible realization $w$ of $W$ such that%
\[
E[g(W,\gamma_{0},\theta_{0})]=0,
\]
where $\theta_{0}$ and $\gamma_{0}$ are the true parameter vector and
function. We will assume that the parameter is identified by these moments,
i.e. that $\theta_{0}$ is the unique solution to $E[g(W,\gamma_{0},\theta)]=0$
over $\theta$ in some set $\Theta.$

The true function $\gamma_{0}$ is unknown so a first step estimator
$\hat{\gamma}$ of $\gamma_{0}$ is used$.$ Let $W_{1},...,W_{n}$ be a sample of
i.i.d. data observations. Estimated sample moment functions can be formed by
plugging in $\hat{\gamma}$ into $g(W_{i},\gamma,\theta)$ and averaging over
data observations to obtain $\sum_{i=1}^{n}g(W_{i},\hat{\gamma},\theta)/n.$
One could form a "plug-in" GMM estimator by minimizing a quadratic form in
these estimated sample moments, but such an estimator will be highly biased by
first step model selection and/or regularization as further detailed in
Section 5. This bias can be reduced by using orthogonal moment functions$.$

Orthogonal moment functions are based on influence functions. To describe them
we need to explain some additional concepts and notation. Let $F$ denote a
possible CDF for a data observation $W$ and suppose that $\hat{\gamma}$ has a
probability limit $\gamma(F)$ when $F$ is the true distribution of $W$. Here
$\gamma(F)$ is the probability limit of $\hat{\gamma}$ under general
misspecification, similar to Newey (1994a), so that $F$ is unrestricted except
for regularity conditions such as existence of $\gamma(F)$ or the expectation
of certain functions of the data. For example if $\hat{\gamma}(x)$ is a
nonparametric estimator of $E[Y|X=x]$ then $\gamma(F)(x)=E_{F}[Y|X=x]$ is the
conditional expectation function when $F$ is the true distribution of $W$,
which is well defined under the regularity condition that $E_{F}[\left\vert
Y\right\vert ]$ is finite.

Next, let $F_{0}$ denote the true distribution of $W,$ $H$ some alternative
distribution, and $F_{\tau}=(1-\tau)F_{0}+\tau H$ for $\tau\in\lbrack0,1].$ We
assume that $H$ is chosen so that $\gamma(F_{\tau})$ exists for $\tau$ small
enough and possibly other regularity conditions are satisfied. We also make
the key assumption that there exists a function $\phi(w,\gamma,\alpha,\theta)$
such that%

\begin{equation}
\frac{d}{d\tau}E[g(W,\gamma(F_{\tau}),\theta)]=\int\phi(w,\gamma_{0}%
,\alpha_{0},\theta)H(dw),\text{ \ \ }E[\phi(W,\gamma_{0},\alpha_{0}%
,\theta)]=0. \label{infdef}%
\end{equation}
Here $\alpha$ is an additional unknown function on which $\phi(w,\gamma
_{0},\alpha_{0},\theta)$ depends and $d/d\tau$ is the derivative from the
right (i.e. for nonnegative values of $\tau$) at $\tau=0.$ This equation is
the well known characterization of the influence function $\phi(w,\gamma
_{0},\alpha_{0},\theta)$ of $\mu(F)=E[g(W,\gamma(F),\theta)]$ as the Gateaux
derivative of $\mu(F),$ as in Von Mises (1947), Hampel (1974), Huber (1981).
The restriction that $\gamma(F_{\tau})$ exists allows $\phi(w,\gamma
_{0},\alpha_{0},\theta)$ to be the influence function when $\gamma(F)$ is only
well defined for certain types of distributions, such as when $\gamma(F)$ is a
conditional expectation or density. The function $\phi(w,\gamma,\alpha
,\theta)$\ will generally exist when $E[g(W,\gamma(F),\theta)]$ has a finite
semiparametric variance bound. Also $\phi(w,\gamma,\alpha,\theta)$ will
generally be unique because we are not restricting $H$ except for regularity
conditions. We will refer to $\phi(w,\gamma,\alpha,\theta)$ as the
\textit{nonparametric influence function} as it characterizes the local effect
of the first step function $\gamma$ on the expected moment $\mu(F).$ The
nonparametric influence function can be calculated from the derivative in
equation (\ref{lrdef}) or as described in Newey (1994a); see Ichimura and
Newey (2017).

Orthogonal moment functions are constructed by adding the nonparametric
influence function to the identifying moment functions to obtain%
\[
\psi(W,\gamma,\alpha,\theta)=g(W,\gamma,\theta)+\phi(W,\gamma,\alpha,\theta).
\]
Estimation of the unknown functions $\gamma$ and $\alpha$ will have no first
order effect on the expected value of $\psi(W,\gamma,\alpha,\theta)$ and
estimation of $\theta$ will not affect the expectation of $\phi(W,\gamma
,\alpha,\theta)$ as we will show. Debiased sample moments can then be
constructed by evaluating at first step estimators of $\gamma,$ $\alpha,$ and
$\theta$ and averaging over the data.

Constructing orthogonal moment functions is greatly facilitated by the wide
variety of known $\phi(w,\gamma,\alpha,\theta).$ The form of $\phi
(w,\gamma,\alpha,\theta)$ for density weighted average derivatives is in
Powell, Stock, and Stoker (1989). For first step least squares projections
(including conditional expectations), density weighted conditional means, and
their derivatives $\phi(W,\gamma,\alpha,\theta)$ is given in Newey (1994a).
Hahn (1998) and Hirano, Imbens, and Ridder (2003) used those results to obtain
$\phi(W,\gamma,\alpha,\theta)$ for treatment effect estimators. Bajari, Hong,
Krainer, and Nekipelov (2010) and Bajari, Chernozhukov, Hong, and Nekipelov
(2009) derived $\phi(W,\gamma,\alpha,\theta)$ for some first steps used in
structural estimation. Hahn and Ridder (2013, 2019) derived $\phi
(W,\gamma,\alpha,\theta)$ for generated regressors that depend on first step
conditional expectations. Ai and Chen (2007) and Ichimura and Newey (2017)
derived $\phi(W,\gamma,\alpha,\theta)$ for first step estimators of functions
satisfying conditional moment restrictions. Semenova (2018) derived
$\phi(W,\gamma,\alpha,\theta)$ for support functions used in partial
identification. Many of these derivations did not directly use equation
(\ref{infdef}) but Ichimura and Newey (2017) show that various $\phi
(W,\gamma,\alpha,\theta)$ solve equation (\ref{infdef}). These wide variety of
known $\phi(W,\gamma,\alpha,\theta)$ can be used to construct orthogonal
moment functions.

We also use cross-fitting, a form of sample splitting, in the construction of
debiased sample moments; e.g. see Schick (1986) and Klaassen (1987). Partition
the observation indices $(i=1,...,n)$ into $L$ groups $I_{\ell},$
$(\ell=1,...,L).$ Let $\hat{\gamma}_{\ell},$ $\hat{\alpha}_{\ell}$, and
$\tilde{\theta}_{\ell}$ be estimators that are constructed using all
observations \textit{not} in $I_{\ell}.$ Debiased sample moments functions are%
\begin{equation}
\hat{\psi}(\theta)=\hat{g}(\theta)+\hat{\phi},\text{ }\hat{g}(\theta)=\frac
{1}{n}\sum_{\ell=1}^{L}\sum_{i\in I_{\ell}}g(W_{i},\hat{\gamma}_{\ell}%
,\theta),\text{ }\hat{\phi}=\frac{1}{n}\sum_{\ell=1}^{L}\sum_{i\in I_{\ell}%
}\phi(W_{i},\hat{\gamma}_{\ell},\hat{\alpha}_{\ell},\tilde{\theta}_{\ell}).
\label{dmom}%
\end{equation}
A debiased GMM estimator is%
\begin{equation}
\hat{\theta}=\arg\min_{\theta\in\Theta}\hat{\psi}(\theta)^{\prime}%
\hat{\Upsilon}\hat{\psi}(\theta), \label{dgmm}%
\end{equation}
where $\hat{\Upsilon}$ is a positive semi-definite weighting matrix and
$\Theta$ is the set of parameter values. A choice of $\hat{\Upsilon}$ that
minimizes the asymptotic variance of $\hat{\theta}$ will be $\hat{\Upsilon
}=\hat{\Psi}^{-1},$ for
\[
\hat{\Psi}=\frac{1}{n}\sum_{\ell=1}^{L}\sum_{i\in I_{\ell}}\hat{\psi}_{i\ell
}\hat{\psi}_{i\ell}^{\prime},\text{ }\hat{\psi}_{i\ell}=g(W_{i},\hat{\gamma
}_{\ell},\tilde{\theta}_{\ell})+\phi(W_{i},\hat{\gamma}_{\ell},\hat{\alpha
}_{\ell},\tilde{\theta}_{\ell}).
\]
There is no need to account for the presence of $\hat{\gamma}$ and
$\hat{\alpha}$ in $\hat{\psi}_{i\ell}$ because the presence of $\phi
(W_{i},\hat{\gamma}_{\ell},\hat{\alpha}_{\ell},\tilde{\theta}_{\ell})$ removes
the first order effect on $\hat{\psi}(\theta)$ of $\hat{\gamma}_{\ell}$ and
$\hat{\alpha}_{\ell}$ under conditions we give. An estimator $\hat{V}$ of the
asymptotic variance of $\sqrt{n}(\hat{\theta}-\theta_{0})$ is%
\begin{equation}
\hat{V}=(\hat{G}^{\prime}\hat{\Upsilon}\hat{G})^{-1}\hat{G}^{\prime}%
\hat{\Upsilon}\hat{\Psi}\hat{\Upsilon}\hat{G}(\hat{G}^{\prime}\hat{\Upsilon
}\hat{G})^{-1},\text{ }\hat{G}=\frac{\partial\hat{\psi}(\hat{\theta}%
)}{\partial\theta}._{{}} \label{varest}%
\end{equation}

The cross-fitting used in this estimator, where $\hat{\psi}(\theta)$ is
averaged over observations not used to form $\hat{\gamma}_{\ell}$ and
$\hat{\alpha}_{\ell},$ eliminates bias due to averaging over observations that
are used to construct the first step$.$ Eliminating such "own observation"
bias helps remainders converge faster to zero, e.g. as in Newey and Robins
(2017), and can be important in practice, e.g. as in Angrist and Krueger
(1995). It also eliminates the need for Donsker conditions for $\hat{\gamma
}_{\ell}$ and $\hat{\alpha}_{\ell},$ which is important for many machine
learner first steps that are not known to satisfy such conditions, as
discussed in Chernozhukov et al. (2018).

The debiased moments require a first step estimator $\hat{\alpha}_{\ell}$ of
unknown functions that appear in $\phi(W,\gamma,\alpha,\theta).$ When the form
of $\alpha_{0}$ is known one can plug-in an estimator $\hat{\alpha}_{\ell}.$
Also, in Section 6 we use the debiased moment functions to construct an
automatic estimator $\hat{\alpha}_{\ell}$ that does not requiring knowing the
form of $\alpha_{0}$. This automatic method generalizes that of Chernozhukov,
Newey, and Singh (2018) beyond functionals of least squares projections.

The efficiency of debiased GMM is entirely determined by the choice of moment
functions, first step, and weighting matrix. The matrix $\hat{\Psi}^{-1}$ is
an optimal choice of weighting matrix as usual for GMM. The efficient choice
of moment functions and first steps will depend on a model and as further
discussed in Section 4. The presence of $\hat{\phi}$ in the orthogonal moment
functions $\hat{\psi}(\theta)$ does not affect identification of $\theta$. The
second, mean zero condition in equation (\ref{infdef}) holds for all possible
distributions of $W$ so that $\hat{\phi}$ will converge in probability to
zero. The sole purpose of including $\hat{\phi}$ is to remove the first-order
effect of $\hat{\gamma}_{\ell}$ on the moment functions. It accomplishes this
in a nonparametric way that does not depend on any model assumptions, as
further discussed in Section 4.

The initial estimator $\tilde{\theta}_{\ell}$ can be based on the identifying
moment conditions and constructed as%
\[
\tilde{\theta}_{\ell}=\arg\min_{\theta\in\Theta}\hat{g}_{\ell}(\theta
)^{\prime}\hat{\Upsilon}_{\ell}\hat{g}_{\ell}(\theta),\text{ }\hat{g}_{\ell
}(\theta)=\frac{1}{n-n_{\ell}}\sum_{\ell^{\prime}\neq\ell}\sum_{i\in
I_{\ell^{\prime}}}g(W_{i},\tilde{\gamma}_{\ell^{\prime}},\theta).
\]
where $\hat{\Upsilon}_{\ell}$ uses only observations not in $I_{\ell},$
$n_{\ell}$ is the number of observations in $I_{\ell}$, and $\tilde{\gamma
}_{\ell^{\prime}}$ uses observations not in $\ell$ and not in $\ell^{\prime}.$
One could iterate on the initial estimator $\tilde{\theta}_{\ell}$ in the
debiased moments $\hat{\psi}(\theta)$ and/or $\hat{\Psi}$ by calculating
$\hat{\theta}$ and/or $\hat{\Psi}$ a second time with $\tilde{\theta}_{\ell}$
being a debiased GMM estimator obtained from a prior iteration.

\bigskip

\textsc{Example 1:} An example that will be used to illustrate the theory has
a data observations $W=(Y,X,Z)$ and $\theta_{0}=E[Z\gamma_{0}(X)]=E[\alpha
_{0}(X)\gamma_{0}(X)]$ for $\gamma_{0}(X)=E[Y|X]$ and $\alpha_{0}(X)=E[Z|X].$
This example is of interest in its own right as the component of the expected
conditional covariance $E[Cov(Z,Y|X)]=E[ZY]-\theta_{0}$ that depends on
unknown functions, which covariance is useful for the analysis of covariance
and for estimation of a partially linear model, Robinson (1988). We specify
the identifying moment function with implied nonparametric influence function
as%
\[
g(w,\gamma,\theta)=z\gamma(x)-\theta,\text{ }\phi(w,\gamma,\alpha
)=\alpha(x)[y-\gamma(x)].
\]
The debiased GMM estimator is%
\[
\hat{\theta}=\frac{1}{n}\sum_{\ell=1}^{L}\sum_{i\in I_{\ell}}\left\{
Z_{i}\hat{\gamma}_{\ell}(X_{i})+\hat{\alpha}_{\ell}(X_{i})[Y_{i}-\hat{\gamma
}_{\ell}(X_{i})]\right\}  =\frac{1}{n}\sum_{\ell=1}^{L}\sum_{i\in I_{\ell}%
}\left\{  Z_{i}\hat{\gamma}_{\ell}(X_{i})+\hat{\alpha}_{\ell}(X_{i})Y_{i}%
-\hat{\alpha}_{\ell}(X_{i})\hat{\gamma}_{\ell}(X_{i})\right\}  .
\]

\subsection{Example 2: Functions of Conditional Quantiles}

The object of interest in this example is an expected linear functional of a
conditional quantile function where%
\[
\theta_{0}=E[m(W,\gamma_{0})],\text{ }\gamma_{0}=\arg\min_{\gamma\in\Gamma
}E[v(Y-\gamma(X))],\text{ }v(u)=[\lambda-1(u<0)]u.
\]
Here $Y$ is a dependent variable of interest and $m(w,\gamma)$ is a linear
functional of $\gamma.$ An example is a weighted average derivative of
$\gamma_{0}$ where $m(w,\gamma)=\int d(x)[\partial\gamma(x)/\partial
x_{1}]dx.$ Here the identifying moment function is $g(w,\gamma,\theta
)=m(w,\gamma)-\theta$. The nonparametric influence function is $\phi
(w,\gamma,\alpha)=\alpha(x)v_{u}(y-\gamma(x)),$ where $v_{u}(u)=\lambda
-1(u<0)$ denotes the derivative of $v(u)$ away from zero, as follows from Ai
and Chen (2007, p. 40) and Ichimura and Newey (2017).

We can construct a debiased GMM estimator of $\theta_{0}$ from any learner of
the conditional quantile $\gamma_{0}$ that converges sufficiently fast in mean
square. Let $\hat{\gamma}_{\ell}(x)$ be a learner of $\gamma_{0}$, computed
from observations not in $I_{\ell}$, and $\hat{\gamma}_{\ell,\ell^{\prime}%
}(x)$ be computed from observations not in $\ell$ or $\ell^{\prime}.$ Also let
$K(u)$ be a bounded, univariate kernel, with $\int K(u)du=1$ and $\int
K(u)udu=0,$ and $h$ a bandwidth. Let $b(x)$ be a $p\times1$ vector of
functions of $x.$ A debiased GMM estimator is
\begin{align}
\hat{\theta}  &  =\frac{1}{n}\sum_{\ell=1}^{L}\sum_{i\in I_{\ell}}%
\{m(W_{i},\hat{\gamma}_{\ell})+\hat{\alpha}_{\ell}(X_{i})v_{u}(Y_{i}%
-\hat{\gamma}_{\ell}(X_{i}))\},\label{quan debia}\\
\hat{\alpha}_{\ell}(x)  &  =b(x)^{\prime}\hat{\rho}_{\ell},\text{ }\hat{\rho
}_{\ell}=\arg\min_{\rho}\{-2\hat{M}_{\ell}^{\prime}\rho+\rho^{\prime}\hat
{Q}_{\ell}\rho+2r_{\lambda}\sum_{j=1}^{p}\left\vert \rho_{j}\right\vert
\},\text{ }\hat{M}_{\ell}=(\hat{M}_{\ell1},...,\hat{M}_{\ell p})^{\prime
},\nonumber\\
\hat{M}_{\ell j}  &  =\frac{1}{n-n_{\ell}}\sum_{i\notin I_{\ell}}m(W_{i}%
,b_{j}),\text{ }\hat{Q}_{\ell}=\frac{1}{n-n_{\ell}}\sum_{\ell^{\prime}\neq
\ell}\sum_{i\in I_{\ell^{\prime}}}\frac{1}{h}K(\frac{Y_{i}-\hat{\gamma}%
_{\ell,\ell^{\prime}}(X_{i})}{h})b(X_{i})b(X_{i})^{\prime}.\nonumber
\end{align}
The $\hat{\alpha}_{\ell}(x)$ is a special case of an automatic Lasso minimum
distance learner given and motivated in Section 6.

This estimator depends on the regularization term $r_{\lambda}$ in the
objective function for $\hat{\rho}_{\ell}$. This $r_{\lambda}$ should be
chosen larger than the conventional Lasso regularization parameter because
$\hat{Q}_{\ell}$ depends on the nonparametric estimator $\hat{\gamma}%
_{\ell,\ell^{\prime}}(X_{i})$, as further discussed in Section 8. The nested
sample splitting used for $\hat{\gamma}_{\ell,\ell^{\prime}}(X_{i})$ requires
that the first step learner be computed for $L^{2}$ subsamples. Use of
$\hat{\gamma}_{\ell}(x)$ as a starting value for computation of each
$\hat{\gamma}_{\ell,\ell^{\prime}}(X_{i})$ may aid in this computation. The
nested sample splitting allows for a very general first step that need only
have a mean square convergence rate. Here $\hat{\gamma}_{\ell}(X_{i})$ could
be used in place of $\hat{\gamma}_{\ell,\ell^{\prime}}(X_{i})$ if$\ \sum
_{i\notin I_{\ell}}\left\vert \hat{\gamma}_{\ell}(X_{i})-\gamma_{0}\left(
X_{i}\right)  \right\vert /(n-n_{\ell})$ converged to zero at some rate that
is a power of $n.$

\section{Example 3: Dynamic Discrete Choice}

For dynamic discrete choice with high dimensional state variables we estimate
structural parameters via learners of conditional choice probabilities. This
approach replaces computation of expected value functions with nonparametric
estimation as suggested in Hotz and Miller (1993). For simplicity we focus
this example on binary choice, providing methods and results that will be
available for the more complicated models employed more widely in practice. In
particular we provide a Lasso estimator of conditional value function
differences where the dependent variable is a function of estimated future
choice probabilities. In Section 8 we provide convergence rate results for
this estimator. We also give the nonparametric influence function for this
kind of first step. This analysis provides a prototype for estimation of
dynamic structural models with high dimensional state variables.

In dynamic binary choice individuals choose between two alternatives $j=1$ and
$j=2$ to maximize the expected present discounted value of per period utility
$U_{tj}=D_{j}(X_{t})^{\prime}\theta_{0}+\varepsilon_{tj},$
$(j=1,2;t=1,...,T),$ where $\varepsilon_{jt}$ is i.i.d. with known CDF,
independent of the entire history $\{X_{s}\}_{s=1}^{\infty}$ of a state
variable vector $X$, and $X_{t}$ is Markov of order $1$ and stationary. The
parameter of interest is $\theta_{0}$. We develop an estimator that allows for
high dimensional $X_{t}$.

We assume that choice $1$ is a renewal choice where the conditional
distribution of $X_{t+1}$ given $X_{t}$ and choice $1$ does not depend on
$X_{it}$. We also assume that $D_{1}(X_{t})=(-1,0^{\prime})^{\prime}$ and
$D_{21}(X_{t})=0$ so that the first element in $\theta$ is a binary choice
constant. Let $Y_{jt}$ equal a dummy variable equal to $1$ when choice $j$ is
made and $\gamma_{10}(X_{t})=\Pr(Y_{2t}=1|X_{t})$ be the conditional choice
probability of alternative 2. Also let $V(X_{t})$ denote the expected value
function. As in Hotz and Miller (1993) there is a known function $H(p)$ such
that for $\gamma_{20}(X_{t})=E[H(\gamma_{10}(X_{t+1}))|X_{t},Y_{2t}=1]$ and
$\gamma_{30}=E[H(\gamma_{10}(X_{t+1}))|Y_{1t}=1],$
\begin{equation}
E[V(X_{t+1})|X_{t},Y_{2t}=1]-E[V(X_{t+1})|Y_{1t}=1]=\gamma_{20}(X_{t}%
)-\gamma_{30}. \label{HotMill}%
\end{equation}
For example when $\varepsilon_{1t}$ and $\varepsilon_{2t}$ are independent
Type I\ extreme value this equation is satisfied for $H(p)=.5227-\ln(1-p).$
Then for the CDF $\Lambda(a)$ of $\varepsilon_{t1}-\varepsilon_{t2},$
$D(X_{t})=D_{2}(X_{t})-D_{1}(X_{t})$, and $\delta$ the discount factor the
conditional choice probability for $j=2$ is%
\begin{equation}
\Pr(Y_{2t}=1|X_{t})=\Lambda(a(X_{t},\theta_{0},\gamma_{20},\gamma
_{30})),\text{ }a(x,\theta,\gamma_{2},\gamma_{3})=D(x)^{\prime}\theta
+\delta\{\gamma_{2}(x)-\gamma_{3}\}. \label{prob}%
\end{equation}

We consider data of i.i.d. observations on individuals each followed for $T$
time periods, that also includes the $T+1$ observation $X_{T+1}$ of the state
variables, where $W=(X_{1}^{\prime},Y_{21},...,X_{T}^{\prime},Y_{2T}%
,X_{T+1}^{\prime})^{\prime}.$ An estimator of $\theta_{0}$ can be obtained by
constructing first step estimators $\hat{\gamma}_{2}(x)$ and $\hat{\gamma}%
_{3}$, substituting these estimators for $\gamma_{2}$ and $\gamma_{3}$ in
$a(x,\theta,\gamma_{2},\gamma_{3})$ in equation (\ref{prob}), and then
maximizing a binary choice log-likelihood as if $\hat{\gamma}_{2}(x)$ and
$\hat{\gamma}_{30}$ were true. We specify as identifying moment functions the
derivative of the pseudo log-likelihood associated with the binary choice
probability in equation (\ref{prob}) with respect to $\theta$,%
\begin{align}
g(W,\gamma,\theta)  &  =\frac{1}{T}\sum_{t=1}^{T}D(X_{t})\pi(a(X_{t}%
,\theta,\gamma_{2},\gamma_{3}))[Y_{2t}-\Lambda(a(X_{t},\theta,\gamma
_{2},\gamma_{3}))],\label{disc choice mom}\\
\pi(a)  &  =\frac{\Lambda_{a}(a)}{\Lambda(a)[1-\Lambda(a)]},\text{ }%
\Lambda_{a}(a)=\frac{d\Lambda(a)}{da}.\nonumber
\end{align}

Estimators of $\gamma_{10}(x)$, $\gamma_{20}(x)$, and $\gamma_{30}$ are needed
as a first step $\hat{\gamma}_{\ell}$ for the identifying moment function. We
will consider any $\hat{\gamma}_{1\ell}(x)$ that converges sufficiently
quickly in mean square. For example $\hat{\gamma}_{1\ell}$ could be logit
Lasso or a linear Lasso estimator with dependent variable $Y_{2t}$. We use
Lasso to construct $\hat{\gamma}_{2\ell}(x)$ in order to control for
estimation error that results from an estimated dependent variable. Let
$\hat{\gamma}_{1\ell,\ell^{\prime}}(x)$ be an estimator of the conditional
choice probability computed from observations not in $I_{\ell}$ or
$I_{\ell^{\prime}}.$ Let $b(x)$ denote a $p\times1$ dictionary of functions of
the state variables $x$. We form $\hat{\gamma}_{2\ell}(x)$ and $\hat{\gamma
}_{3\ell}$
\begin{align}
\hat{\gamma}_{2\ell}(x)  &  =b(x)^{\prime}\hat{\beta}_{2\ell},\text{ }%
\hat{\beta}_{2\ell}=\arg\min_{\beta}\{-2\hat{M}_{2\ell}^{\prime}\beta
+\beta^{\prime}\hat{G}_{2\ell}\beta+2r_{1}\sum_{j=1}^{p}\left\vert \beta
_{j}\right\vert \}\text{,}\label{gam2}\\
\hat{M}_{2\ell}  &  =\frac{1}{(n-n_{\ell})T}\sum_{\ell^{\prime}\neq\ell}%
\sum_{i\in I_{\ell^{\prime}}}\sum_{t=1}^{T}Y_{2it}b(X_{it})H\left(
\hat{\gamma}_{1\ell,\ell^{\prime}}(X_{i,t+1}\right)  ),\text{ }\hat{Q}_{2\ell
}=\frac{1}{(n-n_{\ell})T}\sum_{i\notin I_{\ell}}\sum_{t=1}^{T}Y_{2it}%
b(X_{it})b(X_{it})^{\prime},\nonumber\\
\hat{\gamma}_{3\ell}  &  =\frac{1}{\hat{P}_{1}(n-n_{\ell})T}\sum_{\ell
^{\prime}\neq\ell}\sum_{i\in I_{\ell^{\prime}}}\sum_{t=1}^{T}Y_{1it}H\left(
\hat{\gamma}_{1\ell,\ell^{\prime}}(X_{i,t+1}\right)  ),\text{ }\hat{P}%
_{1}=\frac{1}{(n-n_{\ell})T}\sum_{i\notin I_{\ell}}\sum_{t=1}^{T}%
Y_{1it}.\nonumber
\end{align}
Here $\hat{\gamma}_{2\ell}(x)$ is Lasso with left hand side variable $H\left(
\hat{\gamma}_{1\ell,\ell^{\prime}}(X_{i,t+1}\right)  )$ and right-hand side
variables $b(X_{it})Y_{2it}$ and $\hat{\gamma}_{3\ell}$ is a sample mean
conditional on $Y_{1it}=1.$ Here $\hat{\gamma}_{2\ell}$, and $\hat{\gamma
}_{3\ell}$ use all observations all observations with $i\notin I_{\ell}.$ The
nested sample splitting in $\hat{\gamma}_{1\ell,\ell^{\prime}}(x)$ is useful
in only requiring mean square convergence rates for conditional choice
probabilities although it is somewhat complicated. When $L$ is moderate size
(e.g. $L=5)$ there may be many $\hat{\gamma}_{1\ell,\ell^{\prime}}(x)$ to
compute (e.g. 25), although an estimate for a particular $\ell$ and
$\ell^{\prime}$ could provide a good starting value for other splits. As in
Example 2 the nested cross-fit estimator $\hat{\gamma}_{1\ell,\ell^{\prime}%
}(x)$can be replaced by single sample splitting if $\hat{\gamma}_{\ell}(x)$
converged in a sample absolute value terms.

The three first steps $\hat{\gamma}_{1\ell}(x),$ $\hat{\gamma}_{2\ell}(x)$,
and $\hat{\gamma}_{3\ell}$ result in a nonparametric influence function that
is the sum of three terms, one term for each first step, with%
\begin{align*}
\phi(W,\gamma,\alpha,\theta)  &  =\phi_{1}(W,\gamma,\alpha,\theta)+\phi
_{2}(W,\gamma,\alpha,\theta)+\phi_{3}(W,\gamma,\alpha),\\
\phi_{1}(W,\gamma,\alpha,\theta)  &  =\frac{1}{T}\sum_{t=1}^{T}\alpha
_{1}(X_{t},\theta)[Y_{2t}-\gamma_{1}(X_{t})],\text{ }\phi_{3}(W,\gamma
,\alpha)=\alpha_{3}\frac{1}{T}\sum_{t=1}^{T}Y_{2t}\{H(\gamma_{1}%
(X_{t+1}))-\gamma_{3}\},\\
\phi_{2}(W,\gamma,\alpha,\theta)  &  =\frac{1}{T}\sum_{t=1}^{T}\alpha
_{2}(X_{t},Y_{2t},\theta)[H(\gamma_{1}(X_{t+1}))-\gamma_{2}(X_{t})].
\end{align*}
The form of each $\phi_{j}(w,\gamma,\alpha,\theta)$ follows from Proposition 4
of Newey (1994a) by treating only $\gamma_{j}$ as unknown and holding the
other first steps fixed at their true values, as in Newey (1994a, p. 1357).
The true values $\alpha_{10},$ $\alpha_{20}$, and $\alpha_{30}$ are%
\begin{align}
\alpha_{10}(x,\theta_{0})  &  =\{E[\alpha_{20}(X_{t},Y_{2t},\theta
_{0})|X_{t+1}=x]+\alpha_{30}E[Y_{1t}|X_{t+1}=x]\}H_{p}(\gamma_{10}%
(x)),\label{alpha ddc}\\
\alpha_{20}(x,y_{2},\theta_{0})  &  =-\delta D(x)\pi(a(x))\frac{\Lambda
_{a}(a(x))y_{2}}{\Lambda(a(x))},\text{ }a(x)=a(x,\theta_{0},\gamma_{20}%
,\gamma_{30}),\nonumber\\
\alpha_{30}  &  =-E[\alpha_{20}(X_{t},Y_{2t},\theta_{0})]/P_{10},\text{
}P_{10}=E[Y_{1t}].\nonumber
\end{align}
\qquad\qquad

To construct $\hat{\alpha}_{2\ell},$ $\hat{\alpha}_{3\ell}$, and $\hat{\alpha
}_{1\ell}$, obtain an initial estimator $\tilde{\theta}_{\ell}$ from binary
choice pseudo maximum likelihood over $i\notin I_{\ell}$, $t\leq T$ with
$\gamma_{2}$ and $\gamma_{3}$ replaced by $\hat{\gamma}_{2\ell}$ and
$\hat{\gamma}_{3\ell}$ respectively in the choice probability formula. Also
let $\hat{a}_{it}=a(X_{it},\hat{\theta}_{\ell},\hat{\gamma}_{2\ell}%
,\hat{\gamma}_{3\ell})$ and construct $\hat{\alpha}_{2\ell}(X_{it}%
,Y_{2it},\tilde{\theta}_{\ell})$ by substituting $\hat{a}_{it}$ for $a(x)$,
$X_{it}$ for $x$, $Y_{2it}$ for $y_{2},$ and $\tilde{\theta}_{\ell}$ for
$\theta_{0}$in the formula for $\alpha_{20}(x,y_{2},\theta_{0})$ in equation
(\ref{alpha ddc}). Next obtain $\hat{\alpha}_{3\ell}$ by replacing
$\alpha_{20}(X_{t},Y_{2t},\theta_{0})$ by $\hat{\alpha}_{2}(X_{it}%
,Y_{2it},\tilde{\theta}_{\ell})$ and population expectations by sample
averages over $i\notin I_{\ell},t\leq T$ in the third line of equation
(\ref{alpha ddc}). Also, obtain $\hat{\alpha}_{1\ell}(x,\theta)$ by replacing
$\alpha_{30}$ and $\gamma_{10}(x)$ by $\hat{\alpha}_{3\ell}$ and $\hat{\gamma
}_{1\ell}(x)$ respectively in the first line of equation (\ref{alpha ddc}) and
by replacing the two conditional expectations by the predicted values from
Lasso regressions over $i\notin I_{\ell},t\leq T$ with regressors
$b(X_{i,t+1}),$ dependent variables equal to each element of $\hat{\alpha
}_{2\ell}(X_{it},Y_{2it},\tilde{\theta}_{\ell})$ and $Y_{1it}$ respectively,
and regularization factors $r_{2}$ and $r_{3}$ respectively, analogously to
$\hat{\gamma}_{2\ell}(x)$. Finally substitute $\hat{\alpha}_{1\ell},$
$\hat{\alpha}_{2\ell},$ and $\hat{\alpha}_{3\ell}$ for $\alpha_{10}%
,\ \alpha_{20}$, and $\alpha_{30}$ and $\hat{\gamma}_{1\ell}$, $\hat{\gamma
}_{2\ell}$ and $\hat{\gamma}_{3\ell}$ for $\gamma_{1},$ $\gamma_{2}$, and
$\gamma_{3}$ in the formulas for $\phi_{1}$, $\phi_{2}$, and $\phi_{3}$ and
construct a debiased GMM estimator as in Section 2.

To explore the finite sample properties of this estimator we carried out a
Monte Carlo study for a model similar to that of Rust (1987). The state
variables consisted of a positive variable $x_{1}$ (mileage) and other
variables $x_{2},...,x_{6}$ with transition%
\begin{align*}
X_{1,t+1}  &  =1+1(Y_{2t}=1)X_{1t}+S_{t+1}^{2},\text{ }S_{t+1}|X_{t}\sim
N(.2+\sum_{k=1}^{5}c_{k}X_{t,k+1},1),\\
c  &  =(.1,.025,.0111,.0063,.004);\text{ }%
\end{align*}
where $(X_{2t},...,X_{6t})$ is i.i.d. over $t$, $X_{2t},$ $X_{4t}$, and
$X_{6t}$ are chi-squared with one degree of freedom and $X_{3t}$ and $X_{5t}$
are binary with $\Pr(X_{kt}=1)=1/2$, $k=3,5.$ We specified that $D(x)$ is two
dimensional with $D_{1}(x)=(-1,0)^{\prime}$ and $D_{2}(x)=(0,\sqrt{x_{1}%
})^{\prime}$ and that $\varepsilon_{1t},$ $\varepsilon_{2t}$ are independent
Type I\ extreme value, so that $\Lambda(a)=e^{a}/(1+e^{a})$ corresponds to
binary logit.

To generate the data we solved the Bellman equation on a finite grid using the
fact that the state space has a two dimensional structure in terms of $x_{1}$
and $\sum_{k=1}^{5}c_{k}x_{k+1},$ with linear interpolation between grid
points. We did not enforce this index structure in estimation, so that the
estimation treated the state space as dimension six. We carry out 500 Monte
Carlo replications for $T=10$ and $n=100,$ $300,$ $1000$, and $10,000.$ We
specified five fold cross fitting where $L=5.$ We consider three
specifications of the vector $b(x)$ used by Lasso, consisting of a) the
elements of $x$, b) those from a) and squares of elements of $x$; c) those
from b) and all products of two elements of $x$. The conditional choice
probability estimators $\hat{\gamma}_{1\ell,\ell^{\prime}}(x)$ and
$\hat{\gamma}_{1\ell}(x)$ were logit Lasso truncated to be between $.0001$ and
$.9999.$ We used the MATLAB Lasso and logit Lasso procedures for computation.
The regularization value $r$ for Lasso was chosen by two fold regularization.
Although we do not know whether the resulting $r$ satisfies the conditions in
the asymptotic theory of Section 8, we uses this $r$ so that the estimator in
the Monte Carlo would be based on an "off the shelf" machine learner of
unknown functions.

The results are reported in Tables 1, 2, and 3. The PI labels the GMM
estimator based only on identifying moment functions, DB the debiased GMM,
Bias is the absolute value of bias, Med SE denotes the median of the estimated
standard errors corresponding to equation (\ref{varest}), SD denotes standard
deviation, and Cvg denotes coverage probability of a nominal 95 percent
confidence interval.

\begin{center}%
\begin{tabular}
[c]{cc|c|c|c|c|c|c|c|}%
\multicolumn{9}{c}{Table 1: $b(X)$ Linear}\\\cline{3-9}
&  & PI Bias & DB Bias & Med SE & PI SD & DB SD & PI Cvg & DB Cvg\\\hline
\multicolumn{1}{|c}{$n=100$} & \multicolumn{1}{|c|}{$\theta_{2}$} & .35 &
.00 & .13 & .11 & .10 & .17 & .99\\\cline{2-9}%
\multicolumn{1}{|c}{} & \multicolumn{1}{|c|}{$\theta_{1}$} & .61 & .04 & .21 &
.21 & .18 & .17 & .96\\\hline
\multicolumn{1}{|c}{$n=300$} & \multicolumn{1}{|c|}{$\theta_{2}$} & .33 &
.01 & .08 & .07 & .06 & .00 & .98\\\cline{2-9}%
\multicolumn{1}{|c}{} & \multicolumn{1}{|c|}{$\theta_{1}$} & .58 & .05 & .12 &
.13 & .11 & .00 & .94\\\hline
\multicolumn{1}{|c}{$n=1000$} & \multicolumn{1}{|c|}{$\theta_{2}$} & .33 &
.01 & .04 & .04 & .03 & .00 & .98\\\cline{2-9}%
\multicolumn{1}{|c}{} & \multicolumn{1}{|c|}{$\theta_{1}$} & .58 & .06 & .07 &
.07 & .06 & .00 & .87\\\hline
\multicolumn{1}{|c}{$n=10000$} & \multicolumn{1}{|c|}{$\theta_{2}$} & .32 &
.01 & .01 & .01 & .01 & .00 & .93\\\cline{2-9}%
\multicolumn{1}{|c}{} & \multicolumn{1}{|c|}{$\theta_{1}$} & .57 & .05 & .02 &
.02 & .02 & .00 & .21\\\hline
\end{tabular}

\begin{tabular}
[c]{cc|c|c|c|c|c|c|c|}%
\multicolumn{9}{c}{Table 2: $b(X)$ Linear, Squares}\\\cline{3-9}
&  & PI Bias & DB Bias & Med SE & PI SD & DB SD & PI Cvg & DB Cvg\\\hline
\multicolumn{1}{|c}{$n=100$} & \multicolumn{1}{|c|}{$\theta_{2}$} & .24 &
.03 & .13 & .11 & .33 & .61 & .98\\\cline{2-9}%
\multicolumn{1}{|c}{} & \multicolumn{1}{|c|}{$\theta_{1}$} & .41 & .03 & .21 &
.20 & .38 & .57 & .97\\\hline
\multicolumn{1}{|c}{$n=300$} & \multicolumn{1}{|c|}{$\theta_{2}$} & .24 &
.03 & .08 & .07 & .07 & .08 & .98\\\cline{2-9}%
\multicolumn{1}{|c}{} & \multicolumn{1}{|c|}{$\theta_{1}$} & .41 & .02 & .12 &
.13 & .13 & .11 & .97\\\hline
\multicolumn{1}{|c}{$n=1000$} & \multicolumn{1}{|c|}{$\theta_{2}$} & .24 &
.02 & .04 & .04 & .03 & .00 & .98\\\cline{2-9}%
\multicolumn{1}{|c}{} & \multicolumn{1}{|c|}{$\theta_{1}$} & .42 & .01 & .07 &
.07 & .06 & .00 & .97\\\hline
\multicolumn{1}{|c}{$n=10000$} & \multicolumn{1}{|c|}{$\theta_{2}$} & .24 &
.02 & .01 & .01 & .01 & .00 & .81\\\cline{2-9}%
\multicolumn{1}{|c}{} & \multicolumn{1}{|c|}{$\theta_{1}$} & .42 & .01 & .02 &
.02 & .02 & .00 & .95\\\hline
\end{tabular}

\begin{tabular}
[c]{cc|c|c|c|c|c|c|c|}%
\multicolumn{9}{c}{Table 3: $b(X)$ Linear, Squares, and Interactions}%
\\\cline{3-9}
&  & PI Bias & DB Bias & Med SE & PI SD & DB SD & PI Cvg & DB Cvg\\\hline
\multicolumn{1}{|c}{$n=100$} & \multicolumn{1}{|c|}{$\theta_{2}$} & .16 &
.01 & .13 & .11 & .38 & .85 & .98\\\cline{2-9}%
\multicolumn{1}{|c}{} & \multicolumn{1}{|c|}{$\theta_{1}$} & .26 & .03 & .21 &
.19 & .29 & .82 & .96\\\hline
\multicolumn{1}{|c}{$n=300$} & \multicolumn{1}{|c|}{$\theta_{2}$} & .15 &
.02 & .07 & .07 & .07 & .51 & .98\\\cline{2-9}%
\multicolumn{1}{|c}{} & \multicolumn{1}{|c|}{$\theta_{1}$} & .24 & .01 & .12 &
.12 & .12 & .52 & .97\\\hline
\multicolumn{1}{|c}{$n=1000$} & \multicolumn{1}{|c|}{$\theta_{2}$} & .14 &
.01 & .03 & .03 & .03 & .04 & .99\\\cline{2-9}%
\multicolumn{1}{|c}{} & \multicolumn{1}{|c|}{$\theta_{1}$} & .23 & .01 & .07 &
.07 & .06 & .07 & .97\\\hline
\multicolumn{1}{|c}{$n=10000$} & \multicolumn{1}{|c|}{$\theta_{2}$} & .13 &
.01 & .01 & .01 & .01 & .00 & .98\\\cline{2-9}%
\multicolumn{1}{|c}{} & \multicolumn{1}{|c|}{$\theta_{1}$} & .23 & .01 & .02 &
.02 & .02 & .00 & .94\\\hline
\end{tabular}

\end{center}

In all cases debiased GMM has much smaller bias than the plug-in estimator.
For the richest dictionary $b(X)$ in Table 3 coverage probabilities are quite
close to the nominal value though conservative. In contrast plug-in GMM has
large bias and confidence interval coverage probabilities that are far from
their nominal values in all cases. Remarkably, for larger sample sizes or
smaller dimensional $b(x)$, debiased GMM is no more variable than plug-in GMM,
and in several cases is less variable. Overall, the performance of the
debiased GMM estimator in this example with an "off the shelf" machine learner
suggests that debiased GMM for dynamic discrete choice and other structural
models could be useful in practice.

The low variance of debiased GMM could result partly from the fact that the
bias correction "partials out" the effect of $\hat{\gamma}_{\ell}$ in the
identifying moments, with the effect of changing $\gamma$ in the identifying
moment functions being approximately cancelled by the effect of $\hat{\gamma
}_{\ell}$ in the nonparametric influence function estimator. In the next
Section we explain this "partialling out" effect.

\section{Neyman Orthogonality}

Neyman orthogonality refers to the unknown functions $\gamma$ and $\alpha$
having no first order effect on the moments%
\[
\bar{\psi}(\gamma,\alpha,\theta):=E[\psi(W,\gamma,\alpha,\theta)]\text{.}%
\]
To show Neyman orthogonality let $\alpha(F)$ denote the probability limit of
$\hat{\alpha}$ when $F$ is the CDF of $W,$ similarly to $\gamma(F).$ Because
$\phi(w,\gamma,\alpha,\theta)$ is the nonparametric influence function it will
satisfy the mean zero condition in equation (\ref{infdef}) identically in $F$,
so that $0\equiv E_{F}[\phi(W,\gamma(F),\alpha(F),\theta)].$ Substituting
$F_{\tau}$ for $F$ and differentiating this identity with respect to $\tau$
gives%
\begin{align}
0  &  =\int\phi(w,\gamma_{0},\alpha_{0},\theta)H(dw)+\frac{\partial}%
{\partial\tau}E[\phi(W,\gamma(F_{\tau}),\alpha(F_{\tau}),\theta)]
\label{lrdef}\\
&  =\frac{\partial}{\partial\tau}E[g(W,\gamma(F_{\tau}),\theta)]+\frac
{\partial}{\partial\tau}E[\phi(W,\gamma(F_{\tau}),\alpha(F_{\tau}%
),\theta)]=\frac{\partial}{\partial\tau}\bar{\psi}(\gamma(F_{\tau}%
),\alpha(F_{\tau}),\theta),\nonumber
\end{align}
where the first equality follows by the chain rule, the second equality
follows from the influence function formula in equation (\textit{\ref{infdef}%
}), and the third from the definition of $\bar{\psi}(\gamma,\alpha,\theta).$
This equation shows that the functions $\gamma$ and $\alpha$ have no first
order effect on $\bar{\psi}(\gamma,\alpha,\theta)$ along the path
$(\gamma(F_{\tau}),\alpha(F_{\tau}))$. The second equality shows how the
presence of $E[\phi(W,\gamma,\alpha,\theta)]$ "partials out" the effect of
varying $\tau$ on $E[g(W,\gamma(F_{\tau}),\theta)]$. The zero mean property of
the nonparametric influence function implies that the local effect of $\gamma$
on $E[g(W,\gamma,\theta_{0})]$ along the path $\gamma(F_{\tau})$ is cancelled,
or "partialled out," by the effect of varying $\gamma$ and $\alpha$ on
$E[\phi(W,\gamma,\alpha,\theta)]$ along the path $\gamma(F_{\tau})$ and
$\alpha(F_{\tau}).$ The following result gives precise conditions for equation
(\ref{lrdef}):

\bigskip

\textsc{Theorem 1:}\textit{ If i) equation (\ref{infdef}) is satisfied; ii)
}$\int\phi(w,\gamma(F_{\tau}),\alpha(F_{\tau}),\theta)F_{\tau}(dw)=0$\textit{
for all }$\tau\in\lbrack0,\bar{\tau})$ \textit{with} $\bar{\tau}>0,$\textit{
and iii) }$\int\phi(w,\gamma(F_{\tau}),\alpha(F_{\tau}),\theta)F_{0}%
(dw)$\textit{ and }$\int\phi(w,\gamma(F_{\tau}),\alpha(F_{\tau}),\theta
)H(dw)$\textit{ are continuous at }$\tau=0$\textit{ then equation
(\ref{lrdef}) is satisfied.}

\bigskip

The proofs of this result and others are given in Appendix A.

\bigskip

\textsc{Example 1: }Here $g(w,\gamma,\theta)=z\gamma(x)-\theta$ and
$\phi(w,\gamma,\alpha)=\alpha(x)[y-\gamma(x)]$ so that%
\begin{align*}
\bar{\psi}(\gamma,\alpha,\theta)  &  =E[\psi(W,\gamma,\alpha,\theta
)]=E[Z\gamma(X)]-\theta+E[\alpha(X)\{Y-\gamma(X)\}]\\
&  =E[\alpha_{0}(X)\gamma(X)]-\theta+E[\alpha(X)\{\gamma_{0}(X)-\gamma(X)\}].
\end{align*}
When $\alpha(X)=\alpha_{0}(X)$ the presence of $\gamma$ in the identifying
moment $E[g(W,\gamma,\theta)]$ is exactly cancelled, or partialled out, by the
presence of $\gamma$ in $E[\phi(W,\gamma,\alpha_{0})]=E[\alpha_{0}%
(X)\{\gamma_{0}(X)-\gamma(X)\}].$

\bigskip

Equation (\textit{\ref{lrdef}}) is a total zero derivative condition for joint
variation in $\left(  \gamma,\alpha\right)  $ along the path $(\gamma(F_{\tau
}),\alpha(F_{\tau})).$ In many cases $\gamma(F)$ and $\alpha(F)$ are distinct
objects so that it is possible to choose $F_{\tau}$ so that $\alpha(F_{\tau})$
varies with $\tau$ and $\gamma(F_{\tau})=\gamma_{0}$ remains equal to its true
value. For example $\gamma(F)$ and $\alpha(F)$ may be determined by the
distributions of different random variables and so be distinct objects. In
such cases equation (\ref{infdef}) implies that $\phi(W,\gamma_{0}%
,\alpha,\theta)$ has mean zero even when $\alpha\neq\alpha_{0}.$

\bigskip

\textsc{Theorem 2:}\textit{ For any }$\alpha$ if \textit{i) there is
}$F_{\alpha}$\textit{ such that }$\alpha\left(  F_{\alpha}\right)  =\alpha
$\textit{ and }$\gamma(F_{\tau}^{\alpha})=\gamma_{0}$ \textit{for }$F_{\tau
}^{\alpha}=(1-\tau)F_{\alpha}+\tau F_{0}$ \textit{and all }$\tau\in
\lbrack0,\bar{\tau}),$ $\bar{\tau}>0$; \textit{ii) }$d\int g(w,\gamma(F_{\tau
}^{\alpha}),\theta)F_{\alpha}(dw)/d\tau=\int\phi(w,\gamma_{0},\alpha
,\theta)F_{0}(dw)$ \textit{then}%
\begin{equation}
E[\phi(W,\gamma_{0},\alpha,\theta)]=0. \label{nif zero mean}%
\end{equation}

\bigskip

Noting that hypothesis ii) is just the characterization of the nonparametric
influence function when the true distribution is $F_{\alpha}$, we see that the
nonparametric influence function has zero expectation when the function
$\alpha$ and/or the parameter $\theta$ are not equal to their true values
$\alpha_{0}$ and $\theta_{0}$ and there is some distribution $F_{\alpha}$ such
that $\alpha=\alpha(F_{\alpha})$. In all the examples of which we are aware
equation (\ref{nif zero mean}) is easy to confirm by inspection of
$\phi(W,\gamma_{0},\alpha,\theta)$.

\bigskip

\textsc{Example 1:} Here $E[\phi(W,\gamma_{0},\alpha,\theta)]=E[\alpha
(X)\{Y-\gamma_{0}(X)\}]=0$ for any $\alpha(X)$ by $\gamma_{0}(X)=E[Y|X]$ and
iterated expectations.

\bigskip

Theorem 2 shows that equation (\ref{nif zero mean}) is a general property of
the nonparametric influence function and is not confined to a particular set
of examples.

Neyman orthogonality with respect to only $\gamma(F_{\tau})$ follows from
choosing $F_{\tau}$ so that $\alpha(F_{\tau})=\alpha_{0}$ when $\gamma(F)$ and
$\alpha(F)$ are distinct objects. Then equation (\ref{lrdef}) implies
$\partial\bar{\psi}(\gamma(F_{\tau}),\alpha_{0},\theta_{0})/\partial\tau=0$.
For asymptotic theory it is useful to have a zero derivative with respect to
$\gamma$ that is stronger than this pathwise derivative condition. The
following result shows that equation (\ref{lrdef}) implies that $\bar{\psi
}(\gamma,\alpha_{0},\theta_{0})$ has a zero Hadamard derivative with respect
to $\gamma$ if the set of pathwise derivatives $d\gamma(F_{\tau})/d\tau$ is
rich enough.

\bigskip

\textsc{Theorem 3:} \textit{If there is a norm }$\left\Vert \gamma\right\Vert
$\textit{, a linear set }$\Gamma,$\textit{ and a set }$\mathcal{H}%
$\textit{\ such that for all }$H\in\mathcal{H}$\textit{; i) }$\alpha(F_{\tau
})=\alpha_{0}$\textit{ and equation (\ref{lrdef}) is satisfied; ii) }%
$\bar{\psi}(\gamma,\alpha_{0},\theta_{0})$\textit{ is Hadamard differentiable
at }$\gamma_{0}$\textit{ tangentially to }$\Gamma$\textit{ with derivative}
$\bar{\psi}_{\gamma}(\delta,\alpha_{0},\theta_{0}),$\textit{ }$\delta\in
\Gamma$; \textit{iii) }$\gamma(F_{\tau})$\textit{ is Hadamard differentiable
at }$\tau=0;$ \textit{iv) the closure of }$\{\partial\gamma(F_{\tau}%
)/\partial\tau:H\in$\textit{ }$\mathcal{H}\}$\textit{ is }$\Gamma$\textit{
then }%
\begin{equation}
\bar{\psi}_{\gamma}(\delta,\alpha_{0},\theta_{0})=0,\text{ }\delta\in\Gamma.
\label{orthogonality}%
\end{equation}
\textit{Furthermore, if }$\bar{\psi}(\gamma,\alpha_{0},\theta_{0})$\textit{ is
twice continuously Frechet differentiable in a neighborhood of }$\gamma_{0}%
$\textit{ then there is }$C>0$\textit{ such that for }$\left\Vert
\gamma-\gamma_{0}\right\Vert $\textit{ small enough}%
\begin{equation}
\left\Vert \bar{\psi}(\gamma,\alpha_{0},\theta_{0})\right\Vert \leq
C\left\Vert \gamma-\gamma_{0}\right\Vert ^{2}. \label{orth rem}%
\end{equation}

\bigskip

Hadamard and Frechet differentiability are defined and discussed e.g. in van
der Vaart (1998, 20.2). The first conclusion of this result is that the
expected value $\bar{\psi}(\gamma,\alpha_{0},\theta_{0})$ of $\psi
(w,\gamma,\alpha,\theta)$ moment function has zero Hadamard derivative with
respect to $\gamma$. Theorems 2 and 3 combined show that adding the
nonparametric influence function $\phi(w,\gamma,\alpha,\theta)$ to identifying
moment functions $g(w,\gamma,\theta)$ makes $\psi(w,\gamma,\alpha
,\theta)=g(w,\gamma,\theta)+$ $\phi(w,\gamma,\alpha,\theta)$ Neyman
orthogonal. Equation (\ref{orth rem}) bounds the departure from zero of the
expected moments. This bound is useful for formulating regularity conditions
for root-n consistency of debiased GMM when $\bar{\psi}(\gamma,\alpha
_{0},\theta_{0})$ is nonlinear in $\gamma,$ as we will explain in Section 8.
When formulating regularity conditions for particular moment functions and
first step estimators it may be simpler to directly confirm equation
(\ref{orth rem}). In many cases equation (\ref{orth rem}) will be satisfied
under specific regularity conditions when $\left\Vert \gamma-\gamma
_{0}\right\Vert $ is a mean square norm, because $\bar{\psi}(\gamma
,\alpha,\theta)$ is an expected value. Equation (\ref{orth rem}) with a mean
square norm $\left\Vert \gamma-\gamma_{0}\right\Vert $ has wide applicability
to machine learning first steps where many mean square convergence rates are available.

If $\gamma$ were taken to be the limit of a nonparametric estimator
$\hat{\gamma}$ for fixed bandwidth, number of series terms, or regularization
then we can think of $\left\Vert \gamma-\gamma_{0}\right\Vert $ as bias in
$\gamma$ from nonparametric estimation. Frechet differentiability of
$\bar{\psi}(\gamma,\alpha_{0},\theta_{0})$ in $\gamma$ and equation
(\ref{orthogonality}) then imply that $\bar{\psi}(\gamma,\alpha_{0},\theta
_{0})=o(\left\Vert \gamma-\gamma_{0}\right\Vert )$, so that the orthogonal
moments $\bar{\psi}(\gamma,\alpha_{0},\theta_{0})$ shrink to zero faster than
the nonparametric bias. Thus orthogonal moment functions have the small bias
property considered in Newey, Hsieh, and Robins (1998, 2004). As usual for GMM
the estimator $\hat{\theta}$ will inherit this property of the orthogonal
moment functions.

Equation (\ref{lrdef}) is similar to Theorem 2.2 of Robins et al. (2008) but
different in applying directly to the standard influence function
characterization in equation (\ref{infdef}) without specifying a score
function (derivative of the log-likelihood of a model). Proceeding in this way
allows us to show Neyman orthogonality of $\psi(w,\gamma,\alpha,\theta)$ using
equation (\ref{lrdef}). To the best of our knowledge equation (\ref{lrdef})
has not appeared in this form previously. Also, Theorems 1-3 appear to be
novel in specifying regularity conditions for Neyman orthogonality of moment
functions obtained from adding the nonparametric influence function.

The construction of orthogonal moment functions we consider has antecedents in
the literature on functional estimation, where the identifying moment
conditions are $g(w,\gamma,\theta)=m(\gamma)-\theta$ for some explicit
functional of $m(\gamma)$ of $\gamma.$ Here $\phi(w,\gamma,\alpha)$ is the
influence function of $m(\gamma(F))$ and $\psi(w,\gamma,\alpha,\theta
)=m(\gamma)-\theta+\phi(w,\gamma,\alpha)$. Examples of such moment functions
when $\gamma$ is a density function were given by Hasminskii and Ibragimov
(1978), Pfanzagl and Wefelmeyer (1981), and Bickel and Ritov (1988). Robins
and Rotnitzky (1992) gave an orthogonal moment function with identifying
moment function $m(w,\gamma)-\theta.$ Newey, Hsieh, and Robins (1998, 2004)
suggested adding the nonparametric influence to identifying moment functions
and showed that the remainder will be second order. Robins et al. (2008), and
Van der Vaart (2014) showed that in general the remainder is second order for
functional estimators. GMM with moment functions obtained by adding the
nonparametric influence to identifying moment functions was considered in
Chernozhukov et al. (2018) and Bravo, Escanciano, and van Keilegom (2020). The
framework described in those papers originated in the joint research for this paper.

Robustness of the nonparametric influence function to the additional functions
on which it depends, as in Theorem 2, is not shown in any of the work cited in
the previous paragraph. The absence of a first order effect of $\gamma$ on
$\bar{\psi}(\gamma,\alpha,\gamma)$ was also not shown, although that is
implicit in the remainder analyses previously given. Theorems 2 and 3 provide
key conditions for first step estimation to have no first order effect on the
asymptotic variance of the debiased GMM estimator $\hat{\theta}$ as needed for
the asymptotic theory in Section 8.

The orthogonal moment functions $\psi(w,\gamma,\alpha,\theta)$ could be
constructed as an efficient influence function of a semiparametric model, as
in Robins and Rotnitzky (1992) and many others since. The construction we give
bypasses the semiparametric efficiency framework and is based on the simpler
influence function characterization in equation (\ref{infdef}) and the limit
$\gamma(F)$ for any semiparametric estimator as in Newey (1994a). This
construction highlights the distinct roles of $g(w,\gamma,\theta)$ as
identifying moments and $\phi(w,\gamma,\alpha,\theta)$ as a bias correction
that does not affect identification and is entirely determined by
$g(w,\gamma,\theta)$ and $\gamma(F)$. The distinct role of $\phi
(w,\gamma,\alpha,\theta)$ leads directly to the robustness result of Theorem 2
and motivates its automatic estimation in Section 6. For these reasons we
choose to construct orthogonal moments by adding the nonparametric influence
function to identifying moments rather than finding a semiparametric efficient
influence function.

The orthogonalization given here is nonparametric, i.e. is estimator based
rather than model based, relying only on $g(w,\gamma,\theta)$ and the
nonparametric limit $\gamma(F)$ of $\hat{\gamma}$ and not on the specification
of any model. Consequently, the Neyman orthogonality shown in Theorems 1-3
does not depend on correct specification of any model. For any first step the
orthogonal moments $\bar{\psi}(\gamma,\alpha,\theta)$ will have the properties
given in the conclusions of Theorems 2 and 3 under the stated regularity
conditions. The absence of any model specification from those regularity
conditions demonstrates the nonparametric nature of these orthogonal moment conditions.

There are also model based approaches to orthogonalization that generalize the
efficient influence function. For example Newey (1990), Belloni, Chernozhukov,
and Kato (2015), and Belloni et al. (2017) showed that in a semiparametric
model the residual from the projection of identifying moment functions on the
tangent set is orthogonal for certain kinds of first steps. The approach here
has the advantage that it does not depend on correct specification of a model.

In some cases the identifying moment functions $g(w,\gamma,\theta)$ may
already be orthogonal, so that $\phi(w,\gamma,\alpha,\theta)=0$. An important
class of orthogonal moment functions are those where $g(w,\gamma,\theta)$ is
the derivative with respect to $\theta$ of an objective function where
nonparametric parts have been concentrated out. That is, suppose that there is
a function $q(w,\theta,\zeta)$ such that $g(w,\gamma,\theta)=\partial
q(w,\theta,\zeta(\theta))/\partial\theta$ and $\zeta(\theta)=\arg\max_{\zeta
}E[q(W,\theta,\zeta)]$, where $\gamma$ includes $\zeta(\theta)$ and possibly
additional functions. Proposition 2 of Newey (1994a) and Lemma 2.5 of
Chernozhukov et al. (2018) then imply that $g(w,\gamma,\theta)$ is orthogonal.
This class of moment functions includes those of Robinson (1988), Ichimura
(1993), and various partially linear regression models where $\zeta$
represents a conditional expectation. It also includes the efficient score for
a semiparametric model when the nonparametric component estimates the maximum
of the expected log likelihood; see Severini and Wong (1992), Newey (1994a,
pp. 1358-1359), and van der Vaart (1998, pp. 391-396).

The nonparametric influence function $\phi(W,\gamma,\alpha,\theta)$ is unique
because $H$ in equation (\ref{lrdef}) is unrestricted except for regularity
conditions. This uniqueness provides another way to understand orthogonality
of $\psi(w,\gamma,\alpha,\theta)$. The influence function of $E_{F}%
[g(W,\gamma(F),\theta_{0})]$ at $F_{0}$ is $\psi(W,\gamma_{0},\alpha
_{0},\theta_{0}),$ under the moment condition $E[g(W,\gamma_{0},\theta
_{0})]=0.$ This result follows by differentiating $E_{F_{\tau}}[g(W,\gamma
(F_{\tau}),\theta_{0})]$ with respect to $\tau$, applying the chain rule, and
using (\ref{lrdef}). The sample average of the debiased moment function
$\hat{\psi}(\theta_{0})$ is a nonparametric estimator of $E_{F}[g(W,\gamma
(F),\theta_{0})].$ If $\hat{\psi}(\theta_{0})$ is asymptotically equivalent to
a sample average and locally regular, meaning that for $H$ and data are i.i.d.
with CDF $F_{\tau_{n}}=(1-\tau_{n})F_{0}+\tau_{n}H$ and $\tau_{n}=O(1/\sqrt
{n})$ the limiting distribution of $\hat{\psi}(\theta_{0})-E_{F_{\tau_{n}}%
}[g(W,\gamma(F_{\tau_{n}}),\theta_{0})]$ does not depend on $\tau_{n},$, then
uniqueness of $\phi$ implies%
\begin{equation}
\frac{1}{n}\sum_{\ell=1}^{L}\sum_{i\in I_{\ell}}[g(W_{i},\hat{\gamma}_{\ell
},\theta_{0})+\phi(W_{i},\hat{\gamma}_{\ell},\hat{\alpha}_{\ell},\tilde
{\theta}_{\ell})]=\hat{\psi}(\theta_{0})=\frac{1}{n}\sum_{i=1}^{n}\psi
(W_{i},\gamma_{0},\alpha_{0},\theta_{0})+o_{p}(n^{-1/2}), \label{asym lin}%
\end{equation}
as in Van der Vaart (1991), Newey (1994a), and Ichimura and Newey (2017). In
other words, the \textit{only} sample average of a function of the data that
$\hat{\psi}(\theta_{0})$ can be asymptotically equivalent to, and also be
locally regular, is $\sum_{i=1}^{n}\psi(W,\gamma_{0},\alpha_{0},\theta
_{0})/n.$ Equation (\ref{asym lin}) is precisely an asymptotic version of
orthogonality where the sample average $\hat{\psi}(\theta_{0})$ is
asymptotically equivalent to the sample average of the same function with
estimators $\hat{\gamma}_{\ell},$ $\hat{\alpha}_{\ell},$ and $\tilde{\theta
}_{\ell}$ replaced by their limits. In this way orthogonality is justified by
uniqueness of $\phi(w,\gamma,\alpha,\theta)$.

Equation (\ref{asym lin}) is the key orthogonality property for the asymptotic
theory of Section 8, where primitive regularity conditions for equation
(\ref{asym lin}) are given. The primitive conditions are motivated by the
analysis of this Section, including the conclusions of Theorems 2 and 3.
Equation (\ref{asym lin}) also helps in the comparison of debiased GMM with
plug-in GMM, to which we now turn.

\section{Comparing Debiased and Plug-in GMM}

To highlight the role of orthogonal moment functions we compare the properties
of debiased GMM with a corresponding cross-fit plug-in GMM estimator%
\[
\tilde{\theta}=\arg\min_{\theta\in\Theta}\hat{g}(\theta)^{\prime}\hat
{\Upsilon}\hat{g}(\theta).
\]
Confidence intervals based on the plug-in GMM estimator are invalid with first
step model selection and plug-in GMM is so biased that it is not root(n)
consistent with a Lasso and other regularized first steps. Debiased GMM does
not suffer from these problems. Plug-in GMM is simpler than debiased GMM in
not requiring computation of $\hat{\alpha}_{\ell}$ and $\hat{\phi}$ but that
reduced computational does not justify it when it has large biases due to
model selection and/or regularization. In this Section we discuss these and
other comparative properties of debiased and plug-in GMM.

To compare the properties of debiased and plug-in GMM it is helpful to compare
the key asymptotic property of debiased GMM in equation (\ref{asym lin}) with
a corresponding key property of plug-in GMM,%
\begin{equation}
\hat{g}(\theta_{0})=\frac{1}{n}\sum_{i=1}^{n}\psi(W_{i},\gamma_{0},\alpha
_{0},\theta_{0})+o_{p}(n^{-1/2}). \label{PIE key}%
\end{equation}
This property is very important. When equation (\ref{PIE key}) is not
satisfied plug-in GMM can have invalid confidence intervals or not be root(n)
consistent. As shown in Section 8, the corresponding condition in equation
(\ref{asym lin}) for debiased GMM will be satisfied under general and simple
regularity conditions. In contrast, equation (\ref{PIE key}) requires an
additional condition that is specific to the first step and more complicated.
Let%
\begin{equation}
\tilde{\phi}=\frac{1}{n}\sum_{\ell=1}^{L}\sum_{i\in I_{\ell}}\phi(W_{i}%
,\hat{\gamma}_{\ell},\alpha_{0},\theta_{0}). \label{PIE remainder}%
\end{equation}

\textsc{Theorem 4:} \textit{If equation (\ref{asym lin}) is satisfied for
}$\hat{\alpha}_{\ell}=\alpha_{0}$\textit{ and }$\tilde{\theta}_{\ell}%
=\theta_{0}$\textit{ then equation (\ref{PIE key}) is satisfied if and only if
}$\tilde{\phi}=o_{p}(n^{-1/2})$.

\bigskip

Equation (\ref{asym lin}) for $\hat{\alpha}_{\ell}=\alpha_{0}$ and
$\tilde{\theta}_{\ell}=\theta_{0}$ is even simpler and more general than
(\ref{asym lin}) for a debiased GMM estimator, making $\sqrt{n}\tilde{\phi
}\overset{p}{\longrightarrow}0$ the key regularity condition for plug in GMM.
This condition will fail with first step model selection or regularization.

\bigskip

\textsc{Example 1:} Here $\psi(W,\hat{\gamma}_{\ell},\alpha_{0},\theta
_{0})-\psi(W,\gamma_{0},\alpha_{0},\theta_{0})=[Z-\alpha_{0}(X)][\hat{\gamma
}_{\ell}(X)-\gamma_{0}(X)]$, so that equation (\textit{\ref{asym lin}}) will
be satisfied under a weak, mean square consistency condition for $\gamma_{0},$
as in Assumption 1 of Section 8. Then equation \textit{(\ref{PIE key})} will
be satisfied if and only if%
\begin{equation}
\sqrt{n}\tilde{\phi}=\frac{1}{\sqrt{n}}\sum_{\ell=1}^{L}\sum_{i\in I_{\ell}%
}\alpha_{0}\left(  X_{i}\right)  [Y_{i}-\hat{\gamma}_{\ell}(X_{i}%
)]\overset{p}{\longrightarrow}0. \label{ex1plugin}%
\end{equation}
This condition is approximate sample orthogonality of the cross-fit residual
$Y_{i}-\hat{\gamma}_{\ell}(X_{i})$ with $\alpha_{0}(X_{i}).$ This condition
can fail to hold under model selection, where $\hat{\gamma}_{\ell}(X)$ only
includes some variables on which $\alpha_{0}(X)$ depends, so that the
residuals $Y_{i}-\hat{\gamma}_{\ell}(X_{i})$ are not approximately orthogonal
to $\alpha_{0}(X_{i})$ in the sample. This condition can also fail to hold
when $\hat{\gamma}_{\ell}(X_{i})$ is a regularized estimator where residuals
are not approximately orthogonal to functions of $X_{i}$ on which $\alpha
_{0}(X_{i})$ depends. Thus we see that the properties of the plug-in estimator
will depend heavily on the nature of $\alpha_{0}(X_{i})$ even though no
estimator of $\alpha_{0}(X_{i})$ is explicitly present.

\subsection{Plug-in Estimators with Model Selection Give Invalid Confidence
Intervals}

If the first step $\hat{\gamma}$ incorporates model selection then plug-in GMM
gives invalid asymptotic confidence intervals. This occurs because model
selection that leads to a correct model with probability approaching must also
choose the same model under local alternatives in a root-n neighborhood where
that model is incorrect. This local misspecification has a first order effect
on plug-in GMM, giving a limiting distribution with nonzero mean under local
alternatives, so that the usual asymptotic confidence intervals, based on a
zero mean limiting distribution, are invalid. Debiased GMM does not suffer
from this problem because the bias from the incorrect first step is
second-order. This feature of debiased GMM makes it preferred to plug-in GMM
in the many applications where the first step incorporates model selection.

\bigskip

\textsc{Example 1:} Here plug-in GMM is $\tilde{\theta}=\sum_{\ell=1}^{L}%
\sum_{i\in I_{\ell}}Z_{i}\hat{\gamma}_{\ell}(X_{i})/n$. Suppose that
$\gamma_{0}(x)=\breve{b}(x)^{\prime}\breve{\beta}$ for finite dimensional
vectors $\breve{b}(x)$ and $\breve{\beta},$ where $\breve{\beta}$ has no zero
components. Consider $\hat{\gamma}_{\ell}$ that incorporates model selection
such that $\hat{\gamma}_{\ell}(x)$ is equal to the least squares regression of
$Y_{i}$ on $\breve{b}(X_{i})$ for all $i\in I_{\ell}$ with probability
approaching one. A variety of estimators have this property, including post
Lasso or post Lasso with thresholding, when $\breve{b}(X_{i})$ are included
among the Lasso regressors and are not too correlated with the other Lasso
regressors. Here equation (\ref{ex1plugin}) will not hold when $\alpha_{0}(X)$
is not a linear combination of $\breve{b}(X_{i})$, because $\alpha_{0}(X_{i})$
need not be approximately orthogonal in the sample to the residuals
$Y_{i}-\hat{\gamma}_{\ell}(X_{i})$ in that case. The precise behavior of
$\tilde{\theta}$ is given in the follow result.

\bigskip

\textsc{Theorem 5:}\textit{ If i) }$\gamma_{0}(x)=\breve{b}(x)^{\prime}%
\breve{\beta};$\textit{ ii) }$\breve{G}=:E[\breve{b}(X)\breve{b}(X)^{\prime}%
]$\textit{ is nonsingular; iii) }$\alpha_{0}(X)$\textit{ and }$\breve{b}%
(X)$\textit{ are bounded; and iv) with probability approaching one }%
$\hat{\gamma}_{\ell}(x)$\textit{ is equal to ordinary least squares from
regressing }$Y_{i}$\textit{ on }$\breve{b}(X_{i})$\textit{ over }$i\notin
I_{\ell};$\textit{ then for }$\bar{\alpha}(x)=\breve{b}(x)^{\prime}\breve
{G}^{-1}E[\breve{b}(X)\alpha_{0}(X)]$\textit{,}%
\[
\tilde{\theta}=\theta_{0}+\frac{1}{n}\sum_{i=1}^{n}\zeta(W_{i})+o_{p}%
(n^{-1/2})\text{, }\zeta(W)=Z\gamma_{0}(X)-\theta_{0}+\bar{\alpha}%
(X)[Y-\gamma_{0}(X)].\text{ }%
\]
\textit{Also if iv) the distribution of }$Y$\textit{ conditional on }%
$(X,Z)$\textit{ has a pdf }$f_{0}(y|x,z)$\textit{ such that there is }%
$C>0$\textit{ with }$E[\int\{\sup_{\left\vert a\right\vert \leq C}%
\{[df_{0}(y+a|X,Z)/da]^{2}/f_{0}(y+a|X,Z)\}\}dy]<\infty,$\textit{ then for
}$\bar{\sigma}^{2}=E[\{\alpha_{0}(X)-\bar{\alpha}(X)\}^{2}],$\textit{ any
}$\mu$, \textit{and }$W_{1},...,W_{n}$ \textit{i.i.d. with CDF }$F_{n}$
\textit{having conditional pdf }$f_{0}(\tilde{y}-n^{-1/2}\mu\{\alpha
_{0}(X)-\bar{\alpha}(X)\}|x,z)$ \textit{for }$Y_{i}$ given $(X_{i}%
,Z_{i})=(x,z)$ \textit{and CDF }$F_{0}(x,z)$ \textit{for }$(X_{i},Z_{i})$
\textit{we have}%
\[
\sqrt{n}(\tilde{\theta}-\theta_{n})\overset{d}{\longrightarrow}N(\mu
\bar{\sigma}^{2},V),
\]
\textit{where }$\theta_{n}=E_{F_{n}}[Z\cdot E_{F_{n}}[Y|X]]$\textit{ is the
parameter of interest for }$F_{n}.$

\bigskip

The first conclusion shows that with model selection plug-in GMM is
asymptotically equivalent to a sample average but the function being averaged
is not $\psi(W,\gamma_{0},\alpha_{0},\theta_{0})$ if $\bar{\alpha}%
(X)\neq\alpha_{0}(X)$. As a result the second conclusion follows, with the
limiting distribution of plug-in GMM having a nonzero mean under a local
alternative when $\bar{\alpha}(X)\neq\alpha_{0}(X)$, so that the usual
asymptotic confidence interval based on a zero mean limiting distribution is
invalid. The first conclusion leads to the second conclusion because the
nonparametric influence function is unique. The asymptotic equivalence of
plug-in GMM to a sample average that is not $\psi(W,\gamma_{0},\alpha
_{0},\theta_{0})$ implies that there are local alternatives where the limiting
distribution has non zero mean. In this way the invalidity of standard
asymptotic confidence intervals for plug-in GMM $\tilde{\theta}$ under model
selection can be thought of as resulting from asymptotic equivalence of
$\sqrt{n}(\tilde{\theta}-\theta_{0})$ to a sample average of a function that
is not $\psi(W,\gamma_{0},\alpha_{0},\theta_{0})$.

The reason that model selection creates a problem for plug-in GMM is that
$\psi(W,\gamma_{0},\alpha_{0},\theta_{0})$ depends on two functions, one being
$\gamma_{0}(X)$ and the other $\alpha_{0}(X).$ Model selection for
$\hat{\gamma}(x)$ gives a good estimator of $\gamma_{0}(X)$ but the variables
selected for estimating $\gamma_{0}(X)$ may not include all the variables on
which $\alpha_{0}(X)$ depends, leading to $\bar{\alpha}(X)\neq\alpha_{0}(X)$
and hence $\bar{\sigma}_{{}}^{2}>0.$ As mentioned previously in this Section
$\sqrt{n}\tilde{\phi}\overset{p}{\longrightarrow}0$ fails in this case because
$\alpha_{0}(X_{i})$ will not be approximately asymptotically orthogonal in the
sample to the cross-fit least square residuals.

Plug-in GMM will be asymptotically equivalent to the sample average of
$\psi(W,\gamma_{0},\alpha_{0},\theta_{0})$ when $\bar{\alpha}(X)=\alpha
_{0}(X)$, i.e. when $\alpha_{0}(X)$ is a linear combination of the same
variables on which $\gamma_{0}(X)$. However, there is generally no good reason
to impose that condition a priori. The one case we are aware of where that
condition does hold a priori is where $Z=Y$, so that $\alpha_{0}(X)=\gamma
_{0}(X).$ In this case the parameter of interest is $\theta_{0}=E[\gamma
_{0}(X)^{2}]=E[\{E[Y|X]\}^{2}]$ and plug-in asymptotic confidence intervals
would be valid.

\bigskip

In general plug-in GMM estimators will have similar properties to this
example. The nonparametric influence function $\phi(W,\gamma_{0},\alpha
_{0},\theta_{0})$ depends on $\alpha_{0}$ in addition to the first step
function $\gamma_{0}$, and a first step model selection will only fit
$\gamma_{0}$ well and not $\alpha_{0}$. Consequently model selection for the
first step will make asymptotic confidence intervals invalid for plug-in GMM.
There may be some exceptions where this model selection problem does not hold,
as discussed in Example 1, but these are few and far between.

In some cases the double selection procedure of Belloni, Chernozhukov, and
Hansen (2014) also corrects for model selection. That procedure includes
variables in the estimation of $\gamma_{0}$ if they are also important in
estimation of $\alpha_{0}$. Such a double selection procedure is not available
for many first step machine learners that are not constructed by selecting
variables in a regression. Also, double selection can lead to much less
parsimonious models than would be used by debiased GMM, which allows different
variables to be selected for estimation of $\alpha_{0}$ than are selected for
estimation of $\gamma_{0}$. It is beyond the scope of this paper to compare
debiased GMM with double selection.

Another way to avoid the model selection problem for plug-in GMM is to limit
selection to models that can approximate any unknown function. For example
Newey (1994a) did this for series estimation by requiring that selection is
made only among models that can approximate any function in large samples.
Forcing a flexible approximation in this way is not very feasible in high
dimensional settings and is not needed for debiased GMM, where model selection
can be applied separately to the various first step estimators.

That model selection can make confidence intervals invalid was shown by Leeb
and Potscher (2005, 2008) for least squares and Lasso. The model selection
problem for plug-in GMM was pointed out by Belloni, Chernozhukov, and Kato
(2015) and Chernozhukov, Hansen, and Spindler (2015) who suggested Neyman
orthogonal moment functions as a solution to this problem. Theorem 4 gives
precise asymptotic theory showing the asymptotic equivalence of plug-in GMM to
a sample average of a function other than $\psi(W,\gamma_{0},\alpha_{0}%
,\theta_{0}).$ This feature of plug-in GMM under model selection is
reminiscent of the Hodges estimator, a well known example of an estimator that
is not asymptotically equivalent to the sample average of the influence
function of the mean.

Debiased GMM avoids the model selection problem because it has the small bias
property discussed in Section 4 following Theorem 3. Model selection imparts a
bias to the first step of size $1/\sqrt{n}$ under local alternatives. The
small bias property for debiased GMM means that its bias vanishes faster than
first step bias, i.e. faster than $1/\sqrt{n}$. Hence the asymptotic
distribution of debiased GMM will have zero mean and asymptotic confidence
intervals will be valid. Given the common feature of model selection in
machine learning and other econometric methods the robustness of debiased GMM
to first step model selection motivates its use in practice.

\subsection{Plug-in GMM Will Not Be Root(n) Consistent for Lasso First Steps}

Many machine learners employ regularization to obtain estimators of functions
that approximately balance bias and standard deviation. For nonparametric
estimation or machine learning with large sets of predictors the standard
deviation of the predictor will shrink slower than $1/\sqrt{n},$ and hence so
will be bias. This bias may pass through to plug in GMM and result in
$\tilde{\theta}$ not being root(n) consistent. This bias problem is clearly
present for Lasso where penalization leads to bias for $\tilde{\theta}$ of
size $\sqrt{\ln(p)/n}$. With bias of that size $\sqrt{n}$ times the bias will
be of size $\sqrt{\ln(p)}$ which goes to infinity, so that that plug-in GMM is
not root-n consistent, as we show in this subsection. Debiased GMM has the
small bias property discussed in Section 2 and so will be root-n consistent
under sufficient regularity conditions, with bias being second order (size
$\ln(p)/n$ for Lasso) permitting debiased GMM to be root-n consistent (by
$\sqrt{n}\ln(p)/n\longrightarrow0$ for Lasso).

\bigskip

\textsc{Example 1: }To illustrate the regularization bias problem for plug-in
GMM we give its properties when $\gamma_{0}(X)$ is a linear combination of a
finite number of functions and $\hat{\gamma}_{\ell}$ is Lasso;
\[
\hat{\gamma}_{\ell}(x)=b(x)^{\prime}\hat{\beta}_{\ell},\text{ }\hat{\beta
}_{\ell}=\arg\min_{\beta}\frac{1}{n-n_{\ell}}\sum_{i\notin I_{\ell}}%
[Y_{i}-b(X_{i})^{\prime}\beta]^{2}+r\sum_{j=1}^{p}\left\vert \beta
_{j}\right\vert ,
\]
where $b(x)$ is a $p\times1$ vector of functions.

\bigskip

\textsc{Theorem 6: }\textit{If i) there is an }$s\times1$\textit{ subvector
}$\breve{b}(x)$\textit{ of }$b(x)$\textit{ such that }$\gamma_{0}(X)=\breve
{b}(X)^{\prime}\breve{\beta},$\textit{ all elements of }$\breve{\beta}%
$\textit{ are nonzero, }$\breve{G}=:E[\breve{b}(X)\breve{b}(X)^{\prime}%
]$\textit{ is nonsingular; ii) }$\alpha_{0}(X)$\textit{ is bounded and
}$c=E[\alpha_{0}(X)\breve{b}(X)^{\prime}]\breve{G}^{-1}\breve{e}\neq0$\textit{
for }$\breve{e}=(sgn(\breve{\beta}_{1}),...,sgn(\breve{\beta}_{s}))^{\prime}%
$\textit{; iii) with probability approaching one }$1(\hat{\beta}_{\ell
j}=0)=1(\beta_{j}=0)$\textit{ for all }$j;$\textit{ iv) }$\sqrt{\ln
(p)/n}=O(r)$\textit{ then }$\tilde{\theta}=\sum_{\ell=1}^{L}\sum_{i\in
I_{\ell}}Z_{i}\hat{\gamma}_{\ell}(X_{i})/n$\textit{ satisfies}%
\[
\sqrt{n}\left\vert \tilde{\theta}-\theta_{0}\right\vert =\left\vert
O_{p}(1)-c\sqrt{n}r\right\vert \overset{p}{\longrightarrow}\infty.
\]

\bigskip

Because $\gamma_{0}(X)$ is a linear combination of a finite number of elements
of $b(X)$ condition iii) is known to be satisfied when all coefficients from
regressing each $b_{j}(X_{i})$ on $\breve{b}(X_{i})$ are small enough in
absolute value; see Zhao and Yu (2006).

\bigskip

Here we see that the plug in estimator is root-n consistent for Lasso in
Example 1 when $\gamma_{0}(X)$ is a linear combination of a finite number of
elements of $b(X)$ and condition iii)\ is satisfied. In general plug-in GMM
will not be root-n consistent with a Lasso first step, though it is beyond the
scope of this paper to show this. More generally plug-in GMM will also have
large bias for first step machine learners other than Lasso, e.g. as found for
random forests in a Monte Carlo example in Chernozhukov et al. (2018).

The robustness of debiased GMM to model selection and its low bias relative to
plug-in GMM are both reasons to prefer debiased GMM over plug-in GMM for an
any application where there is first step model selection or where the
regularization bias in the first step estimator passes through to the
estimator of the parameter of interest. Both features are present for many
machine learning first steps making debiased GMM especially useful there. We
note that this preference is based on first order properties which dominate
any second order comparison of debiased and plug-in GMM.

\subsection{Orthogonal Moments are Doubly Robust When They Are Affine in the
First Step}

Doubly robust moment conditions are those that hold when either $\gamma
=\gamma_{0}$ or $\alpha=\alpha_{0}$. Such moment conditions are of wide
interest. We will show in Section 7 that orthogonal moment conditions are
doubly robust if and only if $\bar{\psi}(\gamma,\alpha_{0},\theta_{0})$ is
affine in $\gamma.$ Plug in GMM moment conditions have more limited robustness
properties, being satisfied when $\gamma$ has a form related to $\alpha_{0}$
for some cases, as also discussed in Section 7.

\subsection{Debiased GMM May Have Better 2nd Order Properties Than Plug-In
GMM}

When there is no first step model selection and regularization does not
destroy root(n) consistency of plug-in GMM, debiased GMM has been shown to
have better 2nd order properties than plug-in GMM $\ $in some cases. Newey,
Hsieh, and Robins (1998, 2004) showed that debiased leave one out kernel
estimators of density weighted averages have smaller asymptotic and small
sample mean square error and than corresponding plug-ins, for a wide range of
bandwidths. Newey and Robins (2017) gave weaker conditions for root(n)
consistency of doubly robust estimators than for plug-in estimators. These
comparisons seem specific to the type of estimator, so it may be too much to
expect that debiased GMM\ always has better second order properties than
plug-in GMM. The first order advantages of debiased GMM when the first step
incorporates model selection is regularized give compelling reasons for its
use there.

\subsection{Asymptotic Theory for Debiased GMM is More General and Simple Than
for Plug-In GMM}

From Theorem 4 we see that the key difference between asymptotic theory for
debiased and plug-in GMM is that for debiased GMM it is sufficient that
equation (\ref{asym lin}) holds while for plug-in GMM it is sufficient that
(\ref{asym lin}) holds with $\hat{\alpha}_{\ell}=\alpha_{0}$ and
$\tilde{\theta}_{\ell}=\theta_{0}$ and that equation (\textit{\ref{PIE key}})
holds. Equation (\ref{asym lin}) is a little more involved for debiased GMM
because $\hat{\alpha}_{\ell}$ is estimated but sufficient conditions are
general and simple in only involving mean square consistency and one or two
mean square rate conditions, as described in Section 8. In contrast showing
equation (\textit{\ref{PIE key}}) for plug-in GMM seems very specific to the
nature of $\hat{\gamma}_{\ell}$ and quite complicated. Using cross-fitting
equation (\textit{\ref{PIE key}}) will follow from
\[
\sqrt{n}\int\phi(w,\hat{\gamma}_{\ell},\alpha_{0},\theta_{0})F_{0}%
(dw)\overset{}{\overset{p}{\longrightarrow}0,\text{ }(\ell=1,...,L)}%
\]
but showing this also seems quite specific to the form of $\hat{\gamma}$ and
complicated. In this way regularity conditions for plug-in GMM are less
general and more complicated than for debiased GMM.

\section{Automatic Estimation of $\alpha_{0}$}

Debiased GMM requires the nonparametric influence function $\phi
(w,\gamma,\alpha,\theta)$ and an estimator $\hat{\alpha}$ of the unknown
function $\alpha_{0}$. As previously discussed $\phi(w,\gamma,\alpha,\theta)$
is readily available in many important settings. In this Section we give a
general approach to constructing $\hat{\alpha}$ that is based only on the
orthogonal moment functions, show how that can be used for debiased GMM
estimation of objects that depend on location functions, and apply that to
obtain the $\hat{\alpha}_{\ell}$ of Example 2.

The orthogonality of $\bar{\psi}(\gamma,\alpha,\theta)$ can be used to
estimate $\alpha_{0}.$ Let $\Gamma$ be the linear set from Theorem 2 that is
known to contain $\gamma_{0}$ and $\hat{\gamma}$ and let $\delta\in\Gamma.$
The zero Hadamard derivative in Theorem 3 implies a zero Gateaux derivative
with respect to $\delta$ so that%
\begin{equation}
\bar{\psi}_{\gamma}(\delta,\alpha_{0},\theta_{0})=\frac{d\bar{\psi}(\gamma
_{0}+\tau\delta,\alpha_{0},\theta_{0})}{d\tau}=0. \label{orth 2}%
\end{equation}
This can be thought of as a population moment condition for $\alpha_{0}$. We
can form a sample moment function corresponding to this population moment
condition by replacing the expectation by a sample average, $\gamma_{0}$ by
$\hat{\gamma}_{\ell}$, and $\theta_{0}$ by $\hat{\theta}_{\ell}$ to obtain%
\begin{equation}
\hat{\psi}_{\gamma}(\delta,\alpha)=\frac{d}{d\tau}\frac{1}{n-n_{\ell}}%
\sum_{i\notin I_{\ell}}\psi(W_{i},\hat{\gamma}_{\ell}+\tau\delta,\alpha
,\tilde{\theta}_{\ell}). \label{orth sample}%
\end{equation}
We can then replace $\alpha$ by a sieve (i.e. parametric approximation) and
estimate the sieve parameters using these sample moments for a variety of
choices of $\delta.$ We can also regularize to allow for a high dimensional
specification for $\alpha.$ The sample moments in equation (\ref{orth sample})
depend only on observations not in $I_{\ell}$ so that the resulting
$\hat{\alpha}_{\ell}$ will also, as required for debiased GMM.

\bigskip

\textsc{Example 1:} Here $\alpha_{0}$ is an function of $X$ that has finite
second moment and
\begin{align*}
\hat{\psi}_{\gamma}(\delta,\alpha)  &  =\frac{d}{d\tau}\frac{1}{n-n_{\ell}%
}\sum_{i\notin I_{\ell}}\left\{  Z_{i}[\hat{\gamma}_{\ell}(X_{i})+\tau
\delta(X_{i})]+\alpha(X_{i})[Y_{i}-\hat{\gamma}_{\ell}(X_{i})-\tau\delta
(X_{i})]-\tilde{\theta}_{\ell}\right\} \\
&  =\frac{1}{n-n_{\ell}}\sum_{i\notin I_{\ell}}\left[  Z_{i}-\alpha
(X_{i})\right]  \delta(X_{i}).
\end{align*}
This is a sample moment corresponding to population moment condition
$E[\left\{  Z-\alpha_{0}(X)\right\}  \delta(X)]=0,$ which holds by $\alpha
_{0}(X)=E[Z|X].$ If $\alpha(X)$ was replaced by a linear combination
$\rho^{\prime}b(x)$ of a dictionary $b(x)=(b_{1}(x),...,b_{p}(x))^{\prime}$
and $\delta(X)$ replaced by element $b_{j}(X)$ then the sample moment
function
\[
\hat{\psi}_{\gamma}(b_{j},\rho^{\prime}b)=\frac{1}{n-n_{\ell}}\sum_{i\notin
I_{\ell}}\left[  Z_{i}-\rho^{\prime}b(X_{i})\right]  b_{j}(X_{i}).
\]
The collection of sample moments $\hat{\psi}_{\gamma}(b_{j},\rho^{\prime}b),$
$(j=1,...,p)$ can be used to construct a Lasso or Dantzig estimator of
$\alpha_{0}(X)=E[Z|X].$

\bigskip

\textsc{Example 4:} Consider an additional example where the parameter of
interest is a linear function $\theta_{0}=E[m(W,\gamma_{0})]$ of an unknown
function $\gamma_{0}$ with%
\begin{equation}
\gamma_{0}(X)=\arg\min_{\gamma\in\Gamma}E[v(Y-\gamma(X))],\text{ }%
E[m(W,\gamma)]=E[\bar{\alpha}(X)\gamma(X)],\text{ for all }E[\gamma
(X)^{2}]<\infty, \label{cond loc}%
\end{equation}
$v(u)$ is a convex function, $\Gamma$ is a set of functions that is closed in
mean square, and $\bar{\alpha}(X)$ is the Riesz representer for the linear
functional $E[m(X,\gamma)]$. Here $\gamma_{0}(X)$ measures the conditional
location of $Y$ given $X$. For example $\gamma_{0}(X)=E[Y|X]$ when
$v(u)=u^{2}/2$ and $\gamma_{0}(X)$ is the $\lambda^{th}$ conditional quantile
when $v(u)=[\lambda-1(u<0)]u.$ For identifying moment function $g(w,\gamma
,\theta)=m(w,\gamma)-\theta$ the nonparametric influence function is%
\[
\phi(w,\gamma,\alpha,\theta)=\alpha(x)v_{u}(Y-\gamma(X)),\text{ }\alpha
_{0}(X)=\frac{\bar{\alpha}(X)}{\bar{v}_{%
\operatorname{u}%
\operatorname{u}%
}(X)},\text{ }\bar{v}_{%
\operatorname{u}%
\operatorname{u}%
}(X)=\frac{d}{d\tau}E[v_{u}(Y-\gamma_{0}(X)+\tau)|X].
\]
This influence function appears in Ai and Chen (2007) for a series estimator
of $\gamma_{0}$ and Ichimura and Newey (2017) show it satisfies equation
(\ref{infdef}). Replacing $\alpha(X)$ by a linear combination $\rho^{\prime
}b(X)$ of a dictionary $b(X)$ and choosing $\delta=b_{j}$ gives the sample
moment%
\begin{align*}
\hat{\psi}_{\gamma}(b_{j},\rho^{\prime}b)  &  =\frac{d}{d\tau}\frac
{1}{n-n_{\ell}}\sum_{i\notin I_{\ell}}\left\{  m(W_{i},\hat{\gamma}_{\ell
}+\tau b_{j})+\rho^{\prime}b(X_{i})v_{u}(Y_{i}-\hat{\gamma}_{\ell}(X_{i})-\tau
b_{j}(X_{i}))-\tilde{\theta}_{\ell}\right\} \\
&  =\frac{1}{n-n_{\ell}}\sum_{i\notin I_{\ell}}\left[  m(W_{i},b_{j}%
)-\rho^{\prime}b(X_{i})v_{%
\operatorname{u}%
\operatorname{u}%
}(Y_{i}-\hat{\gamma}_{\ell}(X_{i}))b_{j}(X_{i})\right]  =\hat{M}_{\ell j}%
-\rho^{\prime}\hat{Q}_{\ell}e_{j},\\
\hat{M}_{\ell j}  &  =\frac{1}{n-n_{\ell}}\sum_{i\notin I_{\ell}}m(W_{i}%
,b_{j}),\text{ }\hat{Q}_{\ell}=\frac{1}{n-n_{\ell}}\sum_{i\notin I_{\ell}}v_{%
\operatorname{u}%
\operatorname{u}%
}(Y_{i}-\hat{\gamma}_{\ell}(X_{i}))b(X_{i})b(X_{i})^{\prime},
\end{align*}
where we have assumed that $v_{u}(u)$ has derivative $v_{%
\operatorname{u}%
\operatorname{u}%
}(u)$ and $e_{j}$ is the $j^{th}$ unit vector.

Lasso or Dantzig estimators of $\alpha_{0}(x)$ can be constructed from the
collection of sample moments $\hat{\psi}_{\gamma}(b_{j},\rho^{\prime}b),$
$(j=1,...,p).$ Let $\hat{M}_{\ell}=(\hat{M}_{\ell1},...,\hat{M}_{\ell
p})^{\prime}$. A Lasso estimator of $\alpha_{0}(x)$ is%
\begin{equation}
\hat{\alpha}_{L\ell}(x)=\hat{\rho}_{L\ell}^{\prime}b(x),\text{ }\hat{\rho
}_{L\ell}=\arg\min_{\rho}\left\{  -2\hat{M}_{\ell}^{\prime}\rho+\rho^{\prime
}\hat{Q}_{\ell}\rho+2r\sum_{j=1}^{p}\left\vert \rho_{j}\right\vert \right\}  .
\label{Lasso min dist}%
\end{equation}
This is a Lasso minimum distance estimator of $\alpha_{0}(x)$ that generalize
that of Chernozhukov, Newey, and Singh (2018) for $v(u)=u^{2}/2$ to any twice
differentiable convex function $v(u).$ It has the nice feature that an
explicit estimator of $\bar{v}_{%
\operatorname{u}%
\operatorname{u}%
}(X)$ is not required to be places in a denominator. Instead the presence of
$\bar{v}_{%
\operatorname{u}%
\operatorname{u}%
}(X)$ is accounted for in the weighted second moment estimator $\hat{Q}_{\ell
}.$

\bigskip

\textsc{Example 2:} When $\gamma_{0}(X)$ is a conditional quantile $v_{u}(u)$
is not differentiable so a different $\hat{Q}_{\ell}$ must be used. The $v_{%
\operatorname{u}%
\operatorname{u}%
}(Y_{i}-\hat{\gamma}_{\ell}(X_{i}))$ in $\hat{Q}_{\ell}$ must be replaced with
something such that $\hat{Q}_{\ell}$ estimates $Q=E[\bar{v}_{%
\operatorname{u}%
\operatorname{u}%
}(X)b(X)b(X)^{\prime}].$ For quantiles $\bar{v}_{%
\operatorname{u}%
\operatorname{u}%
}(X)=f(0|X)$ where $f(0|X)$ is the conditional pdf of $U=Y-\gamma_{0}(X)$
given $X$ at $U=0$. The $\hat{Q}_{\ell}$ given in Example 2 in Section 2 takes
account of this in the presence of the kernel term in%
\[
\hat{Q}_{\ell}=\frac{1}{n-n_{\ell}}\sum_{\ell^{\prime}\neq\ell}\sum_{i\in
I_{\ell^{\prime}}}\frac{1}{h}K(\frac{Y_{i}-\hat{\gamma}_{\ell,\ell^{\prime}%
}(X_{i})}{h})b(X_{i})b(X_{i})^{\prime},
\]
where $\hat{\gamma}_{\ell,\ell^{\prime}}(X_{i})$ uses only observations not in
$I_{\ell}$ or $I_{\ell}^{\prime}.$ Here $\hat{Q}_{\ell}$ estimates
$Q=E[f(0|X)b(X)b(X)^{\prime}]$. We also use nested sample splitting to
construct $\hat{\gamma}_{\ell,\ell^{\prime}}(X_{i})$ in order to obtain
asymptotic results for $\hat{\theta}$ in Section 8 using just a mean square
convergence rate for $\hat{\gamma}_{\ell\ell^{\prime}}$.

\bigskip

It would be interesting to use the moment functions (\ref{orth sample}) to
construct $\hat{\alpha}$ for first steps other than a conditional location
function $\gamma_{0}(X)$ in equation (\ref{cond loc}). That is beyond the
scope of this paper and is reserved to future work, including identification
of $\alpha_{0}$ and asymptotic theory for $\hat{\alpha}$.

The construction of estimating equations for $\alpha$ from the Gateaux
derivative of orthogonal sample moments with respect to variation in the first
step away from $\hat{\gamma}$ in equation (\ref{orth sample}) generalizes the
estimating equations of Chernozhukov, Newey, and Singh (2018) for linear
functionals of a conditional expectation to any orthogonal moment function and
first step $\hat{\gamma}$. These also generalize estimators for $\alpha$ for
average treatment effects that were proposed and analyzed in Vermeulen and
Vansteelandt (2015), Avagyan and Vansteelandt (2017), and Tan (2018).

This approach of estimating the nonparametric influence function uses its form
$\phi(w,\gamma,\alpha,\theta)$ to construct an estimator of $\alpha_{0}$.
Using the form of $\phi(w,\gamma,\alpha,\theta)$ seems good in high
dimensional settings where $\phi(w,\gamma,\alpha,\theta)$ may be a high
dimensional function. It is also possible to estimate the entire nonparametric
influence function using just the first step and the identifying moments. Such
estimators are available for first step series and kernel estimation. For
first step series estimation an estimator of $\phi(w,\gamma,\alpha,\theta)$
can be constructed by treating the first step estimator as if it were
parametric and applying a standard formula parametric two-step estimators,
e.g. as in Newey (1984) and Murphy and Topel (1985). Newey (1994a, 1997) and
Ackerberg, Chen, and Hahn (2012) used this approach for estimating the
asymptotic variance of functions of series estimators, while we refer here to
another use in constructing debiased GMM. Debiased GMM estimators can be
constructed by adding the nonparametric influence function obtained by
treating the first step as if it were parametric. For parametric maximum
likelihood the resulting orthogonal moment functions are the basis of Neyman's
(1959) C-alpha test. Wooldridge (1991) generalized such moment conditions to
parametric nonlinear least squares and Lee (2005), Bera, Motes-Rohas, and
Sosa-Escudero. (2010), and Chernozhukov, Hansen, and Spindler to GMM.

For first step kernel estimation one can use the numerical influence function
estimator of Newey (1994b) to estimate $\phi(w,\gamma,\alpha,\theta),$ as
suggested in a previous version of this paper and proven to work in a low
dimensional nonparametric setting in Bravo, Escanciano, and van Keilegom
(2020). The idea is to differentiate with respect to the effect of the
$i^{th}$ observation on sample moments. Kernel estimators are not well suited
to high dimensional settings with machine learning so we do not consider them here.

It is also possible to estimate the nonparametric influence function using a
numerical derivative version of equation (\ref{infdef}). This approach has
been given in Carone, Luedtke, and van der Laan (2016) and Bravo, Escanciano,
and van Keilegom (2020) for construction of orthogonal moment functions. We
focus here on estimating $\alpha_{0}$ where there is a known form
$\phi(w,\gamma,\alpha,\theta)$ because that information is widely available
and using it seems likely to be good in high dimensional settings.

\section{Double Robustness}

The zero derivative condition in equation (\ref{lrdef}) is an appealing
robustness property. This condition can be interpreted as local insensitivity
of the moments to the value of $\gamma$, with the moments remaining close to
zero as $\gamma$ varies away from its true value. Because it is difficult to
get nonparametric functions exactly right, especially in high dimensional
settings, this property is an appealing one.

Such robustness considerations, well explained in Robins and Rotnitzky (2001),
have motivated the development of doubly robust moment conditions. Doubly
robust moment conditions have expectation zero if one first step component is
incorrect. Doubly robust moment conditions allow two chances for the moment
conditions to hold, an appealing robustness feature. Also, doubly robust
moment conditions have simpler conditions for asymptotic normality than
general debiased GMM, as discussed in Section 8.

In this Section we characterize double robustness and derive several novel
classes of doubly robust moment conditions. We construct doubly robust moment
functions by adding to identifying moment functions the nonparametric
influence functions to obtain orthogonal moments. In this way the derivation
of new doubly robust moment functions is aided by the construction of
orthogonal moment functions from adding the nonparametric influence function.

\subsection{Characterizing Double Robustness}

Double robustness is that for all $\gamma\in\Gamma,$ $\alpha$, and $\theta$%
\[
0=\bar{\psi}(\gamma,\alpha_{0},\theta_{0})=\bar{\psi}(\gamma_{0},\alpha
,\theta),
\]
where $\gamma\in\Gamma$ is a set of possible first steps. The second equality
already follows from Theorem 2. The first conclusion of Theorem 3 gives a
local version of the first equality. If $\bar{\psi}(\gamma,\alpha_{0}%
,\theta_{0})$ is affine in $\gamma$ then this local property becomes global so
that double robustness holds. Clearly doubly robust moment conditions are also
affine, so that $\bar{\psi}(\gamma,\alpha_{0},\theta_{0})$ being affine is a
complete characterization of double robustness. The following result gives
this characterization for a zero Gateaux derivative, a condition easier to
check than zero Hadamard derivative in Theorem 3.

\bigskip

\textsc{Theorem 7:}\textit{ If }$\Gamma$\textit{ is linear then}
$\psi(w,\gamma,\alpha,\theta)$\textit{ is doubly robust if and only if for
every }$\gamma\in\Gamma$\textit{ }%
\[
\left.  \frac{\partial\bar{\psi}((1-\tau)\gamma_{0}+\tau\gamma,\alpha
_{0},\theta_{0})}{\partial\tau}\right\vert _{\tau=0}=0,
\]
\textit{ and }$\bar{\psi}(\gamma,\alpha_{0},\theta_{0})$\textit{ is affine in
}$\gamma$\textit{.}

\bigskip

This characterization can be used to construct doubly robust moment conditions
from identifying moment conditions that are affine in $\gamma$. If
$g(W,\gamma,\theta_{0})$ is affine in $\gamma$ and $\phi(W,\gamma,\alpha
_{0},\theta_{0})$ is also affine in $\gamma$ then $\bar{\psi}(\gamma
,\alpha_{0},\theta_{0})$ will also be affine in $\gamma.$ In addition the zero
Gateaux derivative condition in Theorem 6 will hold by Theorem 3, so that
$\bar{\psi}(\gamma,\alpha,\theta)$ will be doubly robust. We use this
construction to obtain doubly robust moment conditions for first steps that
satisfy conditional moment restrictions and for a first step that is a density function.

Robins and Rotnitzky (2001) gave conditions for the existence of doubly robust
moment conditions in semiparametric models. Theorem 5 is complementary to
those results in giving a complete characterization of doubly robust moments
when $\Gamma$ is linear.

\subsection{Double Robustness with First Step Conditional Moment Restriction}

A novel class of doubly robust moment functions are those where the first step
$\gamma_{0}$ satisfies a conditional moment restriction
\begin{equation}
E[\lambda(W,\gamma_{0})|X]=0, \label{cond mom restrict}%
\end{equation}
where $\lambda(W,\gamma)$ is a scalar functional of $\gamma$ that is affine in
$\gamma$ and $X$ are regressors or instrumental variables. Suppose that
$\hat{\gamma}$ is the nonparametric 2SLS estimator of Newey and Powell (1989,
2003) and Newey (1991) where $\gamma(F)=\arg\min_{\Gamma}E_{F}[\{E_{F}%
[\lambda(W,\gamma)|X]\}^{2}].$ It follows from Newey (1994a), Ai and Chen
(2007, p. 40), and Ichimura and Newey (2017) that when equation
(\ref{cond mom restrict}) is satisfied and $\phi(w,\gamma,\alpha,\theta)$
exists there is $\alpha(x,\theta)$ such that%
\[
\phi(w,\gamma,\alpha,\theta)=\alpha(x,\theta)\lambda(W,\gamma).
\]
The following result characterizes double robustness in this setting.

\bigskip

\textsc{Theorem 8:}\textit{\ }$\psi(W,\gamma,\alpha,\theta)=g(W,\gamma
,\theta)+\alpha(X)\lambda(W,\gamma)$\textit{ is doubly robust over }$\gamma
\in\Gamma$ \textit{if and only if}%
\[
E[g(W,\gamma,\theta_{0})]=-E[\alpha_{0}(X)\lambda(W,\gamma)]\text{,
\textit{for all} }\gamma\in\Gamma.
\]

\bigskip

When $E[g(W,\gamma,\theta_{0})]$ is affine in $\gamma$ the condition in
Theorem 8 is an expected outer product representation of $E[g(W,\gamma
,\theta_{0})].$ This characterization of a doubly robust moment function with
a nonparametric 2SLS first step has several interesting special cases. The
form of the doubly robust moment function in Theorem 8 is similar to the
efficient score in some cases in Ai and Chen (2012).

\bigskip

\textsc{Example 5: }An important example is a linear functional $\theta
_{0}=E[m(W,\gamma_{0})]$ of a regression function $\gamma$, where
$\lambda(W,\gamma)=Y-\gamma(X)$ for some outcome variable $Y$ and
$m(w,\gamma)$ is linear in $\gamma$. Here the identifying moment function is
$g(w,\gamma,\theta)=m(w,\gamma)-\theta$, which is affine in $\gamma$. Also
$\lambda(w,\gamma)$ is also affine in $\gamma$, so that the conditions of
Theorem 8 are satisfied. The next result follows from Theorem 8.

\bigskip

\textsc{Corollary 9: }\textit{If }$m(w,\gamma)$\textit{ is linear in }$\gamma
$\textit{ and there is }$\alpha_{0}(x)$ \textit{such that }$E[\alpha
_{0}(X)^{2}]<\infty$\textit{ and }$E[m(W,\gamma)]=E[\alpha_{0}(X)\gamma
(X)]$\textit{ for all }$E[\gamma(X)^{2}]<\infty$\textit{ then }$\psi
(w,\gamma,\alpha,\theta)=m(w,\gamma)-\theta+\alpha(x)[y-\gamma(x)]$\textit{ is
doubly robust.}

\bigskip

In Corollary 9 $\alpha_{0}(x)$ is the Riesz representer of the functional
$E[m(W,\gamma)]=E[\alpha_{0}(X)\gamma(X)]$ for all $\gamma$ with
$E[\gamma(X)^{2}]<\infty$, as in Proposition 4 of Newey (1994a). Many
important doubly robust moment functions are special cases of Corollary 9,
including average treatment effects, policy effects, and average derivatives,
as discussed in Newey and Robins (2017), Chernozhukov, Newey, and Robins
(2018), Hirshberg and Wager (2018), and Chernozhukov, Newey, and Singh (2018).
In these papers Corollary 9 is also used to derive new doubly robust moment functions.

\bigskip

\textsc{Example 6:} An interesting generalization allowing for endogeneity has
$\lambda(W,\gamma)=Y-\gamma(Z)$ where $Z$ need not be equal to $X$. Here the
conditional moment restriction (\ref{cond mom restrict}) is a nonparametric
instrumental variables model as in Newey and Powell (1989, 2003) and Newey
(1991). Theorem 8 can be applied to derive doubly robust moment functions
since $\lambda(W,\gamma)=Y-\gamma(Z)$ is affine in $\gamma.$

\bigskip

\textsc{Corollary 10: }\textit{If }$m(w,\gamma)$\textit{ is linear in }%
$\gamma$\textit{ and there is }$\alpha_{0}(x)$ \textit{such that }%
$E[\alpha_{0}(X)^{2}]<\infty$\textit{ and }$E[m(W,\gamma)]=E[\alpha
_{0}(X)\gamma(Z)]$\textit{ for all }$E[\gamma(Z)^{2}]<\infty$\textit{ then
}$\psi(w,\gamma,\alpha,\theta)=m(w,\gamma)-\theta+\alpha(x)[y-\gamma
(z)]$\textit{ is doubly robust.}

\bigskip

As discussed in Ichimura and Newey (2017), if there is $v(Z)$ with
$E[v(Z)^{2}]<\infty$ and $E[m(W,\gamma)]=E[v(Z)\gamma(Z)]$ then existence of
$\alpha_{0}(X)$ satisfying the condition of Corollary 10 requires
$v(Z)=E[\alpha_{0}(X)|Z]$, which is necessary for root-n consistent
estimability of $\theta_{0}$, as shown by Severini and Tripathi (2012).

\bigskip

\textsc{Example 7:} Many novel examples of doubly robust moment functions can
be derived from Corollary 10, including policy effects and average
derivatives. A weighted average derivative example has $m(w,\gamma)=\bar
{v}(z)\partial\gamma(z)/\partial z_{1}$ for some known $\bar{v}(z).$ A doubly
robust moment function is%
\[
\psi(w,\gamma,\alpha,\theta)=\bar{v}(z)\frac{\partial\gamma(z)}{\partial
z_{1}}-\theta+\alpha(x)[y-\gamma(z)],\text{ }E[\alpha_{0}(X)|Z]=-\frac
{\partial\{f_{0}(Z)\bar{v}(Z)\}/\partial z_{1}}{f_{0}(Z)},
\]
where $f_{0}(z)$ is the marginal pdf of $Z$. This is a doubly robust moment
function that could be used to construct a doubly robust version of the
plug-in estimator of Ai and Chen (2007).

\bigskip

Using Theorem 7 to construct doubly robust moment functions can depend on
specifying $\gamma$ to make $g(W,\gamma,\theta_{0})$ and $\phi(W,\gamma
,\alpha_{0},\theta_{0})$ affine in $\gamma.$ We illustrate with a well known example.

\bigskip

\textsc{Example 8: }Suppose that the object of interest is $\theta
_{0}=E[Y^{\ast}]$ where $Y=1(D=1)Y^{\ast}$ is observed for a observed
completed data indicator $D\in\{0,1\}$ and the data are missing at random with
$E[Y^{\ast}|X,D=1]=E[Y^{\ast}|X]$ for observed covariates $X.$ Inverse
probability weighting gives $\theta_{0}=E[DY/P_{0}(X)]=E[P_{0}(X)^{-1}%
E[DY|X]],$ which is nonlinear in the unknown propensity score $P_{0}%
(X)=\Pr(D=1|X)$. A corresponding affine in $\gamma$ identifying moment
function is $g(w,\gamma,\theta)=g(w,\gamma,\theta)=\gamma(x)dy-\theta$ with
true first step $\gamma_{0}(X)=P_{0}(X)^{-1}.$ This $\gamma_{0}$ satisfies the
conditional moment restriction in equation (\ref{cond mom restrict}) for
$\lambda(w,\gamma)=1-\gamma(x)d$ that is affine in $\gamma.$ Also, for
$\alpha_{0}(X)=E[Y|X,D=1]=E[DY|X]\gamma_{0}(X)$ we have%
\begin{align*}
E[g(W,\gamma,\theta_{0})]  &  =E[E[DY|X]\{\gamma(X)-\gamma_{0}(X)\}]=E[\alpha
_{0}(X)\gamma_{0}(X)^{-1}\{\gamma(X)-\gamma_{0}(X)\}]\\
&  =E[\alpha_{0}(X)\{\gamma(X)P_{0}(X)-1\}]=E[\alpha_{0}(X)\{\gamma
(X)D-1\}]=-E[\alpha_{0}(X)\lambda(W,\gamma)].
\end{align*}
The doubly robust moment function from Theorem 8 is then $\psi(w,\gamma
,\alpha,\theta)=\gamma(x)dy-\theta+\alpha(x)(1-\gamma(x)d),$ which is the
doubly robust moment function of Robins, Rotnitzky, and Zhao (1994). This
example shows how that classic doubly robust moment function is a special case
of Theorem 7, with moment condition that is affine in a first step $\gamma$
for $\gamma_{0}$ equal to the inverse propensity score. The only if part of
Theorem 7 states that every doubly robust moment function will have
expectation that is affine in $\gamma$.

Rotnitzky, Smucler, and Robins (2019) give a general class of robust
estimators that includes interesting examples not treated here.

\subsection{Double Robustness with First Step Probability Density}

Another novel class of doubly moment conditions are those where the first step
$\gamma$ is a pdf of a function $X$ of the data observation $W.$ By
Proposition 5 of Newey (1994a), the first step influence function is
\[
\phi(w,\gamma,\alpha,\theta)=\alpha(x)-\int\alpha(u)\gamma(u)du,
\]
which is affine in $\gamma$. When the identifying moment function is affine
adding this nonparametric influence function gives a doubly robust moment function.

\bigskip

\textsc{Theorem 11:}\ \textit{If there exists }$\alpha_{0}(x)$ \textit{with
}$\int\alpha_{0}(u)^{2}du$ $<\infty$\textit{ and }$E[g(W,\gamma,\theta
_{0})]=\int\alpha_{0}(u)[\gamma(u)-\gamma_{0}(u)]du$ \textit{for all }$\gamma$
with $\int\gamma(u)^{2}du<\infty$\textit{ then }$\psi(W,\gamma,\alpha
,\theta)=g(W,\gamma,\theta)+\alpha(X)-\int\alpha(u)\gamma(u)du$\textit{ is
doubly robust.}

\bigskip

Here $\alpha_{0}(x)$ is the Riesz representer of Proposition 5 of Newey
(1994a) for the Lebesgue inner product

\bigskip

\textsc{Example 9: }An example is the density weighted average derivative of
Powell, Stock, and Stoker (1989), where $g(w,\gamma,\theta)=-2y\cdot
\partial\gamma(x)/\partial x-\theta$ and $\alpha_{0}(x)=\partial
\{E[Y|X=x]\gamma_{0}(x)\}/\partial x.$ Because $g(w,\gamma,\theta)$ is affine
in $\gamma$ Theorem 11 implies%
\[
\psi(W,\gamma,\alpha,\theta)=-2Y\frac{\partial\gamma(X)}{\partial x}%
-\theta+\alpha(X)-\int\alpha(u)\gamma(u)du,
\]
is doubly robust. Double robustness of this moment function seems to be a
novel result.

\subsection{Identification Via Doubly Robust Moment Conditions}

Doubly robust moment conditions can be used to identify parameters of interest.

\bigskip

\textsc{Theorem 12: }\textit{If Assumption 1 is satisfied, }$\alpha_{0}%
$\textit{ is identified, and for some }$\bar{\gamma}$\textit{ the equation
}$E[\psi(W,\bar{\gamma},\alpha_{0},\theta)]=0$\textit{ has a unique solution
at }$\theta=$\textit{ }$\theta_{0}$\textit{ then }$\theta_{0}$\textit{ is
identified as that solution.}

\bigskip

\textsc{Example 10:} Applying this result to the nonparametric instrumental
variables setting of Assumption 6 leads to identification of functionals of
$\gamma_{0}$ without requiring that $\gamma_{0}$ be identified. Focusing on
Example 7, note that $\alpha_{0}(X)$ is identified as a solution to
$-f_{0}(Z)^{-1}\partial\{f_{0}(Z)\bar{v}(Z)\}/\partial z_{1}=E[\alpha
_{0}(X)|Z].$ Setting $\bar{\gamma}=0$ in Theorem 7 then identifies $\theta
_{0}=E[\alpha_{0}(X)Y]$, extending Santos (2011) and Severini and Tripathi
(2006, 2012), to the weighted average derivative.

\subsection{Partial Robustness of Plug-In GMM}

Partial robustness refers to identifying moments where $E[g(W,\theta_{0}%
,\bar{\gamma})]=0$ for some $\bar{\gamma}\neq\gamma_{0}$. This is a weaker
property for identifying moment function than double robustness for the
associated orthogonal moment function. Also, partial robustness for
identifying moments that are affine in $\gamma$ with $E[\phi(W,\gamma
,\alpha_{0},\theta_{0})]$ affine in $\gamma$ can be characterized by the
nonparametric influence function, since double robustness implies%
\[
E[g(W,\theta_{0},\gamma)]=-E[\phi(W,\gamma,\alpha_{0},\theta_{0})].
\]
We give two examples of partial robustness results that follow from double robustness.

\bigskip

\textsc{Example 5: } For a linear functional $\theta_{0}=E[m(W,\gamma_{0})]$
of a regression function $\gamma_{0}(X)=E[Y|X]$, let $b(X)$ be a $p\times1$
vector of functions of $X$ and $\bar{\gamma}(X)=$ $b(X)^{\prime}\delta,$
$\delta=E[b(X)b(X)^{\prime}])^{-1}E[b(X)Y]$, be the best linear predictor of
$\gamma_{0}(X)$ by $b(X)$.

\bigskip

\textsc{Theorem 13:} \textit{If }$E[b(X)b(X)^{\prime}]$\textit{ is nonsingular
and }$\alpha_{0}(X)=\rho_{0}^{\prime}b(X)$\textit{ for some }$\rho_{0}%
$\textit{ then }$\theta_{0}=E[m(W,\bar{\gamma})]$\textit{. }

\bigskip

This result generalizes Stoker's (1986) result that linear regression
coefficients equal average derivatives when the regressors are multivariate
Gaussian to any linear functional $m(w,\gamma)$ and nonlinear $b(X)$.

\bigskip

\textsc{Example 7:} Consider the average derivative $\theta_{0}=E[\partial
\gamma_{0}(Z)/\partial z_{1}]$ where $g(w,\gamma,\theta)=\partial
\gamma(z)/\partial z_{1}-\theta.$ Let $\delta=(E[a(X)b(Z)^{\prime}%
])^{-1}E[a(X)Y]$ be the limit of the linear instrumental variables estimator
with right hand side variables $b(Z)$ and the same number of instruments
$a(X)$, and $\bar{\gamma}(Z)=b(Z)^{\prime}\delta$ the linear instrumental
variables estimand.

\bigskip

\textsc{Theorem 14:} \textit{If }$-\partial\ln f_{0}(Z)/\partial
z_{r}=c^{\prime}b(Z)$\textit{ for a constant vector }$c$\textit{,
}$E[b(Z)b(Z)^{\prime}]$\textit{ is nonsingular, and }$E[a(X)|Z]=\Pi
b(Z)$\textit{ for a square nonsingular }$\Pi$\textit{ then }$\theta
_{0}=E[\partial\bar{\gamma}(Z)/\partial z_{1}].$

\bigskip

This is a generalization to nonparametric instrumental variables of Stoker's
(1986) result.

\section{Asymptotic Theory}

In this Section we give simple and general asymptotic theory for debiased GMM.
The results differ from Chernozhukov et al. (2018) in the use of Theorem 2, a
different remainder decomposition that leads to simpler conditions, and
incorporation of double robustness in the general conditions. We begin with
conditions for the key property%
\begin{equation}
\sqrt{n}\hat{\psi}(\theta_{0})=\frac{1}{\sqrt{n}}\sum_{i=1}^{n}\psi
(W_{i},\theta_{0},\gamma_{0},\alpha_{0})+o_{p}(1). \label{no effec}%
\end{equation}

\textsc{Assumption 1: }$E[\left\Vert \psi(W_{i},\theta_{0},\gamma_{0}%
,\alpha_{0})\right\Vert ^{2}]<\infty$\textit{ and }%
\begin{align*}
&  \mathit{i)\ }\int\left\Vert g(w,\hat{\gamma}_{\ell},\theta_{0}%
)-g(w,\gamma_{0},\theta_{0})\right\Vert ^{2}F_{0}%
(dw)\overset{p}{\longrightarrow}0;\text{ }\mathit{ii)\ }\int\left\Vert
\phi(w,\hat{\gamma}_{\ell},\alpha_{0},\theta_{0})-\phi(w,\gamma_{0},\alpha
_{0},\theta_{0})\right\Vert ^{2}F_{0}(dw)\overset{p}{\longrightarrow}0,\\
&  \mathit{iii)\ }\int\left\Vert \phi(w,\gamma_{0},\hat{\alpha}_{\ell}%
,\tilde{\theta}_{\ell})-\phi(w,\gamma_{0},\alpha_{0},\theta_{0})\right\Vert
^{2}F_{0}(dw)\overset{p}{\longrightarrow}0.
\end{align*}

\bigskip

These are mild mean square consistency conditions for $\hat{\gamma}_{\ell}$
and $(\hat{\alpha}_{\ell},\tilde{\theta}_{\ell})\,$separately. They differ
from Chernozhukov et al. (2018) in being separate conditions for $\hat{\gamma
}_{\ell}$ and $(\hat{\alpha}_{\ell},\tilde{\theta}_{\ell})$ and for
$g(w,\gamma,\theta)$ and $\phi(w,\gamma,\alpha,\theta).$ Let
\[
\hat{\Delta}_{\ell}(w)=\phi(w,\hat{\gamma}_{\ell},\hat{\alpha}_{\ell}%
,\tilde{\theta}_{\ell})-\phi(w,\gamma_{0},\hat{\alpha}_{\ell},\tilde{\theta
}_{\ell})-\phi(w,\hat{\gamma}_{\ell},\alpha_{0},\theta_{0})+\phi(w,\gamma
_{0},\alpha_{0},\theta_{0})
\]

\bigskip

\textsc{Assumption 2: }\textit{For each }$\ell=1,...,L$\textit{, either i)}%
\[
\sqrt{n}\int\hat{\Delta}_{\ell}(w)F_{0}(dw)\overset{p}{\longrightarrow
}0,\text{ }\int\left\Vert \hat{\Delta}_{\ell}(w)\right\Vert ^{2}%
F_{0}(dw)\overset{p}{\longrightarrow}0,
\]
\textit{ or ii) }$\sum_{i\in I_{\ell}}\left\Vert \hat{\Delta}_{\ell}%
(W_{i})\right\Vert /\sqrt{n}\overset{p}{\longrightarrow}0,$ \textit{or iii)
}$\sum_{i\in I_{\ell}}\hat{\Delta}_{\ell}(W_{i})/\sqrt{n}%
\overset{p}{\longrightarrow}0.$

\bigskip

This condition imposes a rate condition on the interaction remainder
$\hat{\Delta}_{\ell}(w)$, that its average must go to zero faster than
$1/\sqrt{n}.$ It differs from Chernozhukov et al. (2018) in applying only to
the nonparametric influence function and allowing for the sample average rate
condition in iii), which is helpful for obtaining weak regularity conditions
of Newey and Robins (2017).

\bigskip

\textsc{Assumption 3:} \textit{For each }$\ell=1,...,L$\textit{, i) }$\int%
\phi(w,\gamma_{0},\hat{\alpha}_{\ell},\tilde{\theta}_{\ell})F_{0}(dw)=0;$
\textit{and either ii) }$\bar{\psi}(\gamma,\alpha_{0},\theta_{0})$ \textit{is
affine in }$\gamma;$ \textit{or iii) }$\left\Vert \hat{\gamma}_{\ell}%
-\gamma_{0}\right\Vert =o_{p}(n^{-1/4})$ \textit{and }$\left\Vert \bar{\psi
}(\gamma,\alpha_{0},\theta_{0})\right\Vert \leq C\left\Vert \gamma-\gamma
_{0}\right\Vert ^{2}$\textit{ for all }$\gamma$ \textit{with }$\left\Vert
\gamma-\gamma_{0}\right\Vert $\textit{ small enough; or} iv)$\sqrt{n}\bar
{\psi}(\hat{\gamma}_{\ell},\alpha_{0},\theta_{0})\overset{p}{\longrightarrow
}0.$

\bigskip

Assumption 3 incorporates Theorem 2 in i) and doubly robust moment functions
through ii), in which case Assumption 3 imposes no conditions additional to
Assumptions 1 and 2. Conditions iii) and iv) are alternative small bias
conditions that are only required to hold for $\hat{\gamma}_{\ell},$ and not
for $\hat{\alpha}_{\ell}.$ Condition iii) requires a faster than $n^{-1/4}$
rate for $\hat{\gamma}$ as is familiar from the semiparametric estimation
literature. In many cases iii) will be satisfied for a mean square norm
$\left\Vert \cdot\right\Vert $ so that Assumptions 1-3 will only require
mean-square convergence rates, as is important in many machine learning
contexts where only mean square rates are available.

\bigskip

\textsc{Lemma 15:} \textit{If Assumptions 1-3 are satisfied then equation
(\ref{no effec}) is satisfied.}

\bigskip

This key asymptotic result differs from previous results of Andrews (1994),
Newey (1994a), Newey and McFadden (1994), Pakes and Olley (1995), Chen,
Linton, and van Keilegom (2003), Ichimura and Lee (2010), Escanciano et al.
(2016), and others in requiring no Donsker conditions. This feature is made
possible by the use of cross-fitting in the moment conditions. It is important
to not impose Donsker conditions for machine learning first steps which
generally do not, or are not known to, satisfy Donsker conditions, as
previously discussed in Chernozhukov et al. (2018).

This result improves upon Chernozhukov et al. (2018) in allowing $\hat{\alpha
}_{\ell}$ to converge slower than $n^{-1/4}$ in general, in Assumption 1
applying separately to $\hat{\gamma}_{\ell}$ and $\hat{\alpha}_{\ell}$, and
having weaker conditions for terms that involve both $\hat{\gamma}_{\ell}$ and
$\hat{\alpha}_{\ell}$ in Assumption 2. These improvements result from Theorem
2 and the structure of orthogonal moments as the sum of identifying moment
functions and the nonparametric influence function.

With additional conditions we obtain consistency of the estimator $\hat{\Psi}$
of the variance of the orthogonal moment functions given in Section 2. Let
$\Psi:=E[\psi(W,\gamma_{0},\alpha_{0},\theta_{0})\psi(W,\gamma_{0},\alpha
_{0},\theta_{0})^{\prime}]$.

\bigskip

\textsc{Lemma 16:} \textit{If Assumption 1 is satisfied and }$\int\left\Vert
g(w,\hat{\gamma},\tilde{\theta}_{\ell})-g(w,\hat{\gamma},\theta_{0}%
)\right\Vert ^{2}F_{0}(dw)\overset{p}{\longrightarrow}0$\textit{ and }%
$\int\left\Vert \hat{\Delta}_{\ell}(w)\right\Vert ^{2}F_{0}%
(dw)\overset{p}{\longrightarrow}0$\textit{ for each }$(\ell=1,...,L),$
\textit{then }$\hat{\Psi}\overset{p}{\longrightarrow}\Psi.$

\bigskip

It is also important to have conditions for convergence of the Jacobian of the
identifying sample moments $\partial\hat{g}(\bar{\theta})/\partial
\theta\overset{p}{\longrightarrow}G=E[\partial g(W,\gamma_{0},\theta
_{0})/\partial\theta]$ for any $\bar{\theta}\overset{p}{\longrightarrow}%
\theta_{0}$. To that end we impose the following condition:

\bigskip

\textsc{Assumption 4: }$G$\textit{\ exists and there is a neighborhood
}$\mathcal{N}$\textit{ of }$\theta_{0}$\textit{ and }$\left\Vert
\cdot\right\Vert $ \textit{such that i) for each }$\ell,$ $\left\Vert
\hat{\gamma}_{\ell}-\gamma_{0}\right\Vert \overset{p}{\longrightarrow}0;$ ii)
for all $\left\Vert \gamma-\gamma_{0}\right\Vert $ small enough $g(W,\gamma
,\theta)$\textit{ is differentiable in }$\theta$ \textit{on }$\mathcal{N}%
$\textit{ with probability approaching }$1$\textit{ and there is }%
$C>0$\textit{ and }$d(W,\gamma)$\textit{ such that for }$\theta\in\mathcal{N}%
$\textit{ and }$\left\Vert \gamma-\gamma_{0}\right\Vert $\textit{ small enough
}%
\[
\left\Vert \frac{\partial g(W,\gamma,\theta)}{\partial\theta}-\frac{\partial
g(W,\gamma,\theta_{0})}{\partial\theta}\right\Vert \leq d(W,\gamma)\left\Vert
\theta-\theta_{0}\right\Vert ^{1/C};\text{ }E[d(W,\gamma)]<C.
\]
\textit{iii) For each }$\ell=1,...,L,$ $j,$ and $k$, $\int\left\vert \partial
g_{j}(w,\hat{\gamma}_{\ell},\theta_{0})/\partial\theta_{k}-\partial
g_{j}(w,\gamma_{0},\theta_{0})/\partial\theta_{k}\right\vert F_{0}%
(dw)\overset{p}{\longrightarrow}0.$

\bigskip

\textsc{Lemma 17:} \textit{If Assumption 4 is satisfied and }$\bar{\theta
}\overset{p}{\longrightarrow}\theta_{0}$\textit{ then }$\partial\hat{g}%
(\bar{\theta})/\partial\theta\overset{p}{\longrightarrow}G.$

\bigskip

With these results in place the asymptotic normality of semiparametric GMM
follows in a standard way.

\bigskip

\textsc{Theorem 18:}\textit{ If Assumptions 1-4 are satisfied, }$\hat{\theta
}\overset{p}{\longrightarrow}\theta_{0},$\textit{ }$\hat{\Upsilon
}\overset{p}{\longrightarrow}\Upsilon$\textit{, and }$G^{\prime}\Upsilon
G$\textit{ is nonsingular, then}%
\[
\sqrt{n}(\hat{\theta}-\theta_{0})\overset{d}{\longrightarrow}N(0,V),\text{
}V=(G^{\prime}\Upsilon G)^{-1}G^{\prime}\Upsilon\Psi \Upsilon G(G^{\prime
}\Upsilon G)^{-1}.
\]
\textit{If also the conditions of Lemma 16 are satisfied then }$\hat{V}%
=(\hat{G}^{\prime}\hat{\Upsilon}\hat{G})^{-1}\hat{G}^{\prime}\hat{\Upsilon
}\hat{\Psi}\hat{\Upsilon}\hat{G}(\hat{G}^{\prime}\hat{\Upsilon}\hat{G}%
)^{-1}\overset{p}{\longrightarrow}V.$

\bigskip

This result and the Lemmas 15-17 are both general and simple. They are general
in applying to any first step estimators $\hat{\gamma}_{\ell}$ and
$\hat{\alpha}_{\ell}.$ They are simple in requiring only a few mean-square
convergence conditions for $\hat{\alpha}$ and $\hat{\gamma}$ when $\left\Vert
\cdot\right\Vert $ denotes the mean square norm. This generality and
simplicity results from the use of orthogonal moment functions and
cross-fitting. The orthogonality of the moment functions leads to Assumption
3, as shown by Theorems 2 and 3. The cross-fitting and a specific remainder
decomposition used in the proof of Lemma 15 motivate Assumptions 1 and 2, with
separate treatment of the identifying moment functions and the nonparametric
influence function.

\subsection{Functionals for a Conditional Moment Restriction}

Functionals of a first step satisfying a conditional moment restriction as in
Section 7.2 are of wide interest, including for Example 2. For an identifying
moment function $m(w,\gamma)-\theta$ and nonparametric influence function
$\alpha(x,\theta)\lambda(W,\gamma)$ debiased GMM $\hat{\theta}$ and $\hat{V}$
are%
\begin{align*}
\hat{\theta}  &  =\frac{1}{n}\sum_{\ell=1}^{L}\sum_{i\in I_{\ell}}\left[
m(W_{i},\hat{\gamma}_{\ell})+\hat{\alpha}_{\ell}(X_{i})\lambda(W_{i}%
,\hat{\gamma}_{\ell})\right]  ,\\
\hat{V}  &  =\frac{1}{n}\sum_{\ell=1}^{L}\sum_{i\in I_{\ell}}\hat{\psi}%
_{i\ell}^{2},\text{ }\hat{\psi}_{i\ell}=m(W_{i},\hat{\gamma}_{\ell}%
)+\hat{\alpha}_{\ell}(X_{i})\lambda(W_{i},\hat{\gamma}_{\ell})-\hat{\theta}.
\end{align*}
Let $\bar{\lambda}(X,\gamma)=E[\lambda(W,\gamma)|X]$, $\bar{\psi}(\gamma
)=\int[m(w,\gamma)+\alpha_{0}(x)\bar{\lambda}(x,\gamma)]F_{0}(dw)-$
$\theta_{0}$, $V=Var(m(W,\gamma_{0})+\alpha_{0}(X)\lambda(W,\gamma_{0}))$, and
$\left\Vert a\right\Vert =\sqrt{\int a(w)^{2}F_{0}^{{}}(w)}$ denote the mean
square norm.

\bigskip

\textsc{Theorem 19 }\textit{If i) }$E[\lambda(W,\gamma_{0})|X]=0$\textit{; ii)
}$\alpha_{0}(X)$\textit{ and }$E[\lambda(W,\gamma_{0})^{2}|X]$\textit{ are
bounded and }$E[m(W,\gamma_{0})^{2}]<\infty$\textit{; for }$(\ell
=1,...,L),$\textit{ iii) }$\int[m(w,\hat{\gamma})-m(w,\gamma_{0})]^{2}%
F_{0}(dw)\overset{p}{\longrightarrow}0$\textit{, }$\int[\lambda(w,\hat{\gamma
}_{\ell})-\lambda(w,\gamma_{0})]^{2}F_{0}(dw)\overset{p}{\longrightarrow}%
0,$\textit{ }$\left\Vert \hat{\alpha}_{\ell}-\alpha_{0}\right\Vert
\overset{p}{\longrightarrow}0;$\textit{ iv) }$\sqrt{n}\bar{\psi}(\hat{\gamma
}_{\ell})\overset{p}{\longrightarrow}0$; \textit{v) Either a) }$\int%
[\hat{\alpha}_{\ell}(x)-\alpha_{0}(x)]^{2}[\lambda(w,\hat{\gamma}_{\ell
})-\lambda(w,\gamma_{0})]^{2}\overset{p}{\longrightarrow}0$\textit{ and
}$\sqrt{n}\left\Vert \hat{\alpha}_{\ell}-\alpha_{0}\right\Vert \left\Vert
\bar{\lambda}(\hat{\gamma}_{\ell})-\bar{\lambda}(\gamma_{0})\right\Vert
\overset{p}{\longrightarrow}0,$ \textit{or b)} $\sqrt{n}\left\Vert \hat
{\alpha}_{\ell}-\alpha_{0}\right\Vert \left\Vert \lambda(\hat{\gamma}_{\ell
})-\lambda(\gamma_{0})\right\Vert \overset{p}{\longrightarrow}0$ \textit{and
}$\hat{\alpha}_{\ell}(x)$\textit{ in }$\hat{\psi}_{i\ell}$\textit{ is replaced
by }$\bar{\alpha}_{\ell}(x)=\hat{\alpha}_{\ell}(x)1(\left\vert \hat{\alpha
}_{\ell}(x)\right\vert \leq M)+sgn(\hat{\alpha}_{\ell}(x))M1(\left\vert
\hat{\alpha}_{\ell}(x)\right\vert >M)$\textit{ and} $M\left\Vert \lambda
(\hat{\gamma}_{\ell})-\lambda(\gamma_{0})\right\Vert
\overset{p}{\longrightarrow}0,$ \textit{then}
\[
\sqrt{n}(\hat{\theta}-\theta_{0})\overset{d}{\longrightarrow}N(0,V)\text{,
}\hat{V}\overset{p}{\longrightarrow}V.
\]

\bigskip

For functionals of a nonparametric 2SLS estimator this result is about
debiased GMM rather than the plug-in series estimator without cross-fitting
considered in Ai and Chen (2007). Theorem 19 is more general than Theorem 4.1
of Ai and Chen (2007) in applying to any first step estimator rather than just
a series estimator. Theorem 19 is simpler in only requiring mean square
convergence rates rather than the many Assumptions 3.1-3.8 and 4.1-4.6 of Ai
and Chen (2007). The estimator $\hat{\theta}$ and Theorem 19 is more
complicated in involving construction of and properties for $\hat{\alpha
}_{\ell},$ but some such $\hat{\alpha}_{\ell}$ is needed in any case for
$\hat{V}$. This comparison would also apply to a host of plug-in estimators
and root-n consistency and asymptotic normality results in the literature,
including early contributions by Powell, Stock, and Stoker (1989) for kernel
density estimators and Newey (1994a) for series estimators of least squares projections.

\subsection{Example 2: Functionals of a High Dimensional Conditional Quantile}

We will specify conditions that allow us to apply Theorem 19 to this problem.
The next condition will be sufficient for condition iv) of Theorem 19 with
$\alpha_{0}(x)$ for $\lambda(w,\gamma)=v_{u}(Y-\gamma).$

\bigskip

\textsc{Assumption 5:}\textit{ i) There exists bounded }$\bar{\alpha}%
(x)$\textit{ such that }$E[m(W,\gamma)]=E[\bar{\alpha}(X)\gamma(X)]$\textit{
for all }$\gamma(X)$\textit{ with }$E[\gamma(X)^{2}]<\infty;$\textit{ ii)
}$U=Y-\gamma_{0}(X)$\textit{ is continuously distributed and there is }%
$C>0$\textit{ such that the conditional density }$f(u|X)$\textit{ of
conditional on }$X$\textit{ satisfies }$C^{-1}\leq f(0|X)\leq C$\textit{ and
is twice continuously differentiable in }$u$\textit{ with probability one with
}$\left\vert \partial^{j}f(u|X)/\partial u^{j}\right\vert \leq C,$ $(j=1,2).$

\bigskip

This condition specifies that $E[m(W,\gamma)]$ is a mean square continuous
functional of $\gamma$ with Riesz representer $\bar{\alpha}(X)$ and imposes
some restrictions on the conditional pdf of $U$ given $X$.

\bigskip

\textsc{Lemma 20: }\textit{If Assumption 5 is satisfied then }$\left\vert
\bar{\psi}(\gamma)\right\vert \leq C\left\Vert \gamma-\gamma_{0}\right\Vert
^{2}$\textit{ for }$\left\Vert \gamma-\gamma_{0}\right\Vert ^{2}=\int%
[\gamma(x)-\gamma_{0}(x)]^{2}F_{0}(dx)$\textit{ and}%
\[
\alpha_{0}(X)=f(0|X)^{-1}\bar{\alpha}(X),\text{ }\bar{\psi}(\gamma
)=E[m(W,\gamma)+\alpha_{0}(X)v_{u}(Y-\gamma(X))]-\theta_{0}.
\]

\bigskip

Here we see that $\alpha_{0}(X)=f(0|X)^{-1}\bar{\alpha}(X)$ is the ratio of
the Riesz representer $\bar{\alpha}(X)$ to the conditional pdf $f(0|X)$. This
formula for $\alpha_{0}(X)$ differs from that of Section 7 for functionals of
conditional means where $\alpha_{0}(X)$ is the Riesz representer of the linear functional.

The $\hat{\alpha}_{\ell}(X)$ given in Section 2 will estimate $\alpha_{0}(X)$
because weighting by $f(0|X)$ is incorporated in the kernel weighting included
in $\hat{Q}_{\ell}.$ This weighting allows us to avoid inverting an estimator
of $f(0|X).$ We obtain a mean square convergence rate for this $\hat{\alpha
}_{\ell}(x)$ by extending the results of Chernozkukov, Newey, and Singh (2018)
to allow kernel weighting in $\hat{Q}_{\ell}$. Because this paper is focused
on the properties of $\hat{\theta}$ we reserve the full conditions to Appendix
B, only stating here the conditions required of the kernel $K(u),$ the
bandwidth $h,$ and the regularization factor $r_{\lambda}$ in the Lasso
minimum distance estimator in equation (\ref{quan debia}).

\bigskip

\textsc{Assumption 6:} i) $K(u)$\textit{ is a symmetric bounded kernel of
order }$\kappa$\textit{ with bounded support; ii) }$h\sqrt{n}\longrightarrow
\infty;$\textit{ iii) for each }$\ell,\ell^{\prime},$\textit{ }$\left\Vert
\hat{\gamma}_{\ell,\ell^{\prime}}-\gamma_{0}\right\Vert =O_{p}(n^{-d_{\gamma}%
})$\textit{; iv) }$\sqrt{\ln(p)/(hn)}+h^{2}+n_{{}}^{-d_{\gamma}}=o(r_{\lambda
});$\textit{ v) }$r_{\lambda}\longrightarrow0.$

\bigskip

\textsc{Lemma 21:} \textit{If Assumptions 6 and B1 are satisfied then
}$\left\Vert \hat{\alpha}_{\ell}-\alpha_{0}\right\Vert =O_{p}(\sqrt
{r_{\lambda}})$\textit{. If Assumption B2 is also satisfied then for the
sparse approximation rate }$\xi\geq1/2$\textit{ from Assumption B2 we have
}$\left\Vert \hat{\alpha}_{\ell}-\alpha_{0}\right\Vert =O_{p}(r_{\lambda
}^{2\xi/(1+2\xi)})$\textit{.}

\bigskip

The following result gives conditions for asymptotic inference for the
estimator of a linear functional of a regression quantile estimator.

\bigskip

\textsc{Theorem 22:} \textit{If i) Assumptions 5 and 6 are satisfied; ii)
}$E[m(W,\gamma_{0})^{2}]<\infty$ \textit{and }$\int[m(w,\hat{\gamma
})-m(w,\gamma_{0})]^{2}F_{0}(dw)\overset{p}{\longrightarrow}0$\textit{; iii)
}$\left\Vert \hat{\gamma}_{\ell}-\gamma_{0}\right\Vert =O_{p}(n^{-d_{\gamma}%
})$\textit{ for }$1/4<d_{\gamma}<1/2;$\textit{ either iv) Assumption B1 is
satisfied and }$\sqrt{n}\sqrt{r_{\lambda}}n^{-d_{\gamma}}\longrightarrow
0$\textit{ or v) Assumptions B1 and B2 are satisfied and }$\sqrt{n}r_{\lambda
}^{2\xi/(1+2\xi)}n^{-d_{\gamma}}\longrightarrow0$\textit{ then}%
\[
\sqrt{n}(\hat{\theta}-\theta_{0})\overset{d}{\longrightarrow}N(0,V)\text{,
}\hat{V}\overset{p}{\longrightarrow}V.
\]

This result depends on the conditional quantile estimator converging at a mean
square rate that is faster than $n^{-1/4}.$ Such a rate for an $L_{1}$
regularized conditional quantile estimator is derived by Belloni and
Chernozhukov (2011).

\subsection{Example 3: Dynamic Discrete Choice}

An result that is important for the properties of $\hat{\theta}$ for dynamic
discrete choice and more generally for economic structural model is a
convergence rate for the estimator $\hat{\gamma}_{2}(x)$ of the value function
term in the choice probability. We continue to let $\left\Vert a\right\Vert
=\sqrt{\int a(w)^{\prime}a(w)F_{0}(dw)}$ denote the mean square norm. We
continue to maintain independence of observations across $i$ but allow
arbitrary dependence across $t$.

\bigskip

\textsc{Assumption 7:} \textit{i) There is }$\varepsilon>0$\textit{ such that
}$\gamma_{10}(X)\in\lbrack\varepsilon,1-\varepsilon]$\textit{, for all }%
$\ell,\ell^{\prime},$ $\hat{\gamma}_{1\ell\ell^{\prime}}(X_{t})\in
\lbrack\varepsilon,1-\varepsilon]$\textit{, and }$H(p)$\textit{ is twice
continuously differentiable on }$[\varepsilon,1-\varepsilon];$\textit{ ii) For
all }$\ell,\ell^{\prime},$ $\left\Vert \hat{\gamma}_{1\ell\ell^{\prime}%
}-\gamma_{0}\right\Vert =O_{p}(n^{-d_{1}})$, $0<d_{1}<1/2;$ \textit{iii)
Assumptions B1 and B2 are satisfied with }$\alpha_{0}(x)=\gamma_{20}(x)$
\textit{and sparse approximation rate }$\xi_{1}>1/2;$\textit{ iv) }$n^{-d_{1}%
}=o(r_{1})$\textit{ and }$r_{1}=O(n^{-d_{1}}\ln(n));$ v) $\gamma_{20}%
(X)=\sum_{j=1}^{\infty}\beta_{j0}b_{j}(X)$ \textit{with} $\sum_{j>p}\left\vert
\beta_{j0}\right\vert =O_{p}(n^{-d_{1}(2\xi_{1}-1)/(2\xi_{1}+1)}\ln(n)).$

\bigskip

In practice condition i) requires fixed trimming where $\hat{\gamma}%
_{1\ell\ell^{\prime}}(X_{t})$ is censored below by $\varepsilon$ and above by
$1-\varepsilon$, with $\varepsilon$ being known. Here and in the Theorem 24
below we impose tighter restrictions on the penalty sizes $r_{1}$, $r_{2}$,
and $r_{3}$ than needed in order to allow smaller sparse approximation rates,
e.g. $\xi_{1}$ in Assumption 7.

\bigskip

\textsc{Lemma 23:} \textit{If Assumption 7 is satisfied then for }$\gamma
_{20}(X_{t})=E[H(\gamma_{10}(X_{t+1}))|X_{t}]$ and $\ell$%
\begin{align*}
\sup_{x}\left\vert \hat{\gamma}_{2\ell}(x)-\gamma_{20}(x)\right\vert  &
=O_{p}(n^{-d_{1}(2\xi_{1}-1)/(2\xi_{1}+1)}\ln(n)),\text{ }\left\vert
\hat{\gamma}_{3\ell}-\gamma_{30}\right\vert =O_{p}(n^{-d_{1}}),\\
\left\Vert \hat{\gamma}_{2\ell}-\gamma_{20}\right\Vert  &  =O_{p}%
(n^{-d_{1}2\xi_{1}/(2\xi_{1}+1)}\ln(n)),\text{ }(\ell=1,...,L).
\end{align*}

\bigskip

We expect this result to be useful more generally for dynamic structural
models to provide a machine learner of expected value differences.

\bigskip

\textsc{Theorem 24: }\textit{If i) Assumption 7 is satisfied, ii) }%
$\Lambda(a)>0$\textit{ for all }$a\in\Re,$\textit{ }$\ln\Lambda_{a}%
(a)$\textit{ is concave, }$\Lambda(a)$\textit{ is twice differentiable with
uniformly bounded derivatives, }$D(x)$ \textit{is bounded,}
$E[D(X)D(X)^{\prime}]$ \textit{is nonsingular; iii) Assumptions B1 and B2 are
satisfied for }$\alpha_{0}(x)$\textit{ equal to each element of }%
$E[D(X_{t})\pi(a_{0}(X_{t}))\Lambda_{a}(a(X_{t})Y_{2t}/\Lambda(a(X_{t}%
))|X_{t+1}=x]$ with sparse approximation rate $\xi_{2}$ and for $E[Y_{1t}%
|X_{t+1}=x]$\textit{ with sparse approximation rate }$\xi_{3};$\textit{ and
iv) }$d_{1}>1/4;$ v)\textit{ }$1+[(2\xi_{1}-1)/(2\xi_{1}+1)]2\xi_{2}/[2\xi
_{2}+1]>1/2d_{1},$ $n^{-d_{1}(2\xi_{1}-1)/(2\xi_{1}+1)}=o(r_{2}),$\textit{ and
}$r_{2}=O(n^{-d_{1}(2\xi_{1}-1)/(2\xi_{1}+1)}\ln(n));$ \textit{vi) }$\xi
_{3}/[2\xi_{3}+1]+d_{1}>1/2$, $\sqrt{\ln(p)/n}=o(r_{3}),$ \textit{and }%
$r_{3}=O(\sqrt{\ln(p)/n}\ln(n))$; \textit{vii) }$(4\xi_{1}-1)/(2\xi
_{1}+1)>1/2d_{1};$\textit{ } \textit{then for} $V=G^{-1}E[\psi_{0}(W)\psi
_{0}(W)^{\prime}]G^{-1}$%
\[
\sqrt{n}(\hat{\theta}-\theta_{0})\overset{d}{\longrightarrow}N(0,V),\text{
}\hat{V}\overset{p}{\longrightarrow}V.
\]

\bigskip

\section{Appendix A: Proofs of Theorems}

\bigskip

\textbf{Proof of Theorem 1:} Let $\phi(w,F_{\tau}):=\phi(w,\gamma(F_{\tau
}),\alpha(F_{\tau}),\theta)$. By ii),
\[
0=(1-\tau)\int\phi(w,F_{\tau})F_{0}(dw)+\tau\int\phi(w,F_{\tau})H(dw).
\]
Dividing by $\tau$ and solving gives%
\[
\frac{1}{\tau}\int\phi(w,F_{\tau})F_{0}(dw)=-\int\phi(w,F_{\tau}%
)H(dw)+\int\phi(w,F_{\tau})F_{0}(w).
\]
Taking limits as $\tau\longrightarrow0$, $\tau>0$ it follows by iii) that%
\begin{equation}
\frac{1}{\tau}\int\phi(w,F_{\tau})F_{0}(dw)\longrightarrow-\int\phi
(w,F_{0})H(dw)+0=-\int\phi(w,F_{0})H(dw). \label{NIF differentiable}%
\end{equation}
By ii) we have $\int\phi(w,F_{0})F_{0}(dw)=0,$ so that $\int\phi(w,F_{\tau
})F_{0}(dw)$ is differentiable in $\tau$ from the right at $\tau=0$ and%
\[
\frac{d}{d\tau}\int\phi(w,F_{\tau})F_{0}(dw)=-\int\phi(w,F_{0})H(dw)=-\frac
{d}{d\tau}\int g(w,\gamma(F_{\tau}),\theta)F_{0}(dw),
\]
where the last equality follows by i). Adding $d\int g(w,\gamma(F_{\tau
}),\theta)F_{0}(dw)/d\tau$ to both sides of this equation gives equation
(\ref{lrdef}). $Q.E.D.$

\bigskip

\textbf{Proof of Theorem 2: }Since $\gamma(F_{\tau}^{\alpha})=\gamma_{0}$ is a
constant hypothesis ii)\ implies that%
\begin{align*}
E[\phi(W,\gamma_{0},\alpha,\theta)]  &  =\int\phi(w,\gamma_{0},\alpha
,\theta)F_{0}(dw)=d\int g(w,\gamma(F_{\tau}^{\alpha}),\theta)F_{\alpha
}(dw)/d\tau\\
&  =d\int g(w,\gamma_{0},\theta)F_{\alpha}(dw)/d\tau=0.\text{ }Q.E.D.
\end{align*}

\bigskip

\textbf{Proof of Theorem 3:} By ii), the chain rule for Hadamard derivatives
(see 20.9 of Van der Vaart, 1998), and by eq. \textit{(\ref{lrdef})} it
follows that for $\delta_{H}=d\gamma(F_{\tau})/d\tau,$%
\[
\bar{\psi}_{\gamma}(\delta_{H},\alpha_{0},\theta_{0})=\frac{\partial\bar{\psi
}(\gamma(F_{\tau}),\alpha_{0},\theta)}{\partial\tau}=0.
\]
Equation (\textit{\ref{orthogonality}}) follows by $\bar{\psi}_{\gamma}%
(\delta,\alpha_{0},\theta_{0})$ being a continuous linear function and iii).
Equation (\ref{orth rem}) follows by Proposition 7.3.3 of Luenberger (1969).
\textit{Q.E.D.}

\bigskip

\textbf{Proof of Theorem 4: }Let $\hat{\psi}=\sum_{i=1}^{n}\psi(W_{i}%
,\gamma_{0},\alpha_{0},\theta_{0})/n$ and $\tilde{\phi}=\sum_{\ell=1}^{L}%
\sum_{i\in I_{\ell}}\phi(W_{i},\hat{\gamma}_{\ell},\alpha_{0},\theta_{0})/n.$
Note that%
\[
\hat{g}(\theta_{0})-\hat{\psi}=\left\{  \sum_{\ell=1}^{L}\sum_{i\in I_{\ell}%
}\psi(W_{i},\hat{\gamma}_{\ell},\alpha_{0},\theta_{0})/n-\hat{\psi}\right\}
-\tilde{\phi}=o_{p}(n^{-1/2})-\tilde{\phi},
\]
by the hypothesis of Theorem 4. Therefore if $\tilde{\phi}=o_{p}(n^{-1/2})$
then $\hat{g}(\theta_{0})-\hat{\psi}=o_{p}(n^{-1/2}).$ Similarly $\tilde{\phi
}=o_{p}(n^{-1/2})-[\hat{g}(\theta_{0})-\hat{\psi}],$ so if $\hat{g}(\theta
_{0})-\hat{\psi}=o_{p}(n^{-1/2})$ then $\tilde{\phi}=o_{p}(n^{-1/2}).$
$Q.E.D.$

\bigskip

\textbf{Proof of Theorem 5:} We consider first the properties of the plug-in
estimator $\tilde{\theta}=\sum_{\ell=1}^{L}\sum_{i\in I_{\ell}}Z_{i}%
\hat{\gamma}_{\ell}(X_{i})/n$ when the distribution of $W_{i}$ is $F_{0}$. Let
$\hat{\beta}_{\ell}$ be the ordinary least squares (OLS) estimator from
regressing $Y_{i}$ on $\breve{b}(X_{i})$ for all observations other then those
in $I_{\ell},$ $\breve{G}=E[\breve{b}(X_{i})\breve{b}(X_{i})^{\prime}],$ and
$\beta_{0}=\breve{G}^{-1}E[\breve{b}(X_{i})Y_{i}]$. Then it follows in a
standard way that for $\varepsilon_{i}=Y_{i}-\breve{b}(X_{i})^{\prime}%
\beta_{0}$%
\[
\hat{\beta}_{\ell}=\beta_{0}+\breve{G}^{-1}\frac{1}{n-n_{\ell}}\sum_{i\notin
I_{\ell}}\breve{b}(X_{i})\varepsilon_{i}+o_{p}(n^{-1/2})=\beta_{0}%
+O_{p}(n^{-1/2}).
\]
Therefore by $\gamma_{0}(x)=\breve{b}(x)^{\prime}\beta_{0}$ we have%
\begin{align*}
&  \frac{1}{n}\sum_{\ell=1}^{L}\sum_{i\in I_{\ell}}Z_{i}[\hat{\gamma}_{\ell
}(X_{i})-\gamma_{0}(X_{i})]\\
&  =\frac{1}{n}\sum_{\ell=1}^{L}\sum_{i\in I_{\ell}}Z_{i}\breve{b}%
(X_{i})^{\prime}(\hat{\beta}_{\ell}-\beta_{0})=E[Z\breve{b}(X)]^{\prime}%
\sum_{\ell=1}^{L}\frac{n_{\ell}}{n}(\hat{\beta}_{\ell}-\beta_{0})\\
&  =E[\alpha_{0}(X)\breve{b}(X)]^{\prime}\breve{G}^{-1}\sum_{\ell=1}^{L}%
\frac{n_{\ell}}{n}\frac{1}{n-n_{\ell}}\sum_{i\notin I_{\ell}}\breve{b}%
(X_{i})\varepsilon_{i}+o_{p}(n^{-1/2})=\frac{1}{n}\sum_{i=1}^{n}\bar{\alpha
}(X_{i})\varepsilon_{i}+o_{p}(n^{-1/2}).
\end{align*}
It then follows that%
\begin{align*}
\tilde{\theta}  &  =\theta_{0}+\frac{1}{n}\sum_{\ell=1}^{L}\sum_{i\in I_{\ell
}}Z_{i}[\hat{\gamma}_{\ell}(X_{i})-\gamma_{0}(X_{i})]+\frac{1}{n}\sum
_{i=1}^{n}Z_{i}\gamma_{0}(X_{i})-\theta_{0}\\
&  =\theta_{0}+\frac{1}{n}\sum_{i=1}^{n}\zeta(W_{i})+o_{p}(n^{-1/2}).
\end{align*}
Therefore, by the Slutzky Theorem,%
\[
\sqrt{n}(\tilde{\theta}-\theta_{0})\overset{d}{\longrightarrow}N(0,V)\text{,
}V=E[\zeta(W)\zeta(W)^{\prime}].
\]

Next let $F_{n}$ be the CDF defined in the Theorem. For this distribution,
since the distribution under $F_{n}$ of $(X,Z)$ is the same as under $F_{0},$%
\begin{align*}
E_{F_{n}}[Y|X]  &  =\gamma_{0}(X)+\left(  \frac{\mu}{\sqrt{n}}\right)
[\bar{\alpha}(x)-\alpha_{0}(x)],\\
\theta_{n}  &  =E[ZE_{F_{n}}[Y|X]]=E[\alpha_{0}(X)E_{F_{n}}[Y|X]]=\theta
_{0}-\left(  \frac{\mu}{\sqrt{n}}\right)  \bar{\sigma}^{2}.
\end{align*}
By $E[\bar{\alpha}(X)\{\alpha_{0}(X)-\bar{\alpha}(X)\}]=0.$ Note also that%
\[
E_{F_{n}}[\zeta(W)]=E[\bar{\alpha}(X)\{E_{F_{n}}[Y|X]-\gamma_{0}(X)\}]=\left(
\frac{\mu}{\sqrt{n}}\right)  E[\bar{\alpha}(X)\{\alpha_{0}(X)-\bar{\alpha
}(X)\}]=0.
\]
By hypothesis iv) the conditions of Proposition 1 of Bickel et al. (1993) are
satisfied for the parametric model $f_{0}(y-\delta\lbrack\bar{\alpha
}(x)-\alpha_{0}(x)]|x,z)$ with parameter $\delta.$ Then by Proposition 3 of
Bickel et al. (1993) the sequence of distributions where $W_{1},...,W_{n}$ are
i.i.d. with CDF $F_{n}$ are contiguous to the sequence where $W_{1},...,W_{n}$
are i.i.d. with CDF $F_{0}.$ Therefore the first conclusion of Theorem 5 holds
under $F_{n}$ also and%
\begin{align*}
\sqrt{n}(\hat{\theta}-\theta_{n})  &  =\sqrt{n}(\hat{\theta}-\theta_{0}%
)+\sqrt{n}(\theta_{0}-\theta_{n})=\frac{1}{\sqrt{n}}\sum_{i=1}^{n}\zeta
(W_{i})+o_{p}(1)+\mu\bar{\sigma}^{2}\\
&  =\frac{1}{\sqrt{n}}\sum_{i=1}^{n}\{\zeta(W_{i})-E_{F_{n}}[\zeta
(W)]\}+\mu\bar{\sigma}^{2}\overset{d}{\longrightarrow}N(\mu\bar{\sigma}%
^{2},V).\text{ }Q.E.D.
\end{align*}

\bigskip

\textbf{Proof of Theorem 6:} Let $\tilde{\alpha}_{\ell}(x)=\alpha_{0}(x)$ so
that $\phi(w,\hat{\gamma}_{\ell},\hat{\alpha}_{\ell})=\alpha_{0}%
(x)[y-\hat{\gamma}_{\ell}(x)].$ Consider
\[
\hat{\phi}_{\ell}=\frac{1}{\sqrt{n}}\sum_{i\in I_{\ell}}\alpha_{0}%
(X_{i})[Y_{i}-\hat{\gamma}_{\ell}(X_{i})],\text{ }\check{\beta}=\arg
\min_{\breve{\beta}}\frac{1}{n-n_{\ell}}\sum_{i\notin I_{\ell}}[Y_{i}%
-\breve{b}(X_{i})^{\prime}\breve{\beta}]^{2}+2r\sum_{j=1}^{s}\left\vert
\breve{\beta}_{j}\right\vert ,\text{ }\check{\gamma}_{\ell}(x)=\breve
{b}(x)^{\prime}\check{\beta}.
\]
By hypothesis iii) we $\hat{\gamma}_{\ell}(x)=\check{\gamma}_{\ell}(x)$ with
probability approaching one (w.p.a.1). Also, by standard proof of consistency
for convex objective functions $\check{\beta}$ $\overset{p}{\longrightarrow
}\breve{\beta},$ so that w.p.a.1 $\check{\beta}$ satisfies the first order
conditions
\[
\frac{1}{n-n_{\ell}}\sum_{i\notin I_{\ell}}\breve{b}(X_{i})[Y_{i}-\breve
{b}(X_{i})^{\prime}\check{\beta}]+r\breve{e}=0.
\]
Let $\breve{\beta}_{L}$ satisfy
\[
E[\breve{b}(X_{i})\{Y_{i}-\breve{b}(X_{i})^{\prime}\breve{\beta}_{L}%
\}+r\breve{e}]=0.
\]
Note that by $E[\breve{b}(X_{i})\{Y_{i}-\breve{b}(X_{i})^{\prime}\breve{\beta
}\}]=0$ it follows that $\breve{\beta}_{L}-\breve{\beta}=r\breve{G}^{-1}%
\breve{e},$ so that
\[
E[\alpha_{0}(X)\breve{b}(X)^{\prime}](\breve{\beta}_{L}-\breve{\beta
})=cr,\text{ }c=E[\alpha_{0}(X)\breve{b}(X)^{\prime}]\breve{G}^{-1}\breve
{e}\neq0.
\]
Also adding and subtracting terms and solving for $\sqrt{n}(\check{\beta
}-\breve{\beta}_{L})$ we obtain%
\[
\sqrt{n-n_{\ell}}(\check{\beta}-\breve{\beta}_{L})=[\frac{1}{n-n_{\ell}}%
\sum_{i\notin I_{\ell}}\breve{b}(X_{i})\breve{b}(X_{i})^{\prime}]^{-1}\left\{
\frac{1}{\sqrt{n-n_{\ell}}}\sum_{i\notin I_{\ell}}\breve{b}(X_{i}%
)[Y_{i}-\breve{b}(X_{i})^{\prime}\breve{\beta}_{L}]+r\breve{e}\right\}
=O_{p}(1).
\]
Therefore%
\begin{align*}
\hat{\phi}_{\ell}  &  =R_{1}+R_{2}+R_{3}+c\frac{n_{\ell}}{\sqrt{n}}r,\\
R_{1}  &  =\frac{1}{\sqrt{n}}\sum_{i\in I_{\ell}}\alpha_{0}(X_{i})\breve
{b}(X_{i})^{\prime}(\breve{\beta}_{L}-\check{\beta}_{L})=O_{p}(1),\\
R_{2}  &  =\frac{1}{\sqrt{n}}\sum_{i\in I_{\ell}}\{\alpha_{0}(X_{i})\breve
{b}(X_{i})-E[\alpha_{0}(X)\breve{b}(X)]\}(\breve{\beta}_{L}-\breve{\beta
})=O_{p}(1),\\
R_{3}  &  =\frac{1}{\sqrt{n}}\sum_{i\in I_{\ell}}\alpha_{0}(X_{i}%
)(Y_{i}-\gamma_{0}(X_{i}))=O_{p}(1).
\end{align*}
Therefore by the triangle inequality we have%
\[
\bar{\phi}=\sum_{\ell=1}^{L}\hat{\phi}_{\ell}=O_{p}(1)+c\frac{n}{\sqrt{n}%
}r=O_{p}(1)+c\sqrt{n}r.
\]
We also have%
\begin{align*}
\sqrt{n}\left(  \tilde{\theta}-\theta_{0}\right)   &  =\frac{1}{\sqrt{n}}%
\sum_{\ell=1}^{L}\sum_{i\in I_{\ell}}Z_{i}[\hat{\gamma}_{\ell}(X_{i}%
)-\gamma_{0}(X_{i})]+\frac{1}{\sqrt{n}}\sum_{i=1}^{n}\left[  Z_{i}\gamma
_{0}(X_{i})-\theta_{0}\right] \\
&  =R_{4}+R_{5}-\bar{\phi}+O_{p}(1),\\
R_{4}  &  =\frac{1}{\sqrt{n}}\sum_{\ell=1}^{L}\sum_{i\in I_{\ell}}\left[
Z_{i}-\alpha_{0}(X_{i})\right]  [\hat{\gamma}_{\ell}(X_{i})-\gamma_{0}%
(X_{i})],\\
R_{5}  &  =\frac{1}{\sqrt{n}}\sum_{i=1}^{n}\alpha_{0}(X_{i})[Y_{i}-\gamma
_{0}(X_{i})].
\end{align*}
Note that $\check{\beta}_{\ell}-\breve{\beta}=O_{p}(1)$ so that%
\[
R_{4}=\frac{1}{\sqrt{n}}\sum_{\ell=1}^{L}\sum_{i\in I_{\ell}}\left[
Z_{i}-\alpha_{0}(X_{i})\right]  [\breve{b}(X_{i})^{\prime}\left\{
\check{\beta}_{\ell}-\breve{\beta}\right\}  ]=O_{p}(1)
\]
by $E[\left\{  Z_{i}-\alpha_{0}(X_{i})\right\}  \breve{b}(X_{i})^{\prime}]=0$
and $\breve{b}(X_{i})$ bounded. Also $R_{5}=O_{p}(1)$. Then by the triangle
inequality we have%
\[
\sqrt{n}\left(  \tilde{\theta}-\theta_{0}\right)  =O_{p}(1)-\bar{\phi}%
=O_{p}(1)-c\sqrt{n}r,
\]
giving the first conclusion. The second conclusion follows by $c\sqrt{n}r\geq
c\sqrt{n}\sqrt{\ln(p)/n}=c\sqrt{\ln(p)}\longrightarrow\infty.$ $Q.E.D.$

\bigskip

\textbf{Proof of Theorem 7:}\ Suppose that $\psi(w,\gamma,\alpha,\theta)$ is
doubly robust. Then for any $\gamma\neq\gamma_{0},\gamma\in\Gamma$ we have%
\[
0=\bar{\psi}(\gamma,\alpha_{0},\theta_{0})=\bar{\psi}(\gamma_{0},\alpha
_{0},\theta_{0})=\bar{\psi}((1-\tau)\gamma_{0}+\tau\gamma,\alpha_{0}%
,\theta_{0}),
\]
for any $\tau.$ Therefore for any $\tau$,%
\[
\bar{\psi}((1-\tau)\gamma_{0}+\tau\gamma,\alpha_{0},\theta_{0})=0=(1-\tau
)\bar{\psi}(\gamma_{0},\alpha_{0},\theta_{0})+\tau\bar{\psi}(\gamma,\alpha
_{0},\theta_{0}),
\]
so that $\bar{\psi}(\gamma,\alpha_{0},\theta_{0})$ is affine in $\gamma.$ Also
by the previous equation $\bar{\psi}((1-\tau)\gamma_{0}+\tau\gamma,\alpha
_{0},\theta_{0})=0$ identically in $\tau$ so that
\[
\frac{d}{d\tau}\bar{\psi}((1-\tau)\gamma_{0}+\tau\gamma,\alpha_{0})=0.
\]

Next suppose that $\bar{\psi}(\gamma,\alpha_{0},\theta_{0})$ is affine
$\gamma$ and $d\bar{\psi}((1-\tau)\gamma_{0}+\tau\gamma,\alpha_{0})/d\tau=0.$
Then by $\bar{\psi}(\gamma_{0},\alpha_{0},\theta_{0})=0$, for any $\gamma
\in\Gamma,$
\begin{align*}
\bar{\psi}(\gamma,\alpha_{0},\theta_{0})  &  =\frac{d[\tau\bar{\psi}%
(\gamma,\alpha_{0},\theta_{0})]}{d\tau}=\frac{d[(1-\tau)\bar{\psi}(\gamma
_{0},\alpha_{0},\theta_{0})+\tau\bar{\psi}(\gamma,\alpha_{0},\theta_{0}%
)]}{d\tau}\\
&  =\frac{\bar{\psi}((1-\tau)\gamma_{0}+\tau\gamma,\alpha_{0})}{d\tau
}=0.\text{ }Q.E.D.
\end{align*}

\bigskip

\textbf{Proof of Theorem 8:} Double robustness implies that for any $\gamma
\in\Gamma,$
\[
0=E[\psi(W,\gamma,\alpha_{0},\theta_{0})]=E[g(W,\gamma,\theta_{0}%
)]+E[\alpha_{0}(X)\lambda(W,\gamma)],
\]
which gives the equation in the statement of the result. Also the statment in
the result implies that for all $\gamma\in\Gamma,$
\[
E[\psi(W,\gamma,\alpha_{0},\theta_{0})]=E[g(W,\gamma,\theta_{0})]+E[\alpha
_{0}(X)\lambda(W,\gamma)]=0.\text{ }Q.E.D.
\]

\bigskip

\textbf{Proof of Corollary 9:} By iterated expectations,%
\begin{align*}
E[\psi(W,\gamma,\alpha_{0},\theta_{0})]  &  =E[m(W,\gamma)]-\theta
_{0}+E[\alpha_{0}(X)\{Y-\gamma(X)\}]\\
&  =E[\alpha_{0}(X)\{\gamma(X)-\gamma_{0}(X)\}]+E[\alpha_{0}(X)\{\gamma
_{0}(X)-\gamma(X)\}]=0.\text{ }Q.E.D.
\end{align*}

\bigskip

\textbf{Proof of Corollary 10:} By $E[Y-\gamma_{0}(Z)|X]=0$ and iterated
expectations,%
\begin{align*}
E[\psi(W,\gamma,\alpha,\theta_{0})]  &  =E[m(W,\gamma)]-\theta_{0}%
+E[\alpha_{0}(X)\{Y-\gamma(Z)\}]\\
&  =E[m(W,\gamma)]-\theta_{0}+E[\alpha_{0}(X)\{Y-\gamma_{0}(Z)+\gamma
_{0}(Z)-\gamma(Z)\}]\\
&  =E[\alpha_{0}(X)\{\gamma(Z)-\gamma_{0}(Z)\}]+E[\alpha_{0}(X)\{\gamma
_{0}(Z)-\gamma(Z)\}]=0.\text{ }Q.E.D.
\end{align*}

\bigskip

\textbf{Proof of Theorem 11:}%
\begin{align*}
E[\psi(W,\gamma,\alpha,\theta_{0})]  &  =E[g(W,\gamma,\theta_{0}%
)]+E[\alpha_{0}(X)]-\int\alpha_{0}(u)\gamma(u)du\\
&  =\int\alpha_{0}(u)\{\gamma(u)-\gamma_{0}(u)\}du+\int\alpha(u)\{\gamma
_{0}(u)-\gamma(u)\}du=0.\text{ }Q.E.D.
\end{align*}

\bigskip

\textbf{Proof of Theorem 12:} If $\alpha_{0}$ is identified then $\psi
(w,\bar{\gamma},\alpha_{0},\theta)$ is identified for every $\theta$. By
double robustness%
\[
E[\psi(W,\bar{\gamma},\alpha_{0},\theta)]=0
\]
at $\theta=\theta_{0}$ and by assumption this is the only $\theta$ where this
equation is satisfied. \textit{Q.E.D.}

\bigskip

\textbf{Proof of Theorem 13: }By orthogonality of the least square projection
and by $\alpha_{0}(X)$ being a linear combination of $b(X)$ it follows that
$E[\alpha_{0}(X)\{Y-\bar{\gamma}(X)\}]=0.$ Then by Corollary 10,
\[
E[m(W,\bar{\gamma})]-\theta_{0}=-E[\alpha_{0}(X)\{Y-\bar{\gamma}%
(X)\}]=0\text{. }Q.E.D.
\]

\bigskip

\textbf{Proof of Theorem 14:} Let $\alpha_{0}(X)=-c^{\prime}\Pi^{-1}a(X)$ so
\ that $E[\alpha_{0}(X)|Z]=-c^{\prime}\Pi^{-1}\Pi p(Z)=-c^{\prime}p(Z).$Then
integration by parts gives%
\begin{align*}
E[g(W,\theta_{0},\tilde{\gamma})]  &  =E[c^{\prime}p(Z)\{\bar{\gamma
}(Z)-\gamma_{0}(Z)\}]=-E[E[\alpha_{0}(X)|Z]\{\bar{\gamma}(Z)-\gamma
_{0}(Z)\}]\\
&  =E[\alpha_{0}(X)\{Y-\bar{\gamma}(Z)\}]=-c^{\prime}\Pi^{-1}E[a(X)\{Y-\bar
{\gamma}(Z)\}]=0.\text{ }Q.E.D.
\end{align*}

\bigskip

\textbf{Proof of Lemma 15: }Define%
\begin{align}
\hat{R}_{1\ell i}  &  =g(W_{i},\hat{\gamma}_{\ell},\theta_{0})-g(W_{i}%
,\gamma_{0},\theta_{0}),\text{ }\hat{R}_{2\ell i}=\phi(W_{i},\hat{\gamma
}_{\ell},\alpha_{0},\theta_{0})-\phi(W_{i},\gamma_{0},\alpha_{0},\theta
_{0}),\label{rem def}\\
\hat{R}_{3\ell i}  &  =\phi(W_{i},\gamma_{0},\hat{\alpha}_{\ell},\hat{\theta
}_{\ell})-\phi(W_{i},\gamma_{0},\alpha_{0},\theta_{0}),\text{ }i\in I_{\ell
}.\nonumber
\end{align}
Then we have%
\begin{equation}
g(W_{i},\hat{\gamma}_{\ell},\theta_{0})+\phi(W_{i},\hat{\gamma}_{\ell}%
,\hat{\alpha}_{\ell},\tilde{\theta}_{\ell})-\psi(W_{i},\gamma_{0},\alpha
_{0},\theta_{0})=\hat{R}_{1\ell i}+\hat{R}_{2\ell i}+\hat{R}_{3\ell i}%
+\hat{\Delta}_{\ell}(W_{i}). \label{rem decomp}%
\end{equation}
Let $\mathcal{W}_{\ell}^{c}$ denote the observations not in $I_{\ell}$, so
that $\hat{\gamma}_{\ell},$ $\hat{\alpha}_{\ell},$ and $\hat{\theta}_{\ell}$
depend only on $\mathcal{W}_{\ell}^{c}$. Therefore by $E[g(W,\gamma_{0}%
,\theta_{0})]=0$,
\begin{align*}
E[\hat{R}_{1\ell i}|\mathcal{W}_{\ell}^{c}]  &  =\int g(w,\hat{\gamma}_{\ell
},\theta_{0})F_{0}(dw),\text{ }E[\hat{R}_{2\ell i}|\mathcal{W}_{\ell}%
^{c}]=\int\phi(w,\hat{\gamma}_{\ell},\alpha_{0},\theta_{0})F_{0}(dw),\\
E[\hat{R}_{3\ell i}|\mathcal{W}_{\ell}^{c}]  &  =\int\phi(w,\gamma_{0}%
,\hat{\alpha}_{\ell},\tilde{\theta}_{\ell})F_{0}(dw)=0,
\end{align*}
where the last equality follows by Assumption 3 i). Also by observations in
$I_{\ell}$ mutually independent conditional on $\mathcal{W}_{\ell}^{c}$ and
Assumption 1 i),
\[
E[\left\{  \frac{1}{\sqrt{n}}\sum_{i\in I_{\ell}}(\hat{R}_{1\ell i}-E[\hat
{R}_{1\ell i}|\mathcal{W}_{\ell}^{c}])\right\}  ^{2}|\mathcal{W}_{\ell}%
^{c}]=\frac{n_{\ell}}{n}Var(\hat{R}_{1\ell i}|\mathcal{W}_{\ell}^{c})\leq
E[\hat{R}_{1\ell i}^{2}|\mathcal{W}_{\ell}^{c}%
]\overset{}{\overset{p}{\longrightarrow}0.}%
\]
By Assumptions 3 ii) and iii) and the same argument with $\hat{R}_{2\ell i}$
and $\hat{R}_{3\ell i}$ replacing $\hat{R}_{1\ell i}$ and by the conditional
Markov inequality it follows that
\[
\frac{1}{\sqrt{n}}\sum_{i\in I_{\ell}}(\hat{R}_{1\ell i}+\hat{R}_{2\ell
i}+\hat{R}_{3\ell i}-E[\hat{R}_{1\ell i}+\hat{R}_{2\ell i}|\mathcal{W}_{\ell
}^{c}])\overset{p}{\longrightarrow}0.
\]
From equation (\ref{rem def}) we have $E[\hat{R}_{1\ell i}+\hat{R}_{2\ell
i}|\mathcal{W}_{\ell}^{c}]=\bar{\psi}(\hat{\gamma}_{\ell},\alpha_{0}%
,\theta_{0}),$ so that by Assumptions 3 ii), iii), or iv),%
\[
\left\Vert \frac{1}{\sqrt{n}}\sum_{i\in I_{\ell}}E[\hat{R}_{1\ell i}+\hat
{R}_{2\ell i}|\mathcal{W}_{\ell}^{c}]\right\Vert =\frac{n_{\ell}}{\sqrt{n}%
}\left\Vert \bar{\psi}(\hat{\gamma}_{\ell},\alpha_{0},\theta_{0})\right\Vert
\leq\sqrt{n}\left\Vert \bar{\psi}(\hat{\gamma}_{\ell},\alpha_{0},\theta
_{0})\right\Vert \overset{p}{\longrightarrow}0.
\]
Then by the triangle inequality%
\[
\frac{1}{\sqrt{n}}\sum_{i\in I_{\ell}}(\hat{R}_{1\ell i}+\hat{R}_{2\ell
i}+\hat{R}_{3\ell i})\overset{p}{\longrightarrow}0.
\]
It follows similarly from Assumption 2 that%
\[
\frac{1}{\sqrt{n}}\sum_{i\in I_{\ell}}\hat{\Delta}_{\ell i}(W_{i}%
)\overset{p}{\longrightarrow}0.
\]
Then by the triangle inequality and equation (\ref{rem decomp}),%
\[
\sqrt{n}\hat{\psi}(\theta_{0})-\frac{1}{\sqrt{n}}\sum_{i=1}^{n}\psi
(W_{i},\theta_{0},\gamma_{0},\alpha_{0})=\sum_{\ell=1}^{L}\frac{1}{\sqrt{n}%
}\sum_{i\in I_{\ell}}[\hat{R}_{1\ell i}+\hat{R}_{2\ell i}+\hat{R}_{3\ell
i}+\hat{\Delta}_{\ell i}(W_{i})]\overset{p}{\longrightarrow}0.\text{ }Q.E.D.
\]

\bigskip

\textbf{Proof of Lemma 16:} Define the remainders $\hat{R}_{1\ell i}$,
$\hat{R}_{2\ell i}$, $\hat{R}_{3\ell i},$ and $\hat{\Delta}_{\ell}(W_{i})$ as
in the proof of Lemma 15. Also, define $\hat{R}_{4\ell i}=g(W_{i},\hat{\gamma
}_{\ell},\tilde{\theta}_{\ell})-g(W_{i},\hat{\gamma}_{\ell},\theta_{0})$. Note
that by hypothesis,
\[
E[\left\Vert \hat{R}_{4\ell i}\right\Vert ^{2}|\mathcal{W}_{\ell}%
^{c}]\overset{p}{\longrightarrow}0.
\]
It similarly follows from Assumption 1 and $\int\left\Vert \hat{\Delta}_{\ell
}(w)\right\Vert ^{2}F_{0}(dw)\overset{p}{\longrightarrow}0$ that for $i\in
I_{\ell}$,
\[
E[\left\Vert \hat{R}_{k\ell i}\right\Vert ^{2}|\mathcal{W}_{\ell}%
^{c}]\overset{p}{\longrightarrow}0,\text{ }k=1,2,3,\text{ }E[\left\Vert
\hat{\Delta}_{\ell}(W_{i})\right\Vert ^{2}|\mathcal{W}_{\ell}^{c}%
]\overset{p}{\longrightarrow}0.
\]
Then it follows that for $\psi_{i}=\psi(W_{i},\gamma_{0},\alpha_{0},\theta
_{0}),$%
\[
E[\frac{1}{n}\sum_{i\in I_{\ell}}\left\Vert \hat{\psi}_{i\ell}-\psi
_{i}\right\Vert ^{2}|\mathcal{W}_{\ell}^{c}]\leq\frac{Cn_{\ell}}{n}\left(
\sum_{k=1}^{4}E[\left\Vert \hat{R}_{k\ell i}\right\Vert ^{2}|\mathcal{W}%
_{\ell}^{c}]+E[\left\Vert \hat{\Delta}_{\ell}(W_{i})\right\Vert ^{2}%
|\mathcal{W}_{\ell}^{c}]\right)  \overset{p}{\longrightarrow}0.
\]
Therefore $\sum_{i\in I_{\ell}}\left\Vert \hat{\psi}_{i\ell}-\psi
_{i}\right\Vert ^{2}/n\overset{p}{\longrightarrow}0$ by the conditional Markov
inequality. It follows by the triangle inequality that for $\tilde{\Psi}%
:=\sum_{i=1}^{n}\psi_{i}\psi_{i}^{\prime}/n,$%
\begin{align*}
\left\Vert \hat{\Psi}-\tilde{\Psi}\right\Vert  &  \leq\sum_{\ell=1}^{L}%
\frac{1}{n}\sum_{i\in I_{\ell}}(\left\Vert \hat{\psi}_{i\ell}-\psi
_{i}\right\Vert ^{2}+2\left\Vert \psi_{i}\right\Vert \left\Vert \hat{\psi
}_{i\ell}-\psi_{i}\right\Vert )\\
&  \leq o_{p}(1)+2\sum_{\ell=1}^{L}\sqrt{\frac{1}{n}\sum_{i\in I_{\ell}%
}\left\Vert \hat{\psi}_{i\ell}-\psi_{i}\right\Vert ^{2}}\sqrt{\frac{1}{n}%
\sum_{i\in I_{\ell}}\left\Vert \psi_{i}\right\Vert ^{2}}=o_{p}(1)(1+O_{p}%
(1))\overset{p}{\longrightarrow}0.
\end{align*}
We also have $\tilde{\Psi}\overset{p}{\longrightarrow}\Psi$ by Khintchine's
law of large numbers, so the conclusion follows by the triangle inequality.
$Q.E.D.$

\bigskip

\textbf{Proof of Lemma 17:} Let $\hat{G}_{\ell}=n_{\ell}^{-1}\sum_{i\in
I_{\ell}}\partial g(W_{i},\hat{\gamma}_{\ell},\bar{\theta})/\partial\theta$
and $\tilde{G}_{\ell}=n_{\ell}^{-1}\sum_{i\in I_{\ell}}\partial g(W_{i}%
,\hat{\gamma}_{\ell},\theta_{0})/\partial\theta.$ By ii) and
\[
E[\frac{1}{n_{\ell}}\sum_{i\in I_{\ell}}d(W_{i},\hat{\gamma}_{\ell
})|\mathcal{W}_{\ell}^{c}]=E[d(W_{i},\hat{\gamma}_{\ell})|\mathcal{W}_{\ell
}^{c}]\leq C\text{,}%
\]
with probability approaching one. Then by the conditional Markov inequality,
$\sum_{i\in I_{\ell}}d(W_{i},\hat{\gamma}_{\ell})=O_{p}(1).$ Then by
conditions i) and ii) and the triangle inequality, with probability
approaching one
\[
\left\Vert \hat{G}_{\ell}-\tilde{G}_{\ell}\right\Vert \leq n_{\ell}^{-1}%
\sum_{i\in I_{\ell}}d(W_{i},\hat{\gamma}_{\ell})\left\Vert \bar{\theta}%
-\theta_{0}\right\Vert ^{1/C}=O_{p}(1)o_{p}(1)\overset{p}{\longrightarrow}0.
\]
Then $\hat{G}_{\ell}-\tilde{G}_{\ell}\overset{p}{\longrightarrow}0$ follows by
the conditional Markov inequality. For $\bar{G}_{\ell}=n_{\ell}^{-1}\sum_{i\in
I_{\ell}}\partial g(W_{i},\gamma_{0},\theta_{0})/\partial\theta$ it follows
similarly from condition iii) that $\tilde{G}_{\ell}-\bar{G}_{\ell
}\overset{p}{\longrightarrow}0.$ By Khintchine's law of large numbers $\bar
{G}_{\ell}\overset{p}{\longrightarrow}G$, so the conclusion follows by the
triangle inequality. $Q.E.D.$

\bigskip

\textbf{Proof of Theorem 18}: Follows in a standard way from Lemmas 15-17.
$Q.E.D.$

\bigskip

\textbf{Proof of Lemma 19}: Let $g(w,\gamma,\theta)=m(w,\gamma)-\theta$ and
$\phi(w,\gamma,\alpha,\theta)=\alpha(x)\lambda(w,\gamma).$ Assumption 1 is
satisfied by conditions ii) and iii) and Assumption 3 is satisfied by
condition v). Also, note that%
\[
\hat{\Delta}_{\ell}(w)=[\hat{\alpha}_{\ell}(x)-\alpha_{0}(x)][\lambda
(w,\hat{\gamma}_{\ell})-\lambda(w,\gamma_{0})].
\]
Therefore by $\alpha_{0}(X)$ bounded, $\sup_{x}\left\vert \hat{\alpha}_{\ell
}(x)\right\vert =O_{p}(1)$, and iii),%
\begin{align*}
\int\hat{\Delta}_{\ell}(w)^{2}F_{0}(dw)  &  =\int[\hat{\alpha}_{\ell
}(x)-\alpha_{0}(x)]^{2}[\lambda(w,\hat{\gamma}_{\ell})-\lambda(w,\gamma
_{0})]^{2}F_{0}(dw)\\
&  \leq O_{p}(1)\int[\lambda(w,\hat{\gamma}_{\ell})-\lambda(w,\gamma_{0}%
)]^{2}F_{0}(dw)\overset{p}{\longrightarrow}0.
\end{align*}
Also by iterated expectations, the Cauchy Schwartz inequality, and condition
iv),%
\begin{align*}
\left\vert \sqrt{n}\int\hat{\Delta}_{\ell}(w)F_{0}(dw)\right\vert  &
=\sqrt{n}\left\vert \int[\hat{\alpha}_{\ell}(x)-\alpha_{0}(x)][\lambda
(w,\hat{\gamma}_{\ell})-\lambda(w,\gamma_{0})]F_{0}(dw)\right\vert \\
&  =\sqrt{n}\left\vert \int[\hat{\alpha}_{\ell}(x)-\alpha_{0}(x)][\bar
{\lambda}(x,\hat{\gamma}_{\ell})-\bar{\lambda}(x,\gamma_{0})]F_{0}%
(dx)\right\vert \\
&  \leq\sqrt{n}\left\Vert \hat{\alpha}_{\ell}-\alpha_{0}\right\Vert \left\Vert
\bar{\lambda}(\hat{\gamma}_{\ell})-\bar{\lambda}(\gamma_{0})\right\Vert
\overset{p}{\longrightarrow}0.
\end{align*}
Therefore Assumption 2 ii) is satisfied, so all the conditions of Lemma 15 are
satisfied, giving the first conclusion. In addition the conditions of Lemma 16
are satisfied so the second conclusion follows from Lemma 16. $Q.E.D.$

\bigskip

\textbf{Proof of Lemma 20: }By $\gamma_{0}(X)$ being the $\lambda^{th}$
conditional quantile of $Y$ given $X$ we have
\[
E[1(U<0)|X]=E[1(Y<\gamma_{0}(X))|X]=\lambda.
\]
Then by Assumption 5 and Taylor expansion with LaGrange remainder%
\begin{align*}
E[v_{u}(Y-\gamma(X))|X]  &  =\lambda-\int_{-\infty}^{\gamma(X)-\gamma_{0}%
(X)}f(u|X)du\\
&  =-f(0|X)[\gamma(X)-\gamma_{0}(X)]-[\partial f(\delta(X)|X)/\partial
u][\gamma(X)-\gamma_{0}(X)]^{2},
\end{align*}
where $\delta(X)$ is between $\gamma(X)-\gamma_{0}(X)$ and zero. Therefore%
\begin{align*}
\left\vert \bar{\psi}(\gamma)\right\vert  &  =|E[\bar{\alpha}(X)\{\gamma
(X)-\gamma_{0}(X)\}]-E[\alpha_{0}(X)f_{0}(0|X)\{\gamma(X)-\gamma_{0}(X)\}]\\
&  +E[\{\partial f(\delta(X)|X)/\partial u\}\{\gamma(X)-\gamma_{0}%
(X)\}^{2}]|\\
&  \leq C\left\Vert \gamma-\gamma_{0}\right\Vert ^{2}.\text{ }Q.E.D.
\end{align*}

\bigskip

\textbf{Proof of Lemma 21:} Given in Appendix B.

\bigskip

\textbf{Proof of Theorem 22: }We proceed by showing that each of the
conditions of Lemma 19 are satisfied for $\lambda(W,\gamma)=\lambda
-1(Y<\gamma(X)),$ $\gamma_{0}(X)$ the $\lambda^{th}$ conditional quantile of
$Y$, and $\alpha_{0}$ given in Lemma 20. Condition i) of Lemma 19 holds by the
definition of $\gamma_{0}(X).$ Condition ii)\ of Lemma 19 holds by Assumption
5, $\lambda(W,\gamma)$ bounded, and hypothesis ii). Condition iii) of Lemma 19
holds by hypothesis, Lemma 21, and%
\begin{align*}
\int[\lambda(w,\hat{\gamma}_{\ell})-\lambda(w,\gamma_{0})]^{2}F_{0}(dw)  &
\leq\int1(\left\vert u\right\vert \leq\left\vert \hat{\gamma}_{\ell}%
(x)-\gamma_{0}(x)\right\vert )f(u|x)F_{0}(dx)\\
&  \leq C\int\left\vert \hat{\gamma}_{\ell}(x)-\gamma_{0}(x)\right\vert
F_{0}(dx)\leq C\left\Vert \hat{\gamma}_{\ell}-\gamma_{0}\right\Vert
\overset{p}{\longrightarrow}0.
\end{align*}
For condition iv), note that by $\left\vert \lambda(w,\gamma)\right\vert
\leq1,$%
\[
\int[\hat{\alpha}_{\ell}(x)-\alpha_{0}(x)]^{2}[\lambda(w,\hat{\gamma}_{\ell
})-\lambda(w,\gamma_{0})]^{2}F_{0}(dw)\leq4\int[\hat{\alpha}_{\ell}%
(x)-\alpha_{0}(x)]^{2}F_{0}(dw)\overset{p}{\longrightarrow}0.
\]
Also, note that by $f(u|x)\leq C$,
\begin{align*}
\left\Vert \bar{\lambda}(\hat{\gamma}_{\ell})-\bar{\lambda}(\gamma
_{0})\right\Vert ^{2}  &  =\int[\int_{-\infty}^{\hat{\gamma}(x)-\gamma_{0}%
(x)}f(u|x)du-\int_{-\infty}^{0}f(u|x)du]^{2}F_{0}(dx)\\
&  \leq\int C|\hat{\gamma}(x)-\gamma_{0}(x)|^{2}F_{0}(dx)=C\left\Vert
\hat{\gamma}-\gamma_{0}\right\Vert ^{2}.
\end{align*}
In addition,
\begin{align*}
\left\Vert \bar{\lambda}(\hat{\gamma}_{\ell})-\bar{\lambda}(\gamma
_{0})\right\Vert ^{2}  &  =\int[\int_{-\infty}^{\hat{\gamma}(x)-\gamma_{0}%
(x)}f(u|x)du-\int_{-\infty}^{0}f(u|x)du]^{2}F_{0}(dx)\\
&  \leq\int C|\hat{\gamma}(x)-\gamma_{0}(x)|^{2}F_{0}(dx)=C\left\Vert
\hat{\gamma}-\gamma_{0}\right\Vert ^{2}.
\end{align*}
Condition iv)\ of Lemma 19 then follows by Lemma 21 and hypothesis iv) or v).
Condition v) of Lemma 19 also follows by Lemma 20. The conclusion then follows
by the conclusion of Lemma 19. $Q.E.D.$

\bigskip

\textbf{Proof of Lemma 23:} Let $Q_{2}=E[b(X)b(X)^{\prime}\gamma_{10}(X)].$ A
standard maximal inequality gives
\[
\left\vert \hat{Q}_{2\ell}-Q_{2}\right\vert _{\infty}=O_{p}(\sqrt{\ln(p)/n}).
\]
Let
\[
\tilde{M}_{2\ell}=\frac{1}{(n-n_{\ell})T}\sum_{\ell^{\prime}\neq\ell}%
\sum_{i\in I_{\ell^{\prime}}}\sum_{t=1}^{T}Y_{2it}b(X_{it})H\left(
\gamma_{10}(X_{i,t+1}\right)  ).
\]
and $\mathcal{W}_{\ell}^{c}$ denote all observations not in $I_{\ell}.$ It
follows by Assumption 7 and the Cauchy-Schwartz inequalities that%
\[
\frac{1}{n_{\ell^{\prime}}T}E[\sum_{i\in I_{\ell^{\prime}}}\sum_{t=1}%
^{T}\left\vert \hat{\gamma}_{1\ell,\ell^{\prime}}(X_{i,t+1})-\gamma
_{10}(X_{i,t+1})\right\vert |\mathcal{W}_{\ell}^{c}]=\int\left\vert
\hat{\gamma}_{1\ell,\ell^{\prime}}(x)-\gamma_{10}(x)\right\vert F_{0}%
(dx)\leq\left\Vert \hat{\gamma}_{1\ell,\ell^{\prime}}-\gamma_{10}\right\Vert
=O_{p}(n^{-d_{1}}).
\]
Then by $H(p)$ having bounded derivative on $[\varepsilon,1-\varepsilon],$
Assumption 7, and the conditional Markov inequality%
\begin{align}
\left\vert \hat{M}_{2\ell}-\tilde{M}_{2\ell}\right\vert _{\infty}  &
\leq\frac{C}{(n-n_{\ell})}\sum_{\ell^{\prime}\neq\ell}\sum_{i\in
I_{\ell^{\prime}}}\sum_{t=1}^{T}\left\vert \hat{\gamma}_{1\ell,\ell^{\prime}%
}(X_{i,t+1})-\gamma_{10}(X_{i,t+1})\right\vert \label{gamma1 rem}\\
&  \leq\frac{C}{(n-n_{\ell})}\sum_{\ell^{\prime}\neq\ell}n_{\ell^{\prime}%
}\frac{1}{n_{\ell^{\prime}}T}\sum_{i\in I_{\ell^{\prime}}}\sum_{t=1}%
^{T}\left\vert \hat{\gamma}_{1\ell,\ell^{\prime}}(X_{i,t+1})-\gamma
_{10}(X_{i,t+1})\right\vert \nonumber\\
&  \leq\frac{C}{(n-n_{\ell})}\sum_{\ell^{\prime}\neq\ell}n_{\ell^{\prime}%
}O_{p}(n^{-d_{1}})=O_{p}(n^{-d_{1}}).\nonumber
\end{align}
For $M_{2}=E[Y_{2t}b(X_{t})H(\gamma_{10}(X_{t+1}))]$ it follows by a standard
maximal inequality that $\left\vert \tilde{M}_{2\ell}-M_{2}\right\vert
_{\infty}=O_{p}(\sqrt{\ln(p)/n}).$ Then by the triangle inequality we have
$\left\vert \hat{M}_{2\ell}-M_{2}\right\vert _{\infty}=O_{p}(n^{-d_{1}}).$
Therefore it follows that%
\[
\left\vert \hat{M}_{2\ell}-M_{2}\right\vert _{\infty}=O_{p}(\varepsilon
_{n}),\text{ }\left\vert \hat{Q}_{2\ell}-Q_{2}\right\vert _{\infty}%
=O_{p}(\varepsilon_{n}),\text{ }\varepsilon_{n}=n^{-d_{1}}.
\]
Let $\beta_{L}=\arg\min\left\Vert \gamma_{20}-\beta^{\prime}b\right\Vert
^{2}+2r_{2}\left\vert \beta\right\vert _{1}.$ It then follows as in the proof
of Theorem 3 of Chernozhukov, Newey, and Singh (2018) that%
\begin{align*}
\sup_{x}|\hat{\gamma}_{2\ell}(x)-\beta_{L}^{\prime}b(x)|  &  \leq\max_{j}%
\sup_{x}\left\vert b_{j}(x)\right\vert \sum_{j=1}^{p}\left\vert \hat{\beta
}_{j}-\beta_{Lj}\right\vert =O_{p}((\varepsilon_{n}^{2})^{-1/(2\xi_{1}%
+1)}r_{2})\\
&  =O_{p}(n^{-d_{1}(2\xi_{1}-1)/(2\xi_{1}+1)}\ln(n)).
\end{align*}
It also follows similarly to this result that%
\[
\sup_{x}|\beta_{L}^{\prime}b(x)-\beta_{0}^{\prime}b(x)|\leq\max_{j}\sup
_{x}\left\vert b_{j}(x)\right\vert \sum_{j=1}^{p}\left\vert \beta_{Lj}%
-\beta_{j}\right\vert =O_{p}(n^{-d_{1}(2\xi_{1}-1)/(2\xi_{1}+1)}\ln(n)).
\]
The first conclusion then follows by hypothesis v) and the triangle inequality.

For the second conclusion let%
\[
\tilde{\gamma}_{3\ell}=\frac{1}{\hat{P}_{1\ell}(n-n_{\ell})T}\sum
_{\ell^{\prime}\neq\ell}\sum_{i\in I_{\ell^{\prime}}}\sum_{t=1}^{T}%
Y_{1it}H\left(  \gamma_{10}(X_{i,t+1}\right)  ).
\]
It follows similarly to equation (\ref{gamma1 rem}) that $\left\vert
\hat{\gamma}_{3\ell}-\tilde{\gamma}_{3\ell}\right\vert =O_{p}(n^{-d_{1}})$.
Also by standard arguments $\left\vert \tilde{\gamma}_{3\ell}-\gamma
_{30}\right\vert =O_{p}(1/\sqrt{n})=O_{p}(n^{-d_{1}}),$ so the second
conclusion follows by the triangle inequality.

The third conclusion follows from Theorem 3 of Chernozhukov, Newey, and Singh
(2018) similarly to the first conclusion. $Q.E.D.$

\bigskip

Before proving Theorem 24 we first prove some useful Lemmas.

\bigskip

\textsc{Lemma A1:} \textit{If Assumption 7 and the hypotheses of Theorem 24
are satisfied then }%
\[
\tilde{\theta}_{\ell}=\theta_{0}+O_{p}(n^{-d_{1}[(2\xi_{1}-1)/(2\xi_{1}%
+1)]}).
\]
regressors $b(X_{i,t+1}),$ dependent variables equal to each element of
$\hat{\alpha}_{2\ell}(X_{it},Y_{2it},\tilde{\theta}_{\ell})$ and

\bigskip

\textbf{Proof of Lemma A1:} Follows from Lemma 23 by standard arguments for
maximum likelihood with $\hat{\gamma}_{2}(x)$ and $\hat{\gamma}_{3}$
plugged-in. $Q.E.D.$

\bigskip

\textsc{Lemma A2:} \textit{If Assumption 7 and the hypotheses of Theorem 24
are satisfied then}%
\begin{align*}
\left\Vert \hat{\alpha}_{1\ell}-\alpha_{10}\right\Vert  &  =O_{p}%
(n^{-d_{1}[(2\xi_{1}-1)/(2\xi_{1}+1)]2\xi_{2}/[2\xi_{2}+1]}[\ln(n)]^{2}%
+n^{-\xi_{3}/(2\xi_{3}+1)}\ln(n)+n^{-d_{1}}),\\
\left\Vert \hat{\alpha}_{2\ell}-\alpha_{20}\right\Vert  &  =O_{p}%
(n^{-d_{1}[(2\xi_{1}-1)/(2\xi_{1}+1)]}\ln(n)),\text{ }\left\Vert \hat{\alpha
}_{3\ell}-\alpha_{30}\right\Vert =O_{p}(n^{-d_{1}(2\xi_{1}-1)/(2\xi_{1}%
+1}),\text{ }\left\Vert H(\hat{\gamma}_{1\ell})-H(\gamma_{10})\right\Vert
=O_{p}(n^{-d_{1}}).
\end{align*}

\textbf{Proof of Lemma A2: }First, note that $a(x)=a(x,\theta_{0},\gamma
_{20},\gamma_{30})$ is bounded by $D(x)$ and $H(\gamma_{10}(x))$ is bounded
$\gamma_{10}(x)\in(\varepsilon,1-\varepsilon)$, so that by Assumption 7 and
the fixed trimming, with $\Lambda(a)>0$ and twice continuous differentiability
of $\Lambda(a)$,
\begin{align*}
\sup_{x}\left\vert \hat{\alpha}_{2\ell k}(x,y_{2})-\alpha_{02k}(x,y_{2}%
)\right\vert  &  \leq C\sup_{x}\left\vert \hat{a}(x)-a(x)\right\vert \leq
C(\left\Vert \tilde{\theta}_{\ell}-\theta_{0}\right\Vert +\sup_{x}\left\vert
\hat{\gamma}_{2\ell}(x)-\gamma_{20}(x)\right\vert \\
+\left\vert \hat{\gamma}_{3\ell}-\gamma_{30}\right\vert )  &  =O_{p}%
(n^{-d_{1}(2\xi_{1}-1)/(2\xi_{1}+1)}),\text{ \ }a(x,\tilde{\theta}_{\ell}%
,\hat{\gamma}_{2\ell},\hat{\gamma}_{3\ell})=\hat{a}(x),
\end{align*}
giving the second conclusion. The third conclusion follows in a standard way.

Let $\hat{\zeta}_{1\ell k}(x)$ denote Lasso with regressors $b(X_{i,t+1}),$
dependent variable equal to the $k_{th}$ element $\hat{\alpha}_{2\ell itk}$ of
$\hat{\alpha}_{2\ell it}=\hat{\alpha}_{2\ell}(X_{it},Y_{2it}),$ and
penalization $r_{2}$ for $i\notin I_{\ell}.$ This estimator has the same form
given in eq. (\ref{gam2}) with
\[
\hat{M}_{\ell}=\frac{1}{(n-n_{\ell})T}\sum_{i\notin I_{\ell}}\sum_{t=1}%
^{T}\hat{\alpha}_{2\ell itk}b(X_{i,t+1}),\text{ }\hat{Q}_{\ell}=\frac
{1}{(n-n_{\ell})T}\sum_{i\notin I_{\ell}}\sum_{t=1}^{T}b(X_{i,t+1}%
)b(X_{i,t+1})^{\prime},
\]
and $r_{2}$ replacing $r_{1}.$ By Assumption 7, Lemma 23, Lemma A1, uniform
boundedness of the elements of $b(x)$ (from Assumption B1), and Bernstein's
inequality (using independence across $i)$,
\begin{align*}
\left\vert \hat{M}_{\ell}-\tilde{M}_{\ell}\right\vert _{\infty}  &  \leq
C\sup_{x}\left\vert \hat{\alpha}_{2\ell k}(x,y_{2})-\alpha_{02k}%
(x,y_{2})\right\vert =O_{p}(n^{-d_{1}(2\xi_{1}-1)/(2\xi_{1}+1)}\ln(n)),\\
\left\vert \tilde{M}_{\ell}-M\text{ }\right\vert _{\infty}  &  =O_{p}%
(\sqrt{\ln(p)/n}),\text{ }\tilde{M}_{\ell}=\frac{1}{(n-n_{\ell})T}%
\sum_{i\notin I_{\ell}}\sum_{t=1}^{T}\alpha_{20k}(X_{it},Y_{2it}%
)b(X_{i,t+1})\\
M  &  =E[\alpha_{20k}(X_{it},Y_{2it})b(X_{i,t+1})].
\end{align*}
Then by the triangle inequality and another application of Bernstein's
inequality,%
\[
\left\vert \hat{M}_{\ell}-M\right\vert _{\infty}=O_{p}(n^{-d_{1}(2\xi
_{1}-1)/(2\xi_{1}+1)}\ln(n)),\text{ }\left\vert \hat{Q}_{\ell}-Q\right\vert
_{\infty}=O_{p}(\sqrt{\ln(p)/n}).
\]
It then follows analogously to the proof of Lemma 23 that%
\[
\left\Vert \hat{\zeta}_{1\ell}-\zeta_{10}\right\Vert =O_{p}(n^{-d_{1}%
[(2\xi_{1}-1)/(2\xi_{1}+1)]2\xi_{2}/[2\xi_{2}+1]}[\ln(n)]^{2}).
\]
It follows similarly that for $\hat{\zeta}_{2\ell}(x)$ denoting Lasso with
regressors $b(X_{i,t+1})$ and dependent variable $Y_{1it},$%
\[
\left\Vert \hat{\zeta}_{2\ell}-\zeta_{20}\right\Vert =O_{p}(n^{-\xi_{3}%
/(2\xi_{1}+1)}\ln(n)).
\]
Also note that $H_{p}(\hat{\gamma}_{1\ell}(x))$ and $H_{p}(\gamma_{10}(x))$
are bounded by the fixed trimming, which together with $\left\Vert \hat
{\gamma}_{1\ell}-\gamma_{10}\right\Vert =O_{p}(n^{-d_{1}})$ also gives%
\[
\left\Vert H_{p}(\hat{\gamma}_{1\ell})-H_{p}(\gamma_{10})\right\Vert
=O_{p}(n^{-d_{1}}).
\]
Then by the triangle inequality and boundedness of $\zeta_{10}(x)$ and
$\zeta_{20}(x),$
\begin{align*}
\left\Vert \hat{\alpha}_{1\ell}-\alpha_{10}\right\Vert  &  =\left\Vert
(\hat{\zeta}_{1\ell}+\hat{\alpha}_{3\ell}\hat{\zeta}_{2\ell})H_{p}(\hat
{\gamma}_{1\ell})-(\zeta_{10}+\alpha_{30}\zeta_{20})H_{p}(\gamma
_{10})\right\Vert \\
&  \leq\left\Vert \left[  (\hat{\zeta}_{1\ell}-\zeta_{10})+\hat{\alpha}%
_{3}(\hat{\zeta}_{2\ell}-\zeta_{20})\right]  H_{p}(\hat{\gamma}_{1\ell
})\right\Vert \\
&  +\left\Vert (\zeta_{10}+\hat{\alpha}_{3\ell}\zeta_{20})\left[  H_{p}%
(\hat{\gamma}_{1\ell})-H_{p}(\gamma_{10})\right]  \right\Vert +\left\Vert
(\hat{\alpha}_{3\ell}-\alpha_{30})\zeta_{20}H_{p}(\gamma_{10})\right\Vert \\
&  =O_{p}(n^{-d_{1}[(2\xi_{1}-1)/(2\xi_{1}+1)]2\xi_{2}/[2\xi_{2}+1]}%
[\ln(n)]^{2}+n^{-\xi_{3}/(2\xi_{1}+1)}\ln(n))+O_{p}(n^{-d_{1}}).
\end{align*}
The first conclusion follows by the triangle inequality. The last conclusion
follows similarly to $\left\Vert H_{p}(\hat{\gamma}_{1\ell})-H_{p}(\gamma
_{10})\right\Vert =O_{p}(n^{-d_{1}}).$ $Q.E.D.$

\bigskip

\textbf{Proof of Theorem 24:} We proceed by verifying Assumptions 1-4 in
Section 8 and the conditions of Lemma 16 for $\gamma=(\gamma_{1},\gamma
_{2},\gamma_{3}).$ Assumption 1 follows by Lemmas 23 and A1 and by $a(x)$,
$Y_{2t},$ $H(\gamma_{10}(X_{t})),$ $\gamma_{20}(X_{t}),$ $\alpha_{10}(X_{t})$,
$\alpha_{20}(X_{t},Y_{2t}),$ $\alpha_{30}(X_{t})$ all bounded, similarly to
the proof of Lemma A2.

To show Assumption 2 ii), note that%
\begin{align*}
\hat{\Delta}_{\ell}(w)  &  =\hat{\Delta}_{\ell1}(w)+\hat{\Delta}_{\ell
2}(w)+\hat{\Delta}_{\ell3}(w),\\
\hat{\Delta}_{\ell1}(w)  &  =-\frac{1}{T}\sum_{t=1}^{T}[\hat{\alpha}_{1\ell
}(x_{t})-\alpha_{10}(x_{t})][\hat{\gamma}_{1\ell}(x_{t})-\gamma_{10}%
(x_{t})],\\
\hat{\Delta}_{\ell2}(w)  &  =\frac{1}{T}\sum_{t=1}^{T}[\hat{\alpha}_{2\ell
}(x_{t},y_{2t})-\alpha_{20}(x_{t},y_{2t})][H(\hat{\gamma}_{1\ell}(x_{t}%
))-\hat{\gamma}_{2\ell}(x_{t})-H(\gamma_{10}(x_{t}))+\gamma_{20}(x_{t})],\\
\hat{\Delta}_{\ell3}(w)  &  =\frac{1}{T}\sum_{t=1}^{T}(\hat{\alpha}_{3\ell
}-\alpha_{30})y_{2t}\{H(\hat{\gamma}_{1}(x_{t+1}))-\hat{\gamma}_{3}%
-H(\hat{\gamma}_{1}(x_{t+1}))+\gamma_{30}\}.
\end{align*}
Then by the first conclusion of Lemma A2 and conditions iv), v), and vi),
\begin{align*}
\sqrt{n}\int\left\Vert \hat{\Delta}_{\ell1}(w)\right\Vert F_{0}(dw)  &
\leq\sqrt{n}\left\Vert \hat{\alpha}_{1\ell}-\alpha_{10}\right\Vert \left\Vert
\hat{\gamma}_{1\ell}-\gamma_{10}\right\Vert \\
&  =O_{p}(\sqrt{n}\{n^{-d_{1}[(2\xi_{1}-1)/(2\xi_{1}+1)]2\xi_{2}/[2\xi_{2}%
+1]}[\ln(n)]^{2}+n^{-\xi_{3}/(2\xi_{3}+1)}\ln(n)+n^{-d_{1}}\}n^{-d_{1}})\\
&  =o_{p}^{{}}(1).
\end{align*}
Next, by the second conclusion of Lemma A2 and condition vii)
\begin{align*}
\sqrt{n}\int\left\Vert \hat{\Delta}_{\ell2}(w)\right\Vert F_{0}(dw)  &
\leq\sqrt{n}\left\Vert \hat{\alpha}_{2\ell}-\alpha_{20}\right\Vert \left\Vert
H(\hat{\gamma}_{1\ell})-\hat{\gamma}_{2\ell}-H(\gamma_{10})+\gamma
_{20}\right\Vert \\
&  =O_{p}(\sqrt{n}[n^{-d_{1}(2\xi_{1}-1)/(2\xi_{1}+1)}\ln(n)]n^{-d_{1}2\xi
_{1}/(2\xi_{1}+1)}\ln(n))=o_{p}^{{}}(1).
\end{align*}
Next, by the third conclusion of Lemma A2, condition v), and $2\xi_{1}%
/(2\xi_{1}+1)<1$ we have
\begin{align*}
\sqrt{n}\int\left\Vert \hat{\Delta}_{\ell2}(w)\right\Vert F_{0}(dw)  &
\leq\sqrt{n}\left\Vert \hat{\alpha}_{3\ell}-\alpha_{30}\right\Vert \left\Vert
H(\hat{\gamma}_{1\ell})-\hat{\gamma}_{3\ell}-H(\gamma_{10})+\gamma
_{30}\right\Vert \\
&  =O_{p}(\sqrt{n}[n^{-d_{1}(2\xi_{1}-1)/(2\xi_{1}+1)}\ln(n)]n^{-d_{1}}%
)=o_{p}^{{}}(1).
\end{align*}
Assumption 2 ii) then follows by the triangle and conditional Markov inequality.

Next, Assumption 3 i) follows by the form of $\phi_{1},$ $\phi_{2},$ and
$\phi_{3}$ given in Section 3 and
\begin{align*}
E[Y_{2t}-\gamma_{10}(X_{t})|X_{t}]  &  =0,\text{ }E[Y_{2t}\{H(\gamma
_{10}(X_{t+1}))-\gamma_{20}(X_{t})\}|X_{t}]=0,\\
E[Y_{1t}\{H(\gamma_{10}(X_{t+1}))-\gamma_{30}\}]  &  =0.
\end{align*}

We now proceed to verify that Assumption 3 iv) is satisfied. For ease of
exposition we suppress the $\ell$ subscript. Let $\hat{a}(x)=a(x,\theta
_{0},\hat{\gamma}_{2},\hat{\gamma}_{3})$ and $\hat{\pi}(x)=\pi(a(x,\theta
_{0},\hat{\gamma}_{2},\hat{\gamma}_{3})).$ Then%

\begin{align*}
\bar{\psi}(\hat{\gamma},\alpha_{0},\theta_{0})  &  =T+T_{1}+T_{2}+T_{3},\text{
}T=\int D(x_{t})\hat{\pi}(x_{t})\{y_{2t}-\Lambda(\hat{a}(x_{t}))\}F_{0}(dw),\\
T_{1}  &  =\int\alpha_{10}(x_{t})[y_{2t}-\hat{\gamma}_{1}(x_{t})]F_{0}%
(dw),\text{ }T_{3}=\alpha_{03}\int y_{1t}[H(\hat{\gamma}_{1\ell}%
(x_{t+1}))-\hat{\gamma}_{3}]F_{0}(dw),\\
T_{2}  &  =\int\alpha_{20}(x_{t},y_{2t})[H(\hat{\gamma}_{1}(x_{t+1}%
))-\hat{\gamma}_{2}(x_{t})]F_{0}(dw).
\end{align*}
Note that%
\begin{align*}
T  &  =\bar{T}_{1}+\bar{T}_{2}+R_{1}+R_{2},\text{ }\\
\bar{T}_{1}  &  =-\delta\int D(x_{t})\pi(x_{t})\Lambda_{a}(a(x_{t}))\left[
\hat{\gamma}_{2}(x_{t})-\gamma_{20}(x_{t})\right]  F_{0}(dw),\\
\bar{T}_{2}  &  =A\left(  \hat{\gamma}_{3}-\gamma_{30}\right)  ,\text{
}A=\delta\int D(x_{t})\pi(x_{t})\Lambda_{a}(a(x_{t}))F_{0}(dw),\\
R_{1}  &  =-\int D(x_{t})\hat{\pi}(x_{t})\Lambda_{aa}(\bar{a}(x_{t}%
))\left\vert \hat{a}(x_{t})-a(x_{t})\right\vert ^{2}F_{0}(dw),\\
R_{2}  &  =-\int D(x_{t})[\hat{\pi}(x_{t})-\pi(x_{t})]\Lambda_{a}%
(a(x_{t}))\left[  \hat{a}(x_{t})-a(x_{t})\right]  F_{0}(dw).
\end{align*}
Also,
\begin{align*}
\left\Vert R_{1}\right\Vert  &  \leq C\left(  \left\Vert \hat{\gamma}_{2\ell
}-\gamma_{20}\right\Vert ^{2}+\left\vert \hat{\gamma}_{3\ell}-\gamma
_{30}\right\vert ^{2}\right)  =O_{p}(n^{-4d_{1}\xi_{1}/(2\xi_{1}+1)}%
)=o_{p}(n^{-1/2}),\\
\left\Vert R_{2}\right\Vert  &  \leq C\left(  \left\Vert \hat{\gamma}_{2\ell
}-\gamma_{20}\right\Vert ^{2}+\left\vert \hat{\gamma}_{3\ell}-\gamma
_{30}\right\vert ^{2}\right)  =o_{p}(n^{-1/2}),
\end{align*}
so that
\[
T=\bar{T}_{1}+\bar{T}_{2}+o_{p}(n^{-1/2}).
\]
Next, note that by the definition of $\alpha_{20}(x,y_{2}),$%
\[
\bar{T}_{1}=\int\alpha_{20}(x_{t},y_{2t})\{\hat{\gamma}_{2}(x_{t})-\gamma
_{20}(x_{t})\}F_{0}(dw).
\]
Therefore%
\[
\bar{T}_{1}+T_{2}=\tilde{T}_{2},\text{ }\tilde{T}_{2}=\int\alpha_{10}%
(x_{t},y_{2t})[H(\hat{\gamma}_{1}(x_{t+1}))-\gamma_{20}(x_{t})]F_{0}(dw).
\]
Note that by $\gamma_{20}(x)=E[H(\gamma_{10}(X_{t+1}))|X_{t}=x,Y_{2t}=1]$ and
subtracting and adding the expression $\int\alpha_{20}(x_{t},y_{2t}%
)H(\gamma_{10}(x_{t+1}))F_{0}(dw)$ we obtain%
\begin{align*}
\tilde{T}_{2}  &  =\int\alpha_{20}(x_{t},y_{2t})[H(\hat{\gamma}_{1}%
(x_{t+1}))-H(\gamma_{10}(x_{t+1}))]F_{0}(dw)+\int\alpha_{20}(x_{t}%
,y_{2t})[H(\gamma_{10}(x_{t+1}))-\gamma_{20}(x_{t})]F_{0}(dw)\\
&  =\int\alpha_{20}(x_{t},y_{2t})[H(\hat{\gamma}_{1}(x_{t+1}))-H(\gamma
_{10}(x_{t+1}))]F_{0}(dw).
\end{align*}
Expanding then gives%
\begin{align*}
\tilde{T}_{2}  &  =\breve{T}_{2}+R_{3},\text{ }\breve{T}_{2}=\int\alpha
_{20}(x_{t},y_{2t})H_{p}(\gamma_{10}(x_{t+1}))[\hat{\gamma}_{1}(x_{t+1}%
)-\gamma_{10}(x_{t+1})]F_{0}(dw),\\
R_{3}  &  =\int\alpha_{20}(x_{t},y_{2t})H_{pp}(\bar{\gamma}_{1}(x_{t+1}%
))[\hat{\gamma}_{1}(x_{t+1})-\gamma_{10}(x_{t+1})]^{2}F_{0}(dw),
\end{align*}
where $\bar{\gamma}_{1}(x)$ is between $\hat{\gamma}_{1}(x_{t})$ and
$\gamma_{10}(x_{t}).$ It follows similarly to previous arguments that
$\left\Vert R_{3}\right\Vert \leq C\left\Vert \hat{\gamma}_{1}-\gamma
_{10}\right\Vert ^{2}=O_{p}(n^{-2d_{1}})=o_{p}(n^{-1/2}),$ so that%
\[
\tilde{T}_{2}=\breve{T}_{2}+o_{p}(n^{-1/2}).
\]
Also,
\[
\breve{T}_{2}=\int\zeta_{10}(x_{t})H_{p}(\gamma_{10}(x_{t}))[\hat{\gamma}%
_{1}(x_{t})-\gamma_{10}(x_{t})]F_{0}(dw),\text{ }\zeta_{10}(x)=E[\alpha
_{20}(X_{t},Y_{2t})|X_{t+1}=x].
\]
Note that from eq. (\ref{alpha ddc}),
\[
\alpha_{10}(x)=[\zeta_{10}(x)+\alpha_{30}\zeta_{20}(x)]H_{p}(\gamma
_{10}(x)),\,\zeta_{20}(x)=E[Y_{1t}|X_{t+1}=x].
\]
Then by $\int\zeta_{10}(x_{t})H_{p}(\gamma_{10}(x_{t}))[y_{2t}-\gamma
_{10}(x_{t})]F_{0}(dw)=0$ we have
\begin{align*}
\breve{T}_{2}+T_{1}  &  =\alpha_{30}\int\zeta_{20}(x_{t})H_{p}(\gamma
_{10}(x_{t}))[y_{2t}-\hat{\gamma}_{1}(x_{t})]F_{0}(dw)\\
&  =\alpha_{30}\int\zeta_{20}(x_{t})H_{p}(\gamma_{10}(x_{t}))[\gamma
_{10}(x_{t})-\hat{\gamma}_{1}(x_{t})]F_{0}(dw).
\end{align*}

Next, note that by iterated expectations and $\alpha_{30}=A/P_{1}$,
\[
T_{3}=\alpha_{03}\int y_{1t}[H(\hat{\gamma}_{1\ell}(x_{t+1}))-\hat{\gamma}%
_{3}]F_{0}(dw)=\alpha_{30}\int\zeta_{20}(x_{t})H(\hat{\gamma}_{1}(x_{t}%
))F_{0}(dw)-A\hat{\gamma}_{3}.
\]
Note also that
\[
A\gamma_{3}=AE[y_{1t}H(\gamma_{10}(x_{t+1})]/P_{1}=\alpha_{30}\int\zeta
_{20}(x_{t})H(\gamma_{10}(x_{t}))F_{0}(dw)
\]
Then by an expansion
\begin{align*}
\bar{T}_{2}+T_{3}  &  =\alpha_{03}\int\zeta_{20}(x_{t})H(\hat{\gamma}%
_{1}(x_{t}))F_{0}(dw)-A\gamma_{3}=\alpha_{03}\int\zeta_{20}(x_{t}%
)[H(\hat{\gamma}_{1}(x_{t}))-H(\gamma_{10}(x_{t}))]F_{0}(dw)\\
&  =\alpha_{03}\int\zeta_{20}(x_{t})[H(\hat{\gamma}_{1}(x_{t}))-H(\gamma
_{10}(x_{t}))]F_{0}(dw)\\
&  =\alpha_{30}\int\zeta_{20}(x_{t})H_{p}(\gamma_{10}(x_{t}))[\hat{\gamma}%
_{1}(x_{t})-\gamma_{10}(x_{t})]F_{0}(dw)+R_{4}=-(\breve{T}_{2}+T_{1})+R_{4},\\
R_{4}^{{}}  &  =\alpha_{30}\int\zeta_{20}(x_{t})H_{pp}(\bar{\gamma}_{1}%
(x_{t}))[\hat{\gamma}_{1}(x_{t})-\gamma_{10}(x_{t})]^{2}F_{0}(dw),
\end{align*}
where $\bar{\gamma}_{1}(x)$ is between $\hat{\gamma}_{1}(x_{t})$ and
$\gamma_{10}(x_{t}).$ It follows similarly to previous arguments that
$\left\Vert R_{3}\right\Vert \leq C\left\Vert \hat{\gamma}_{1}-\gamma
_{10}\right\Vert ^{2}=O_{p}(n^{-2d_{1}})=o_{p}(n^{-1/2}).$ Therefore%
\[
\bar{T}_{2}+T_{3}=-(\breve{T}_{2}+T_{1})+o_{p}(n^{-1/2}).
\]
Summarizing, it follows from what has been shown that%
\begin{align*}
\bar{\psi}(\hat{\gamma},\alpha_{0},\theta_{0})  &  =T+T_{1}+T_{2}+T_{3}%
=\bar{T}_{1}+\bar{T}_{2}+T_{1}+T_{2}+T_{3}+o_{p}(n^{-1/2})\\
&  =\tilde{T}_{2}+\bar{T}_{2}+T_{1}+T_{3}+o_{p}(n^{-1/2})=\breve{T}_{2}%
+\bar{T}_{2}+T_{1}+T_{3}+o_{p}(n^{-1/2})\\
&  =(\breve{T}_{2}+T_{1})+(\bar{T}_{2}+T_{3})+o_{p}(n^{-1/2})=(\breve{T}%
_{2}+T_{1})-(\breve{T}_{2}+T_{1})+o_{p}(n^{-1/2})\\
&  =o_{p}(n^{-1/2}),
\end{align*}
giving Assumption 3 iv).

Next, note that by the fixed trimming $H\left(  \hat{\gamma}_{1}(x)\right)  $
and $\hat{\gamma}_{1}(x)$ are uniformly bounded. Also, by Theorem 23,
$\hat{\gamma}_{2}(x)$ and $\hat{\gamma}_{3}$ are uniformly bounded with
probability approaching one, so%
\[
\left\Vert \hat{\Delta}_{\ell}(w)\right\Vert \leq C(\left\Vert \hat{\alpha
}_{1\ell}(x)-\alpha_{10}(x)\right\Vert +\left\Vert \hat{\alpha}_{2\ell
}(x)-\alpha_{20}(x)\right\Vert +\left\Vert \hat{\alpha}_{3\ell}-\alpha
_{30}\right\Vert ).
\]
The second condition of Lemma 16 then follows by Lemma A2. The first condition
of Lemma 16 also follows in a straightforward manner from uniform boundedness
of $\hat{\gamma}_{2}(x)$ and $\hat{\gamma}_{3}$ with probability approaching one.

Finally, Assumption 4 follows in a straightforward manner from the same
boundedness properties, so the conclusion follows by Theorem 18. $Q.E.D.$

\bigskip

\section{Appendix B: Convergence Rate for $\hat{\alpha}_{\ell}$ in Example 2.}

\bigskip

In this Appendix we state the conditions from Chernozhukov, Newey, and Singh
(2018) that are determine the convergence rates of various Lasso learners
$\hat{\alpha}$ in the paper$.$ We also prove Lemma 21, and refer to that proof
for results for the various $\hat{\alpha}$ in Example 3.

\textsc{Assumption B1:} \textit{There is }$C>0$\textit{ such that i) with
probability one }$\max_{1\leq j\leq p}|b_{j}(X)|\leq C$; ii) for every $n$
there is \textit{a }$p\times1$\textit{ vector }$\rho_{n}$\textit{ such that
}$\left\vert \rho_{n}\right\vert _{1}\leq C$\textit{ and} $\Vert\alpha
_{0}-b^{\prime}\rho_{n}\Vert^{2}=O(\varepsilon_{n}).$

\bigskip

This Assumption will suffice for a convergence rate for $\hat{\alpha}_{\ell}$
of $\sqrt{r_{\lambda}}.$ We can speed up this convergence rate under a
stronger approximate sparsity condition and a sparse eigenvalue condition. For
any $\rho=(\rho_{1},...,\rho_{p})$ let $\mathcal{J}=\{1,...,p\},$
$\mathcal{J}_{\rho}$ be the subset of $\mathcal{J}$ with $\rho_{j}\neq0$, and
$\mathcal{J}_{\rho}^{c}$ be the complement of $\mathcal{J}_{\rho}$ in
$\mathcal{J}$.

\bigskip

\textsc{Assumption B2:} \textit{i) there exists }$C,$\textit{ }$\xi>0$\textit{
such that for all }$\bar{s}$\textit{ with }$\bar{s}\leq C(\varepsilon_{n}%
^{2})^{-1/(1+2\xi)}$\textit{ there is }$\bar{\rho}$\textit{ with }$\bar{s}%
$\textit{ nonzero elements such that}%
\[
\left\Vert \alpha_{0}-b^{\prime}\bar{\rho}\right\Vert \leq C(\bar{s})^{-\xi}%
\]
\textit{ii) }$\bar{G}$\textit{ is nonsingular and has largest eigenvalue
uniformly bounded in }$n$\textit{; iii) for }$\rho=\bar{\rho}$\textit{ and
}$\rho=\arg\min_{\rho}\{\left\Vert \alpha_{0}-b^{\prime}\bar{\rho}\right\Vert
^{2}+2r_{\lambda}\sum_{j=1}^{p}\left\vert \rho_{j}\right\vert \}$\textit{
there is }$k>3$\textit{ such such that }
\[
\inf_{\{\delta:\delta\neq0,\sum_{j\in\mathcal{J}_{\rho_{L}}^{c}}|\delta
_{j}|\leq k\sum_{j\in\mathcal{J}_{\rho_{L}}}|\delta_{j}|\}}\frac
{\delta^{\prime}G\delta}{\sum_{j\in\mathcal{J}_{\rho_{L}}}\delta_{j}^{2}}>0.
\]

Before proving Lemma 21 we derive a convergence rate for $\hat{Q}_{\ell}.$ Let%
\[
\varepsilon_{n}=\sqrt{\ln(p)/(hn)}+h^{2}+n_{{}}^{-d_{\gamma}}.
\]
Also let $M_{j}=E[m(W,b_{j})]=E[\bar{\alpha}(X)b_{j}(X)]=E[\alpha_{0}%
(X)b_{j}(X)f(0|X)],$ $M=(M_{1},...,M_{p})^{\prime},$ and $\hat{M}_{\ell}%
=(\hat{M}_{\ell1},...,\hat{M}_{\ell p})^{\prime}.$

\bigskip

\textsc{Lemma B1:} \textit{If there is }$C$\textit{ such that }$\left\vert
b_{j}(X)\right\vert \leq C$\textit{ for all }$j$\textit{ then for
}$Q=E[f(0|X)b(X)b(X)^{\prime}],$%
\[
\left\vert \hat{Q}_{\ell}-Q\right\vert _{\infty}=O_{p}(\varepsilon_{n}),\text{
}\left\vert \hat{M}_{\ell}-M\right\vert _{\infty}=O_{p}(\varepsilon_{n}).
\]

\bigskip

\textbf{Proof of Lemma B1: }Consider%
\[
\hat{Q}_{\ell\ell^{\prime}}:=\frac{1}{n_{\ell^{\prime}}}\sum_{i\in
I_{\ell^{\prime}}}\frac{1}{h}K(\frac{Y_{i}-\hat{\gamma}_{\ell,\ell^{\prime}%
}(X_{i})}{h})b(X_{i})b(X_{i})^{\prime}.
\]
For notational convenience we drop the $\ell$ and $\ell^{\prime}$ subscripts
and replace $\sum_{i\in I_{\ell^{\prime}}}$ with $\sum_{i=1}^{n}$ while
retaining independence of $\hat{\gamma}$ from the data being averaged over.
Let $K_{h}(u)=h^{-1}K(u/h),$ $f(u|x)$ denote the conditional pdf of
$U=Y-\gamma_{0}(X)$ given $X=x,$ and $\hat{\Delta}(X)=\gamma_{0}%
(X)-\hat{\gamma}\left(  X\right)  .$ Note by two change of variables
$v=(u+\hat{\Delta})/h,$ and $\tilde{v}=u/h$%
\begin{align*}
\left\vert E[K_{h}(Y-\hat{\gamma}(X))|X]-E[K_{h}(U)|X]\right\vert  &
=\left\vert \int K_{h}(u-\hat{\Delta}(X))f(u|X)du-\int K_{h}%
(u)f(u|X)du\right\vert \\
&  =\left\vert \int K(v)[f(hv+\hat{\Delta}(X)|X)-f(hv|X)]dv\right\vert \\
&  \leq\int\left\vert K(v)\right\vert \left\vert f(hv+\hat{\Delta
}(X)|X)-f(hv|X)\right\vert dv\leq C\hat{\Delta}(X).
\end{align*}
Also, note that by a mean value expansion for small enough $h$,
\begin{align*}
f(hv|X)  &  =\sum_{k=0}^{1}v^{k}h^{k}\frac{d^{k}f(0|X)}{du^{k}}+h^{2}R(v,X),\\
\left\vert E[K_{h}(U)|X]-f(0|X)\right\vert  &  =\left\vert \int
K(v)[f(hv|X)-f(0|X)]dv\right\vert \leq Ch^{2}%
\end{align*}
Note also that conditional on $\hat{\gamma}$, by the conditional Markov
inequality and $\int\left\vert \hat{\Delta}(x)\right\vert F_{0}(x)\leq
\left\Vert \hat{\gamma}-\gamma_{0}\right\Vert =O_{p}(n^{-d_{\gamma}})$ we have%
\[
\frac{1}{n}\sum_{i=1}^{n}\left\vert \hat{\Delta}(X_{i})\right\vert
=O_{p}(n^{-d_{\gamma}}).
\]
Therefore by $\left\vert b(X_{i})\right\vert _{\infty}\leq C$ we have%
\begin{align*}
&  \left\vert \frac{1}{n}\sum_{i=1}^{n}b(X_{i})b(X_{i})^{\prime}%
\{E[K_{h}(Y_{i}-\hat{\gamma}(X_{i}))|X_{i}]-f(0|X_{i})\}\right\vert _{\infty
}\\
&  \leq C\frac{1}{n}\sum_{i=1}^{n}\left\vert \hat{\Delta}(X_{i})\right\vert
+Ch^{2}=O_{p}(n^{-d_{\gamma}}+h^{2}),
\end{align*}
Note also that by a change of variable $v=(U+\hat{\Delta}(X))/h^{{}}$%
\begin{align*}
\left\vert b_{j}(X)b_{j^{\prime}}(X)K_{h}(Y-\hat{\gamma}(X_{i}))\right\vert
&  \leq Ch^{-1},\text{ }\\
E[b_{j}(X)^{2}b_{j^{\prime}}(X)^{2}K_{h}(Y-\hat{\gamma}(X))^{2}]  &  \leq
CE[K_{h}(U+\hat{\Delta}(X))^{2}]\leq Ch^{-1}E[\int K(v)^{2}f(hv-\hat{\Delta
}(X)|X)]\leq Ch^{-1}.
\end{align*}
It then follows by Lemma 19.32 of Van der Vaart (1998) and a standard argument
that%
\[
\left\vert \hat{Q}-\frac{1}{n}\sum_{i=1}^{n}b(X_{i})b(X_{i})^{\prime}%
E[K_{h}(Y_{i}-\hat{\gamma}(X_{i}))|X_{i}]\right\vert _{\infty}=O_{p}%
(\sqrt{\frac{\ln(p)}{hn}}).
\]
It then follows by the triangle inequality that%
\[
\left\vert \hat{Q}-\frac{1}{n}\sum_{i=1}^{n}b(X_{i})b(X_{i})^{\prime}%
f(0|X_{i})\right\vert _{\infty}=O_{p}(\varepsilon_{n}).
\]
In addition it follows by a standard application of Hoeffding's inequality
that%
\[
\left\vert \frac{1}{n}\sum_{i=1}^{n}b(X_{i})b(X_{i})^{\prime}f(0|X_{i}%
)-Q\right\vert _{\infty}=O_{p}(\sqrt{\frac{\ln(p)}{n}})=O_{p}(\varepsilon
_{n}),
\]
so $\left\vert \hat{Q}-Q\right\vert _{\infty}=O_{p}(\varepsilon_{n})$ follows
by the triangle inequality. The first conclusion follows by another
application of the triangle inequality. The second conclusion follows by
another application of Hoeffding's inequality since $\sqrt{\ln(p)/n}%
\leq\varepsilon_{n}$ for $n$ large enough. $Q.E.D.$

\bigskip

\textbf{Proof of Lemma 21: }Using Lemma B1 we apply the results of
Chernozhukov, Newey, and Singh (2018) for the distribution of $X$ with
expectation $\bar{E}$ given by $\bar{E}[a(X)]=E[a(X)f(0|X)]/E[f(0|X)]$. Since
$f(0|X)$ is bounded and bounded away from zero, a convergence rate for this
expectation will imply a convergence rate in the original expectation. The
first and second conclusions of Lemma 21 then follow by Theorems 1 and 3
respectively in Chernozhukov, Newey, Singh (2018). $Q.E.D.$

\bigskip

Acknowledgements

Escanciano acknowledges research support by Spanish grant PGC
2018-096732-B-100 and Newey by NSF Grant 1757140. Helpful comments were
provided by M. Cattaneo, X. Chen, B. Deaner, J. Hahn, M. Jansson, Z. Liao, O.
Linton, R. Moon, A. Pakes, A. de Paula, P. Phillips, V. Semenova, Y. Zhao and
participants in seminars at Cambridge, Columbia, Cornell, Harvard-MIT, UCL,
USC, Yale, and Xiamen. B. Deaner provided capable research assistance.

\setlength{\parindent}{-.5cm} \setlength{\parskip}{.1cm}

\begin{center}
\textbf{REFERENCES}
\end{center}

\textsc{Ackerberg, D., X. Chen, and J. Hahn} (2012): "A Practical Asymptotic
Variance Estimator for Two-step Semiparametric Estimators," \textit{The Review
of Economics and Statistics} 94: 481--498.

\textsc{Ackerberg, D., X. Chen, J. Hahn, and Z. Liao} (2014): "Asymptotic
Efficiency of Semiparametric Two-Step GMM," \textit{The Review of Economic
Studies} 81: 919--943.

\textsc{Ai, C. {\small and} X. Chen} (2003): \textquotedblleft Efficient
Estimation of Models with Conditional Moment Restrictions Containing Unknown
Functions,\textquotedblright\ \textit{Econometrica }71, 1795-1843.

\textsc{Ai, C. {\small and} X. Chen} (2007): "Estimation of Possibly
Misspecified Semiparametric Conditional Moment Restriction Models with
Different Conditioning Variables," \textit{Journal of Econometrics} 141, 5--43.

\textsc{Ai, C. {\small and} X. Chen} (2012): "The Semiparametric Efficiency
Bound for Models of Sequential Moment Restrictions Containing Unknown
Functions," \textit{Journal of Econometrics} 170, 442--457.

\textsc{Andrews, D.W.K.} (1994): \textquotedblleft Asymptotics for
Semiparametric Models via Stochastic Equicontinuity,\textquotedblright%
\ \textit{Econometrica} 62, 43-72.

\textsc{Angrist, J.D. and A.B. Krueger} (1995): "Split-Sample Instrumental
Variables Estimates of the Return to Schooling," \textit{Journal of Business
and Economic Statistics} 13, 225-235.

\textsc{Athey, S., G. Imbens, and S. Wager} (2018): "Approximate residual
balancing: debiased inference of average treatment effects in high
dimensions," \textit{Journal of the Royal Statistical Society, Series B,} 80, 597-623.

\textsc{Avagyan, V. and S. Vansteelandt} (2017): "Honest Data-adaptive
Inference for the Average Treatment Effect Using Penalised Bias-Reduced
Double-Robust Estimation," 

arxiv.org/pdf/1708.03787.pdf.

\textsc{Bajari, P., V. Chernozhukov, H. Hong, and D. Nekipelov} (2009):
"Nonparametric and Semiparametric Analysis of a Dynamic Discrete Game,"
working paper, Stanford.{}

\textsc{Bajari, P., H. Hong, J. Krainer, and D. Nekipelov} (2010): "Estimating
Static Models of Strategic Interactions," \textit{Journal of Business and
Economic Statistics} 28, 469-482.

\textsc{Belloni, A., and V. Chernozhukov} (2011): "$\ell1$-Penalized Quantile
Regression in High-Dimensional Sparse Models," \textit{Annals of Statistics}
9, 82-130.

\textsc{Belloni, A., D. Chen, V. Chernozhukov, and C. Hansen} (2012):
\textquotedblleft Sparse Models and Methods for Optimal Instruments with an
Application to Eminent Domain,\textquotedblright\ \textit{Econometrica} 80, 2369--2429.

\textsc{Belloni, A., V. Chernozhukov, and Y. Wei }(2013): \textquotedblleft
Honest Confidence Regions for Logistic Regression with a Large Number of
Controls,\textquotedblright\ arXiv preprint arXiv:1304.3969.

\textsc{Belloni, A., V. Chernozhukov, and C. Hansen} (2014): "Inference on
Treatment Effects after Selection among High-Dimensional Controls,"
\textit{The Review of Economic Studies} 81, 608--650.

\textsc{Belloni, A., V. Chernozhukov, and K. Kato} (2015): "Uniform
Post-Selection Inference for Least Absolute Deviation Regression and Other
Z-Estimation Problems," \textit{Biometrika }102, 77--94.

\textsc{Belloni, A., V. Chernozhukov, I. Fernandez-Val, and C. Hansen} (2017):
"Program Evaluation and Causal Inference with High-Dimensional Data,"
\textit{Econometrica }85, 233-298.

\textsc{Bera, A.K., G. Montes-Rojas, and W. Sosa-Escudero} (2010): "General
Specification Testing with Locally Misspecified Models," \textit{Econometric
Theory} 26, 1838--1845.

\textsc{Bickel, P.J. and Y. Ritov} (1988): "Estimating Integrated Squared
Density Derivatives: Sharp Best Order of Convergence Estimates,"
\textit{Sankhy\={a}: The Indian Journal of Statistics, Series A }238, 381-393. \ 

\textsc{Bickel, P.J., C.A.J. Klaassen, Y. Ritov, {\small and} J.A. Wellner}
(1993): \textit{Efficient and Adaptive Estimation for Semiparametric Models},
Springer-Verlag, New York.

\textsc{Bonhomme, S., and M. Weidner} (2018): "Minimizing Sensitivity to
Misspecification," arxiv.org/abs/1807.02161v1.

\textsc{Bravo, F., J.C. Escanciano, and I. van Keilegom} (2020): "Two-step
Semiparametric Likelihood Inference," \textit{Annals of Statistics }48, 1-26.

\textsc{Carone, M., A.R. Luedtke, and M.J. van der Laan} (2016): "Toward
Computerized Efficient Estimation in Infinite Dimensional Models," 

arXiv:1608.08717v1.

\textsc{Cattaneo, M.D., and M. Jansson }(2018): "Kernel-Based Semiparametric
Estimators: Small Bandwidth Asymptotics and Bootstrap Consistency,"
\textit{Econometrica }86, 955--995.

\textsc{Cattaneo, M.D., M. Jansson, and X. Ma }(2018): "Two-step Estimation
and Inference with Possibly Many Included Covariates," \textit{The Review of
Economic Studies} 86, 1095--1122.

\textsc{Chen, X. and X. Shen }(1997): \textquotedblleft Sieve Extremum
Estimates for Weakly Dependent Data,\textquotedblright\ \textit{Econometrica}
66, 289-314.

\textsc{Chen, X., O.B. Linton, {\small and} I. {van Keilegom}} (2003):
\textquotedblleft Estimation of Semiparametric Models when the Criterion
Function Is Not Smooth,\textquotedblright\ \textit{Econometrica} 71, 1591-1608.

\textsc{Chernozhukov, V., C. Hansen, and M. Spindler} (2015): "Valid
Post-Selection and Post-Regularization Inference: An Elementary, General
Approach," \textit{Annual Review of Economics} 7: 649--688.

\textsc{Chernozhukov, V., J.C. Escanciano, H. Ichimura, W.K. Newey} (2016):
"Locally Robust Semiparametric Estimation," 

arxiv.org/pdf/1608.00033v1.pdf.

\textsc{Chernozhukov, V., D. Chetverikov, M. Demirer, E. Duflo, C. Hansen, W.
Newey, J. Robins }(2018): "Debiased/Double Machine Learning for Treatment and
Structural Parameters,\textit{Econometrics Journal }21, C1-C68.

\textsc{Chernozhukov, V., J.A. Hausman, and W.K. Newey }(2018): "Demand
Analysis with Many Prices," working paper, MIT.

\textsc{Chernozhukov, V., W.K. Newey, and J. Robins }(2018): "Double/De-Biased
Machine Learning Using Regularized Riesz Representers," arxiv.org/abs/1802.08667v1.

\textsc{Chernozhukov, V., W.K. Newey, and V. Semenova }(2019): "Inference on
Average Welfare with High Dimensional State Space," 

https://arxiv.org/pdf/1908.09173.pdf.

\textsc{Escanciano, J-C., D. Jacho-Chavez, and A. Lewbel }(2016):
\textquotedblleft Identification and Estimation of Semiparametric Two Step
Models\textquotedblright, \textit{Quantitative Economics} 7, 561-589.

\textsc{Farrell, M. }(2015): "Robust Inference on Average Treatment Effects
with Possibly More Covariates than Observations," \textit{Journal of
Econometrics} 189, 1--23.

\textsc{Firpo, S. and C. Rothe} (2019): \textquotedblleft Properties of Doubly
Robust Estimators when Nuisance Functions are Estimated
Nonparametrically,\textquotedblright\ \textit{Econometric Theory} 35, 1048--1087.

\textsc{Foster, D.F. and V. Syrgkanis (2019):} "Orthogonal Statistical
Learning," 

arxiv.org/pdf/1901.09036.pdf.

\textsc{Graham, B.W.} (2011): "Efficiency Bounds for Missing Data Models with
Semiparametric Restrictions," \textit{Econometrica }79, 437--452.

\textsc{Hahn, J. (1998):} "On the Role of the Propensity Score in Efficient
Semiparametric Estimation of Average Treatment Effects," \textit{Econometrica}
66, 315-331.

\textsc{Hahn, J. and G. Ridder} (2013): "Asymptotic Variance of Semiparametric
Estimators With Generated Regressors," \textit{Econometrica} 81, 315-340.

\textsc{Hahn, J. and G. Ridder} (2019): \textquotedblleft Three-stage
Semi-Parametric Inference: Control Variables and Differentiability,"
\textit{Journal of Econometrics} 211, 262-293.

\textsc{Hampel, F.R. (1974):} "The Influence Curve and Its Role in Robust
Estimation," \textit{Journal of the American Statistical Association }69, 383-393.

\textsc{Hasminskii, R.Z. and I.A. Ibragimov} (1978): "On the Nonparametric
Estimation of Functionals," \textit{Proceedings of the 2nd Prague Symposium on
Asymptotic Statistics}, 41-51.

\textsc{Hirano, K., G. Imbens, and G. Ridder} (2003): "Efficient Estimation of
Average Treatment Effects Using the Estimated Propensity Score,"
\textit{Econometrica} 71: 1161--1189.

\textsc{Hirshberg, D.A. and S. Wager} (2019): "Augmented Minimax Linear
Estimation," arxiv.org/pdf/1712.00038.pdf.

\textsc{Hotz, V.J. and R.A. Miller} (1993): "Conditional Choice Probabilities
and the Estimation of Dynamic Models," \textit{Review of Economic Studies} 60, 497-529.

\textsc{Huber, P. (1981):} \textit{Robust Statistics, }New York:\ Wiley.

\textsc{Ichimura, H. }(1993): "Estimation of Single Index Models,"
\textit{Journal of Econometrics} 58, 71-120.

\textsc{Ichimura, H., {\small and} S. Lee} (2010): \textquotedblleft
Characterization of the Asymptotic Distribution of Semiparametric
M-Estimators,\textquotedblright\ \textit{Journal of Econometrics} 159, 252--266.

\textsc{Ichimura, H. and W.K. Newey} (2017): "The Influence Function of
Semiparametric Estimators," CEMMAP Working Paper, CWP06/17.

\textsc{Kandasamy, K., A. Krishnamurthy, B. P%
\'{}%
oczos, L. Wasserman, J.M. Robins }(2015): "Influence Functions for Machine
Learning: Nonparametric Estimators for Entropies, Divergences and Mutual
Informations," arxiv.

\textsc{Klaassen, C.A.J.} (1987): "Consistent Estimation of the Influence
Function of Locally Asymptotically Linear Estimators," \textit{Annals of
Statistics} 15, 1548-1562.

\textsc{Lee, Lung-fei }(2005): \textquotedblleft A $C(\alpha)$-type Gradient
Test in the GMM Approach,\textquotedblright\ working paper.

\textsc{Leeb, H. and B.M. Potscher} (2005): "Model Selection and Inference:
Facts and Fiction," \textit{Econometric Theory} 21, 21-59.

\textsc{Leeb, H. and B.M. Potscher} (2008): "Sparse Estimators and the Oracle
Property, or the Return of Hodges' Estimator," \textit{Journal of
Econometrics} 142, 201-211.

\textsc{Luedtke, A.R. and M.J. Van der Laan} (2016): "Statistical Inference
For the Mean Outcome Under A Possibly Non-unique Optimal Treatment Strategy,"
\textit{Annals of Statistics} 44 713.

\textsc{Luenberger, D.G.} (1969): \textit{Optimization by Vector Space
Methods}, New York: Wiley.

\textsc{Murphy, K.M. and R.H. Topel} (1985): "Estimation and Inference in
Two-Step Econometric Models," \textit{Journal of Business and Economic
Statistics }3, 370-379.

\textsc{Newey, W.K.} (1984): "A Method of Moments Interpretation of Sequential
Estimators," \textit{Economics Letters }14, 201-206.

\textsc{Newey, W.K.} (1990): "Semiparametric Efficiency Bounds,"
\textit{Journal of Applied Econometrics }5, 99-135.

\textsc{Newey, W.K.} (1991): \textquotedblright Uniform Convergence in
Probability and Stochastic Equicontinuity,\textquotedblright%
\ \textit{Econometrica} 59, 1161-1167.

\textsc{Newey, W.K.} (1994a): "The Asymptotic Variance of Semiparametric
Estimators," \textit{Econometrica} 62, 1349-1382.

\textsc{Newey, W.K.} (1994b): \textquotedblright Kernel Estimation of Partial
Means and a General Variance Estimator,\textquotedblright\ \textit{Econometric
Theory} 10, 233-253.

\textsc{Newey, W.K.} (1997): \textquotedblright Convergence Rates and
Asymptotic Normality for Series Estimators,\textquotedblright\ \textit{Journal
of Econometrics} 79, 147-168.

\textsc{Newey, W.K., {\small and} D. McFadden }(1994): \textquotedblleft Large
Sample Estimation and Hypothesis Testing," in \textit{Handbook of
Econometrics}, Vol. 4, ed. by R. Engle, and D. McFadden, pp. 2113-2241. North Holland.

\textsc{Newey, W.K., {\small and} J.L. Powell }(1989): "Instrumental Variable
Estimation of Nonparametric Models," presented at Econometric Society winter
meetings, 1988.

\textsc{Newey, W.K., {\small and} J.L. Powell }(2003): "Instrumental Variable
Estimation of Nonparametric Models," \textit{Econometrica} 71, 1565-1578.

\textsc{Newey, W.K., F. Hsieh, {\small and} J.M. Robins} (1998):
\textquotedblleft Undersmoothing and Bias Corrected Functional Estimation,"
MIT Dept. of Economics working paper\ 72, 947-962.

\textsc{Newey, W.K., F. Hsieh, {\small and} J.M. Robins} (2004):
\textquotedblleft Twicing Kernels and a Small Bias Property of Semiparametric
Estimators,\textquotedblright\ \textit{Econometrica} 72, 947-962.

\textsc{Newey, W.K., and J. Robins} (2017): "Cross Fitting and Fast Remainder
Rates for Semiparametric Estimation," CEMMAP Working paper WP41/17.

\textsc{Neyman, J.} (1959): \textquotedblleft Optimal Asymptotic Tests of
Composite Statistical Hypotheses,\textquotedblright\ \textit{Probability and
Statistics, the Harald Cramer Volume}, ed., U. Grenander, New York, Wiley.

\textsc{Pfanzagl, J., and W. Wefelmeyer} (1982): "Contributions to a General
Asymptotic Statistical Theory. Springer Lecture Notes in Statistics.

\textsc{Pakes, A. and G.S. Olley} (1995): "A Limit Theorem for a Smooth Class
of Semiparametric Estimators," \textit{Journal of Econometrics} 65, 295-332.

\textsc{Powell, J.L., J.H. Stock, and T.M. Stoker }(1989): "Semiparametric
Estimation of Index Coefficients," \textit{Econometrica} 57, 1403-1430.

\textsc{Robins, J.M., and Rotnitzky, A.} (1992): "Recovery of Information and
Adjustment for Dependent Censoring Using Surrogate Markers," \textit{AIDS
Epidemiology - Methodological Issues}, Eds: Jewell N., Dietz K., Farewell V.
Boston, MA: Birkh\"{a}user. pp. 297-331.

\textsc{Robins, J.M., A. Rotnitzky, and L.P. Zhao }(1994): "Estimation of
Regression Coefficients When Some Regressors Are Not Always Observed,"
\textit{Journal of the American Statistical Association} 89: 846--866.

\textsc{Robins, J.M. and A. Rotnitzky }(1995): "Semiparametric Efficiency in
Multivariate Regression Models with Missing Data," \textit{Journal of the
American Statistical Association} 90:122--129.

\textsc{Robins, J.M., A. Rotnitzky, and L.P. Zhao }(1995): "Analysis of
Semiparametric Regression Models for Repeated Outcomes in the Presence of
Missing Data," \textit{Journal of the American Statistical Association} 90,106--121.

\textsc{Robins, J.M.,and A. Rotnitzky (2001):} Comment on \textquotedblleft
Semiparametric Inference: Question and an Answer,\textquotedblright\ by P.A.
Bickel and J. Kwon, \textit{Statistica Sinica }11, 863-960.

\textsc{Robins, J.M., A. Rotnitzky, and M. van der Laan} \ (2000): "Comment on
'On Profile Likelihood'\ by S. A. Murphy and A. W. van der Vaart,
\textit{Journal of the American Statistical Association} 95, 431-435.

\textsc{Robins, J.M., L. Li, E. Tchetgen, and A. van der Vaart} (2008):
"Higher Order Influence Functions and Minimax Estimation of Nonlinear
Functionals," \textit{IMS Collections Probability and Statistics: Essays in
Honor of David A. Freedman, Vol 2, }335-421.

\textsc{Robins, J., P. Zhang, R. Ayyagari, R. Logan, E. Tchetgen, L. Li, A.
Lumley, and}

\textsc{A. van der Vaart} (2013): "New Statistical Approaches to
Semiparametric Regression with Application to Air Pollution Research,"
Research Report Health E Inst.

\textsc{Robinson, P.M. }(1988): "`Root-N-consistent Semiparametric
Regression," \textit{Econometrica} 56, 931-954.

\textsc{Rotnitzky, A., E. Smucler, and J.M. Robins} (2019): "Characterization
of Parameters with a Mixed Bias Property," arXiv, 

https://arxiv.org/pdf/1904.03725.pdf

\textsc{Rust, J. }(1987): "Optimal Replacement of GMC Bus Engines: An
Empirical Model of Harold Zurcher," \textit{Econometrica }55, 999-1033.

\textsc{Santos, A.} (2011): "Instrumental Variable Methods for Recovering
Continuous Linear Functionals,"\ \textit{Journal of Econometrics}, 161, 129-146.

\textsc{Scharfstein D.O., A. Rotnitzky, and J.M. Robins (1999): }Rejoinder to
\textquotedblleft Adjusting For Nonignorable Drop-out Using Semiparametric
Non-response Models,\textquotedblright\ \textit{Journal of the American
Statistical Association }94, 1135-1146.

\textsc{Schick, A.} (1986): "On Asymptotically Efficient Estimation in
Semiparametric Models," \textit{Annals of Statistics} 14, 1139-1151.

\textsc{Semenova, V.} (2018): "Machine Learning for Set-Identified Linear
Models," https://arxiv.org/pdf/1712.10024.pdf.

\textsc{Severini, T.A. and W.H. Wong} (1992): "Profile Likelihood and
Conditionally Parametric Models," \textit{Annals of Statistics} 20, 1768-1802.

\textsc{Severini, T. and G. Tripathi (2006): "}Some Identification Issues in
Nonparametric Linear Models with Endogenous Regressors," \textit{Econometric
Theory} 22, 258-278.

\textsc{Severini, T. and G. Tripathi (2012):} "Efficiency Bounds for
Estimating Linear Functionals of Nonparametric Regression Models with
Endogenous Regressors," \textit{Journal of Econometrics} 170, 491-498.

\textsc{Shen, X. }(1997): "On Methods of Sieves and Penalization,"
\textit{Annals of Statistics} 25, 2555-2591.

\textsc{Singh, R., and L. Sun} (2019): "De-biased Machine Learning in
Instrumental Variable Models for Treatment Effects," 

https://arxiv.org/pdf/1909.05244.pdf.

\textsc{Stoker, T. }(1986): "Consistent Estimation of Scaled Coefficients,"
\textit{Econometrica} 54, 1461-1482.

\textsc{Tan, Z. }(2018): Model-Assisted Inference for Treatment Effects Using
Regularized Calibrated Estimation with High-Dimensional Data," 

arxiv.org/pdf/1801.09817.pdf.

\textsc{Toth, B., and M van der Laan} (2016): "TMLE for Marginal Structural
Models Based On An Instrument," Technical Report, UC Berkeley Division of Biostatistics.

\textsc{van der Laan, M.J. and S. Rose }(2011): \textit{Targeted Learning:
Causal Inference for Observational and Experimental Data,} Springer Science \&
Business Media.

\textsc{van der Laan, M.J. and D. Rubin }(2006): "Targeted Maximum Likelihood
Learning," \textit{The International Journal of Biostatistics} 2.

\textsc{{van der Vaart}, A.W.} (1991): \textquotedblleft On Differentiable
Functionals,\textquotedblright\ \textit{The Annals of Statistics,} 19, 178-204.

\textsc{{van der Vaart}, A.W.} (1998): \textit{Asymptotic Statistics,
}Cambridge University Press, Cambridge, England.

\textsc{{van der Vaart}, A.W.} (2014): "Higher Order Tangent Spaces and
Influence Functions," \textit{Statistical Science} 29, 679--686.

\textsc{Vermeulen, K. and S. Vansteelandt} (2015): "Bias-Reduced Doubly Robust
Estimation," \textit{Journal of the American Statistical Association} 110, 1024-1036.

\textsc{Von Mises, R. }(1947), "On the Asymptotic Distribution of
Differentiable Statistical Functions," \textit{Annals of Mathematical
Statistics} 18, 309-348.

\textsc{Wooldridge, J.M. }(1991): \textquotedblleft On the Application of
Robust, Regression-Based Diagnostics to Models of Conditional Means and
Conditional Variances,\textquotedblright\ \textit{Journal of Econometrics} 47, 5-46.

\textsc{Zhao, P., and B. Yu} (2006): "On Model Selection Consistency of
Lasso," \textit{Journal of Machine Learning Research} 7, 2541-2563.

\bigskip

\setlength{\parindent}{.0cm} \setlength{\parskip}{.1cm}

\end{document}